\documentclass[10pt]{amsart}
\usepackage[T1]{fontenc}
\usepackage[utf8]{inputenc}
\usepackage{amsmath,amsfonts,amssymb}
\usepackage{graphicx}
\usepackage{amsthm}
\usepackage[english]{babel}
\usepackage{enumitem}
\allowdisplaybreaks
\usepackage{xcolor}
\usepackage{verbatim}
\usepackage{tikz}
\usetikzlibrary{arrows}
\usepackage[font=small,labelfont=bf]{caption}
\usepackage[a4paper, margin=2.75 cm]{geometry}
\usepackage[colorlinks=true, linkcolor=blue, citecolor=red, urlcolor=blue]{hyperref}

\usepackage{booktabs}
\newcommand{\R}{{\mathbb{R}}}
\newcommand{\C}{{\mathbb{C}}}
\newcommand{\Z}{{\mathbb{Z}}}
\newcommand{\N}{{\mathbb{N}}}

\newcommand{\Cc}{\mathcal{C}}
\newcommand{\Ec}{\mathcal{E}}
\newcommand{\Fc}{\mathcal{F}}
\newcommand{\Gc}{\mathcal{G}}
\newcommand{\Hc}{\mathcal{H}}
\newcommand{\Nc}{\mathcal{N}}
\newcommand{\Sc}{\mathcal{S}}
\newcommand{\Vc}{\mathcal{V}}

\numberwithin{equation}{section}

\newtheorem{theorem}{Theorem}[section]
\newtheorem{proposition}[theorem]{Proposition}
\newtheorem{lemma}[theorem]{Lemma}
\newtheorem{remark}[theorem]{Remark}
\newtheorem{example}[theorem]{Example}
\newtheorem{definition}[theorem]{Definition}

\newtheorem{corollary}[theorem]{Corollary}

\begin{document}
\title[Hodge Laplacians \& ESA]{Hodge Laplacians on Weighted Simplicial Complexes: Forms, Closures, and Bounded Realizations}


 \author[\tiny Marwa Ennaceur]{Marwa Ennaceur$^{1}$}
 \address{$^{1}$Department of Mathematics, College of Science, University of Ha'il, Ha'il 81451, Saudi Arabia}
\email{\tt mar.ennaceur@uoh.edu.sa}
\author[\tiny Amel Jadlaoui]{Amel Jadlaoui$^{2}$ }
 \address{$^{2}$Department of Mathematics, Faculty of Sciences of Sfax, University
of Sfax, 3000 Sfax, Tunisie}
\email{\tt amel.jadlaoui@yahoo.com}
\begin{abstract}
We establish operator-norm bounds for discrete Hodge Laplacians on weighted flag complexes of a fixed dimension $n$, whose $k$-simplices are the $(k+1)$-cliques of a weighted graph; essential self-adjointness on natural cores follows, with no completeness or curvature assumption. Dual up/down degrees give Schur-type bounds in every degree. At top degree a unitary conjugation turns $\Delta_{n}$ into an adjacency operator on the line-complex plus a diagonal potential, to which Schur's test applies directly. At $n=1$ this is the Anderson--Morley edge-degree bound $\|\Delta_{1}\|\le\max\{\deg u+\deg v\}$, obtained here on infinite and on weighted graphs, whence $\|\Delta_{1}\|\le2d$ in the $d$-regular case. The sharpness analysis is carried out at $n=1$, for the edge block; in higher degrees we prove boundedness and essential self-adjointness, not sharpness. We give a spectral criterion for the attainment of $2d$, together with sufficient conditions: it is attained on every finite $d$-regular bipartite graph, and on every infinite one that is amenable. Amenability cannot be dropped: the $d$-regular tree, $d\ge3$, has exact norm $d+2\sqrt{d-1}<2d$. The standard periodic lattices being amenable, we compute their exact norms from the Bloch symbols: $2d$ in the bipartite cases, against $9$ for the triangular and $16$ for the face-centered cubic lattice, on which it is strict. Finally, an ordered coloring  which any countable complex admits reformulates the skew model unitarily on unoriented simplices, the orientation signs becoming a function of the colors alone.
\end{abstract}
\maketitle
\tableofcontents 
\section{Introduction}\label{sec:intro}

Throughout this introduction we use the following notation, which will be made fully precise in Section~\ref{sec:prelim-comb} and Section~\ref{section3}.
The skew convention is used throughout this paper, and we therefore carry \emph{no} convention subscript in the notation: $d^k$, $\delta^k$, $\Delta_k$ always denote the skew objects. On a colorable complex, every statement transfers verbatim, through an explicit unitary, to the canonical-basis reformulation of Section~\ref{sec:color} (Corollary~\ref{cor:transfer}).
For each $k\in\{1,\ldots,n\}$, the (formal) coboundary on compactly supported $(k-1)$-cochains is denoted by
\[
d^{k}\colon \Cc^{k-1}_{c}(\Vc)\longrightarrow \Cc^{k}(\Vc),
\]
the image consisting of compactly supported cochains as soon as the complex is locally finite (Section~\ref{sec:prelim-comb}),
and its formal adjoint with respect to the weighted inner product is denoted by $\delta^{k}$ (this $\delta^k$ is sometimes written $(d^k)^{\ast}$ in the literature). In the skew convention the relation $d^{k+1}\,d^k=0$ holds unconditionally; see Remark~\ref{rem:dsq-zero}.
Throughout this paper, we work with flag complexes (also called clique complexes), i.e., simplicial complexes whose $k$-simplices are precisely the $(k+1)$-cliques of a weighted underlying graph; see Definitions~\ref{def:Fk-clique}--\ref{def:simplicial-complex}. This is a substantive restriction: it excludes, e.g., ``hollow'' configurations in which the edges of a triangle are present without the corresponding $2$-simplex. The hypothesis is used essentially in two places: the line-complex reduction (Theorem~\ref{thm:kn-esa}) and the Schur-type bounds of Section~\ref{sec:normalization}. In both, the cofaces of a face are identified with the vertices of the graph-neighborhood intersection $F_{(x_1,\dots,x_k)}$. That identification may fail on a general abstract simplicial complex, so the sharp bounds obtained here do not automatically transfer beyond the flag setting.

The bounds established here are operator-norm bounds for the Hodge Laplacians on weighted simplicial complexes, from which essential self-adjointness (ESA) on a natural compactly supported core follows as a corollary. ESA is obtained from \emph{boundedness} on $\ell^2$, and not through geometric completeness, $\chi$-completeness, curvature, or capacity. This is a weaker goal than ESA in the unbounded regime, but the constants we get are geometry-free.

We work with the (Friedrichs) energy block
\[
\Delta_{k}\;=\;d^{k}\,\delta^{k}+\delta^{k+1}\,d^{k+1}\quad\text{on }\ell^2(\text{$k$-faces},m_k),
\]
and its \emph{normalized} counterpart on the unweighted $\ell^2$,
\[
\widetilde\Delta_{k}:=M_k^{1/2}\,\Delta_{k}\,M_k^{-1/2},
\]
where $M_k$ is multiplication by $m_k$. Since $M_k^{-1/2}$ is a \emph{unitary} from the unweighted $\ell^2$ onto $\ell^2(m_k)$, this is a unitary transport and no ellipticity of the weights is needed for it: $\widetilde\Delta_k$ and the Friedrichs realization $\Delta_k^{\mathrm F}$ are isospectral, and of equal norm, in complete generality (Lemma~\ref{lem:elliptic-similarity}).
On an unweighted complex $M_k=I$ and $\widetilde\Delta_{k}=\Delta_{k}$, so the two notations are interchangeable in every unweighted statement below --- and in the abstract, where the tilde is omitted for readability. Most of our bounds are stated for $\widetilde\Delta_{k}$, which has a uniform spectral interpretation across weights. The top-degree results of Sections~\ref{sec:esa} and~\ref{sharp} are stated for $\Delta_{k}$ itself; by the preceding remark the same numerical bounds hold for $\widetilde\Delta_{k}$.

The main contributions are as follows.
\begin{itemize}[leftmargin=2em]
  \item \emph{Schur-type bounds in every degree, and the top-degree reduction.}
  For every degree $k$, Schur-type criteria expressed in terms of dual up/down degrees show that $\widetilde\Delta_{k}$ is bounded---and hence essentially self-adjoint on the natural core---whenever these degrees are finite (estimate~\eqref{eq:form-bounds}).

  At top degree $k=n$ the down term is the only one present. After a unitary conjugation, Theorem~\ref{thm:kn-esa} recasts it as an adjacency operator on the line-complex, with an explicit signed kernel, plus a diagonal potential, a form to which Schur's test applies directly. Neither route dominates the other. On an unweighted complex the line-complex bound never exceeds the Schur one, and the two agree exactly when some top simplex has all of its $n+1$ faces of maximal coface count, or, on an infinite complex, when a sequence of top simplices does so asymptotically --- which is the case at $n=1$ on a $d$-regular graph, where both give exactly $2d$. On a weighted complex the two are genuinely incomparable, the geometric mean $\sqrt{m_n(\sigma)m_n(\sigma')}$ produced by the symmetrization being smaller or larger than $m_n(\sigma')$ according to the local weight ratio; see Example~\ref{ex:two-routes} for both phenomena on a two-triangle complex. In all cases the main role of the reduction is to make the sharpness analysis of Section~\ref{sharp} possible. On the quality of these constants, see Remark~\ref{rem:form-bounds-sharpness}.
  \item \emph{The edge level, and weighted extensions.}
  Specializing the previous item to $n=k=1$: on any unweighted $d$-regular graph ($d\ge 2$), regarded as a $1$-complex, the edge block satisfies
  \[
  \|\widetilde\Delta_{1}\|\ \le\ 2d ,
  \]
  without any geometric or completeness hypothesis. This bound is the classical edge-degree estimate of Anderson and Morley in disguise (Remark~\ref{rem:anderson-morley}); what is new here is its derivation inside the Hodge framework, in the infinite weighted setting, and the sharpness analysis that accompanies it. Weighted extensions are controlled by the elementary constant $C_1/c_0$, through $\|\Delta_1\|\le (C_1/c_0)(\Delta(L(G))+2)$, where $0<c_0\le m_0$ and $m_1\le C_1<\infty$ are uniform comparability bounds for the vertex and edge weights (Proposition~\ref{prop:weighted-extension}).
  \item \emph{Orientation signs made explicit.}
  An ordered coloring yields an explicit unitary identification of the skew model with a model on unoriented simplices, built from canonical color-ordered representatives (Theorem~\ref{thm:unitary-color}). Its interest is computational: the orientation signs, abstract in the general theory, become a function of the colors alone. They cancel exactly in the bipartite top-degree case, which is also the case in which our edge bound turns out to be sharp --- provided the graph is amenable, bipartiteness alone being neither necessary nor sufficient for sharpness (Remark~\ref{rem:signs-bipartite}). See Section~\ref{subsec:bake-comparison} for the comparison with~\cite{BaKe}, and Remark~\ref{rem:unitary-correction} for the exact status of the color sign, which \emph{is} the intertwiner in a fixed reference basis, but towards the canonical-basis coboundary and not towards the unsigned symmetric one.
  \item \emph{Sharp constants via Floquet--Bloch.}
  We first settle when the bound $\|\widetilde\Delta_{1}\|\le 2d$ is attained: on every \emph{finite} unweighted $d$-regular bipartite graph, and on every infinite such graph that is moreover \emph{amenable}. Amenability is needed: on the $d$-regular tree, bipartite but not amenable, the conclusion fails (see the paragraph on bipartiteness at the end of Section~\ref{app:bloch}).
 The standard periodic lattices being amenable, the bipartite ones realize $2d$; on the three non-bipartite lattices of Table~\ref{tab:lattices-combined} the bound is strict, as the Floquet--Bloch symbols of Section~\ref{app:bloch} show by giving the exact norms. Non-bipartiteness alone does not force strictness on an infinite graph; see Section~\ref{app:bloch}.
\end{itemize}

\paragraph{Scope of the quantitative results.} It is worth stating plainly how these contributions are distributed across degrees, since they are not uniform. In every degree $k$ we obtain a boundedness criterion with an explicit constant (Lemma~\ref{lem:form-bounds}); no optimality is claimed for it, and Remark~\ref{rem:form-bounds-sharpness} says where its slack lies. At top degree $k=n$ the estimate acquires structure --- adjacency plus potential on the line-complex (Theorem~\ref{thm:kn-esa}) --- which is what makes a sharpness analysis possible at all. The sharpness analysis itself, together with the exact constants, is carried out for the edge block, that is at $n=1$: Proposition~\ref{prop:sharpness-bipartite}, Table~\ref{tab:lattices-combined} and Proposition~\ref{prop:top-bounded-ESA} all live there, and Remark~\ref{rem:anderson-morley} records that the resulting degree-$1$ bound is the classical estimate of Anderson and Morley, obtained here in the infinite and weighted setting. In degrees $k\ge2$ we prove boundedness and essential self-adjointness, not sharpness; the only quantitative illustration beyond degree~$1$ is Example~\ref{ex:two-routes}. Determining the exact constant in intermediate degrees, and the higher-degree analogue of Table~\ref{tab:lattices-combined}, are left open.

The graph case is well developed~\cite{KLW-monograph,KellerLenzESA,Wojciechowski,Woj-survey,HKLW,Jor-Pearse,Mas,Tor-Hamza,DM,Gol2,Go}; the higher-dimensional case of \emph{simplicial complexes}~\cite{AnTo,Ay,AAcT,BGJ,BSJ1,Che,HJ,EckmannHodge} is less so. The most directly comparable reference is the work of Bartmann and Keller~\cite{BaKe}; see Section~\ref{sec:related} for the comparison.

\subsection*{Cochain conventions: the skew default, and a colorable-complex reformulation}

The \emph{skew convention} is the only framework used in this paper --- in the abstract, the introduction, the preliminaries (Section~\ref{sec:prelim-comb}) and the whole body; full definitions are given there and in Section~\ref{sec:color}.

\paragraph{The skew convention (the default, and the only one used).}
A $k$-cochain $f$ is \emph{alternating}:
\[
f([x_{\pi(0)},\dots,x_{\pi(k)}])=(-1)^{\varepsilon(\pi)}f([x_0,\dots,x_k])\quad\text{for every permutation $\pi$.}
\]
The coboundary $d^k_{\mathrm{skew}}$ carries the alternating signs $(-1)^i$ when the $i$-th vertex is omitted. This is the classical cochain convention; the identity $d^{k+1}_{\mathrm{skew}}\circ d^k_{\mathrm{skew}}=0$ holds unconditionally (Remark~\ref{rem:dsq-zero}(i)). Every Hodge-theoretic notion in this paper (Laplacians, forms, Friedrichs realizations, boundedness and ESA criteria) is defined and proved for this convention, with no restriction on the complex beyond the standing assumptions recorded in Section~\ref{sec:prelim-comb}.

\paragraph{A cautionary remark on the unsigned symmetric coboundary.}
One may also consider \emph{fully symmetric} $k$-cochains, $f([x_{\pi(0)},\dots,x_{\pi(k)}])=f([x_0,\dots,x_k])$ for all $\pi$, with an unsigned coboundary $d^k_{\mathrm{sym}}$ (no $(-1)^i$ factors). This unsigned symmetric coboundary is \emph{not} used anywhere in the paper: in general $d^{k+1}_{\mathrm{sym}}\circ d^k_{\mathrm{sym}}\neq 0$, so the symmetric cochains do not form a cochain complex and the associated ``Laplacian'' is not a genuine Hodge Laplacian (Remark~\ref{rem:dsq-zero}(ii)). We mention it here only to distinguish it sharply from the reformulation below, with which it should not be confused.

\paragraph{Canonical-basis reformulation (colorable complexes only).}
Section~\ref{sec:color} introduces, exclusively for \emph{colorable} complexes (those whose $1$-skeleton admits a proper vertex coloring), a unitary
\[
U_k:\ell^2_{\mathrm{skew}}(m_k)\to\ell^2_{\mathrm{can}}(m_k),\qquad (U_k f)(|\sigma|):=f(\sigma_{\mathrm{can}}),
\]
that transports the skew theory onto the space $\ell^2_{\mathrm{can}}(m_k)$ of functions on \emph{unoriented} simplices, evaluated at a canonical color-ordered representative $\sigma_{\mathrm{can}}$ (Definition~\ref{def:canonical-basis}). Read on canonical representatives, the transported coboundary $d^k_{\mathrm{sym},\mathrm{can}}:=U_k\,d^k_{\mathrm{skew}}\,U_{k-1}^{-1}$ takes the form $\sum_i(-1)^i f((\tau_i)_{\mathrm{can}})$. This looks like the expression written on symmetric cochains, but differs from it twice over: its arguments are canonical representatives, and its alternating signs are indexed by the color rank rather than by an arbitrary ordering. It satisfies $d^{k+1}_{\mathrm{sym},\mathrm{can}}\circ d^k_{\mathrm{sym},\mathrm{can}}=0$ \emph{because} it is unitarily equivalent to the skew coboundary (Theorem~\ref{thm:unitary-color}(1)). This is a purely computational reformulation of the skew theory, not a second convention on a par with the skew one:
\(
U_k\,\Delta_{k,\mathrm{skew}}\,U_k^{-1}=\Delta^{\mathrm{can}}_{k,\mathrm{sym}}
\)
(Theorem~\ref{thm:unitary-color}(2)), so every spectral and ESA result obtained for the skew Laplacian transfers verbatim to $\Delta^{\mathrm{can}}_{k,\mathrm{sym}}$ on colorable complexes.

\medskip\noindent
\emph{Notational convention:} from Section~\ref{sec:prelim-comb} onward, the unadorned symbols $d^k,\delta^k,\Delta_k,\mathfrak q_k,\Cc^k_c(\Vc),\dots$ always denote the \emph{skew} objects; no convention subscript is carried. The canonical-basis operators of Section~\ref{sec:color} keep their explicit labels $d^k_{\mathrm{sym},\mathrm{can}}$ and $\Delta^{\mathrm{can}}_{k,\mathrm{sym}}$, and are never abbreviated; by Corollary~\ref{cor:transfer} every result below applies to them on a colorable complex.

\section{Related work}
\label{sec:related}

\subsection{Self-adjoint extensions of graph Laplacians}

\paragraph{Form-theoretic foundations.}
The treatment of regular Dirichlet forms on weighted graphs goes back to Keller and Lenz~\cite{KellerLenzESA}, who developed the correspondence between Dirichlet forms on locally finite weighted graphs and weighted graph Laplacians, and dealt with stochastic completeness of subgraphs. Haeseler, Keller, Lenz and Wojciechowski~\cite{HKLW} analyzed Dirichlet and Neumann boundary conditions in the infinite-graph context and characterized when the formal Laplacian admits a unique self-adjoint extension. Wojciechowski~\cite{Wojciechowski,Woj-survey} systematically studied stochastic completeness of graphs --- equivalently, conservativeness of the associated heat semigroup --- a property distinct from, but closely tied to, the uniqueness questions considered here. The monograph~\cite{KLW-monograph} gives a contemporary account of this material, including self-adjointness criteria and semigroup properties, and we refer to it for further references.

Jorgensen and Pearse~\cite{Jor-Pearse} introduced, from a different angle, a Hilbert-space approach to effective resistance metrics and graph Laplacians, focused on the transience/recurrence dichotomy. Masamune~\cite{Mas} obtained Liouville-type theorems for the formal Laplacian and applied them to essential self-adjointness.

\paragraph{Geometric criteria for ESA.}
Another line of work concerns \emph{geometric} criteria for essential self-adjointness, in regimes where boundedness fails but completeness-type conditions still ensure ESA. The first contributions are due to Torki-Hamza~\cite{Tor-Hamza}, who introduced an analogue of geodesic completeness for graphs (\emph{$\chi$-completeness}), and to Colin de Verdi\`ere, Torki-Hamza and Truc~\cite{CTT,CTT2}, who proved ESA results for combinatorial Schr\"odinger operators (with and without magnetic fields) under metric-completeness assumptions. The abstract prototype of arguments of this type is Chernoff's criterion~\cite{Ch}, which derives essential self-adjointness, of a generator and of its powers, from finite propagation speed for the associated hyperbolic evolution.

Earlier related work of Dodziuk and Mathai~\cite{DM} treated Kato's inequality and the asymptotic spectral properties of discrete magnetic Laplacians. Ann\'e and Torki-Hamza~\cite{AnTo} and Ayadi~\cite{Ay} extended these geometric criteria to operators on $k$-forms, and Ann\'e, Ayadi, Chebbi and Torki-Hamza~\cite{AAcT} to magnetic Laplacians on triangulations.

\paragraph{Adjacency and unbounded operators.}
A third line is concerned with the boundedness/unboundedness dichotomy itself. Gol\'enia~\cite{Go} characterized when adjacency matrices of locally finite graphs are unbounded, and studied in~\cite{Gol2} Hardy-type inequalities and asymptotic eigenvalue distributions. Gol\'enia and Truc~\cite{SF} treated the magnetic Laplacian on discrete cusps. Spectral properties on specific lattices --- limiting absorption, non-propagation, transport, including triangular and graphene-like geometries --- have been analyzed by Athmouni and coauthors~\cite{A,AD,ABDE,AEG,AEG1,AEGJ1,AEGJ2}.

\paragraph{Standard textbook references.}
Chung's monograph~\cite{Chung} is the standard reference for spectral graph theory in the finite case, and Colin de Verdi\`ere's lecture notes~\cite{CdV} develop the spectral theory of graph Laplacians with an emphasis on its geometric aspects. Of particular relevance here is the bound $\mathrm{spec}(\widetilde\Delta_0)\subset[0,2]$ of~\cite{Chung}, valid for the normalizing vertex weight $m_0(x)=\sum_y m_1(x,y)$; for a general vertex weight the corresponding bound is $2\sup_x d_{\Vc}(x)$, see Section~\ref{sec:normalization}. 

\subsection{Edge-degree bounds in spectral graph theory}

A distinct line of work, in the spectral theory of \emph{finite} graphs, meets the present paper at degree $1$. For an unweighted finite graph ($m_0\equiv m_1\equiv1$, so that $\Delta_0$ is the combinatorial Laplacian $D-A$ and \emph{not} the normalized operator of the preceding paragraph), Anderson and Morley~\cite{AM} proved the edge-degree estimate
\[
\lambda_{\max}(\Delta_0)\ \le\ \max\bigl\{\deg u+\deg v:\ \{u,v\}\ \text{an edge of }G\bigr\},
\]
whose relation to Theorem~\ref{thm:kn-esa} at $n=1$ is set out in Remark~\ref{rem:anderson-morley}.

Rojo, Soto and Rojo~\cite{RSR} and Das~\cite{Das} subsequently sharpened the estimate by subtracting the number of common neighbours of the two endpoints. The gain is strict precisely when the edges realizing the maximum lie in a triangle, which is the case on the non-bipartite lattices of Section~\ref{sharp}; see again Remark~\ref{rem:anderson-morley}.

\subsection{Hodge Laplacians on simplicial complexes}

The subject goes back to Eckmann~\cite{EckmannHodge}, whose paper of 1944/45 introduced the discrete Hodge decomposition underlying modern simplicial Laplacian theory. For finite weighted complexes, Horak and Jost~\cite{HJ} provide the standard reference, clarifying the role of normalization choices. The infinite case is more recent; it has been studied in~\cite{AnTo} (Gauss--Bonnet operator on infinite graphs), \cite{Ay} (spectra on forms), \cite{AAcT} (magnetic Laplacians on triangulations), \cite{BGJ} (adjacency matrix and discrete Laplacian on forms), \cite{BSJ1} ($2$-forms), and \cite{Che} ($2$-simplicial complexes).

In the finite case, the higher-degree analogue of the edge-degree bounds recalled above has been developed recently. Fan, Wu and Wang~\cite{FWW} bound the largest eigenvalue of the $i$-th up-Laplacian by comparison with the corresponding signless up-Laplacian, improving the estimate of~\cite{HJ} and generalizing~\cite{AM}, with an equality criterion formulated through the balancedness of the incidence-signed graph. The normalized counterpart is treated by Song, Wu and Fan~\cite{SWF}. Their signed graph plays, in that finite setting, the role of the signed adjacency operator $\mathcal A_\eta$ of Theorem~\ref{thm:kn-esa}, and their balancedness criterion is the higher-degree analogue of the sign computation of Remark~\ref{rem:signs-bipartite}.

\subsection{Comparison with~\cite{BaKe} and position of the present work}\label{subsec:bake-comparison}

Bartmann and Keller~\cite{BaKe} develop a form-theoretic framework for Hodge Laplacians on simplicial complexes, possibly non-locally finite, and give general criteria for uniqueness of self-adjoint realizations. In their terminology, the Laplacians studied here are the \emph{Dirichlet realizations} of~\cite{BaKe}, denoted $\Delta^H_D$ and associated with the closed form $Q^H_D$ obtained by closing $Q^H$ on the finitely supported functions. This is the Friedrichs extension of the formal Hodge Laplacian on the core of compactly supported cochains; more precisely, our block $\Delta_k^{\mathrm F}$ is the restriction of $\Delta^H_D$ to degree $k$, the operator of~\cite{BaKe} acting across all degrees at once. The present paper was prepared independently, and we learned of~\cite{BaKe} only at a late stage; we thank its authors for the correspondence that followed.

\begin{itemize}[leftmargin=2em]
\item \emph{Boundedness criteria.} Our ESA statements assume boundedness of the relevant degree quantities: the up/down degrees $D_k^{\downarrow},D_k^{\uparrow}$ in Lemma~\ref{lem:form-bounds}, the line-complex Schur degree $D$ and the potential $V_n$ in Theorem~\ref{thm:kn-esa}, and $\Delta(L(G))$ in Proposition~\ref{prop:top-bounded-ESA}. In each such regime the operators are bounded on $\ell^2(m_k)$ and self-adjointness is automatic. Boundedness is also addressed in~\cite[Theorem~4.1]{BaKe}, which gives an equivalence for finite-dimensional complexes such as ours and a sufficient criterion in general.
\item \emph{How the two compare.} The criteria obtained here are stated \emph{degree by degree}, for each block $\Delta_k$ separately, whereas~\cite[Theorem~4.1]{BaKe} characterizes boundedness of the full direct sum $\bigoplus_k\Delta_k$ at once. Since $n<\infty$ here and $\|\bigoplus_{k\le n}A_k\|=\max_k\|A_k\|$, the two criteria are qualitatively equivalent in our setting; what the results of this paper add is the explicit constant, not the boundedness criterion itself. See Remark~\ref{rem:BaKe-thm63}, which also relates Theorem~\ref{thm:kn-esa} to the signed Schr\"odinger-operator representation of~\cite[Section~3.4]{BaKe}.
\item \emph{The unbounded regime.} The ESA criteria of~\cite[Theorems~4.9 and~4.11]{BaKe} rest respectively on a uniformly positive weight with lower-bounded Forman curvature, and on a Gaffney-type condition of finiteness of metric balls. The framework of~\cite{BaKe} is in this respect broader than ours, requiring neither local finiteness nor boundedness; we do not treat the unbounded case.
\item \emph{What we add.} The present contribution is quantitative. It consists of Schur-type estimates in every degree; the adjacency--potential decomposition and line-complex reduction at top degree; the explicit determination of the orientation cocycle from an ordered coloring; a weighted extension under uniform comparability; and Floquet--Bloch computations giving the exact constants on periodic lattices, together with the sharpness criterion of Proposition~\ref{prop:sharpness-bipartite}. These give concrete geometry-free numerical bounds, for instance $\|\widetilde{\Delta}_{1}\|\le 2d$ for unweighted $d$-regular graphs with $d\ge 2$, complementing the abstract form-based theory of~\cite{BaKe}.

  At degree $1$, however, that bound is the classical estimate of Anderson and Morley~\cite{AM}, and it is not sharper than the refinements of~\cite{RSR,Das}; see Remark~\ref{rem:anderson-morley}. The novelty in that degree lies in the exact constants and in the sharpness criterion, not in the bound itself.
\end{itemize}

Compared with~\cite{HJ}, who work on \emph{finite} complexes, we carry the analysis over to infinite ones, where the operators act on infinite-dimensional Hilbert spaces. (The phrase ``finite-dimensional complex'' used in the first bullet above refers, by contrast, to the dimension $n$ of the complex, which is finite throughout this paper.)

A word on what is, and what is not, new here. We do not claim novelty for the \emph{existence} of a unitary identification between the skew model and a model on unoriented simplices: such an identification is available on an arbitrary complex, by pointwise multiplication with a suitable sign~\cite[Theorem~2.8]{BaKe}. What Theorem~\ref{thm:unitary-color} contributes is an explicit and computable one, whose canonical representatives are read off from an ordered coloring (Remark~\ref{rem:unitary-correction}); it is not what produces the constants of Section~\ref{sharp}, which come from Lemma~\ref{lem:iso-01} and the Bloch symbols.

Nor are the numerical values themselves new as band-structure data: the extrema of the adjacency symbol on $\Z^r$, on the triangular, BCC, FCC and Kagome lattices are classical. What we claim is their transfer to the Hodge edge block, together with the sharpness criterion of Proposition~\ref{prop:sharpness-bipartite}, in the simplicial setting.

\section{Preliminaries}\label{sec:prelim-comb}
\subsection{Combinatorial setting}

\subsubsection{Graphs and simplicial complexes}

We begin with some standard definitions from graph theory and the simplicial structures they induce.

\begin{definition}[Weighted graph]\label{def:weighted-graph}
A (non-oriented) weighted graph is a triple $\Gc=(\Vc,m_0,m_1)$ where
\begin{itemize}
  \item $\Vc$ is a countable vertex set;
  \item $m_0:\Vc\to(0,\infty)$ is a vertex weight;
  \item $m_1:\Vc\times\Vc\to[0,\infty)$ is a \emph{symmetric} edge weight vanishing on the diagonal, i.e., $m_1(x,y)=m_1(y,x)$ for all $x,y$ and $m_1(x,x)=0$ for all $x$ (so that $\Gc$ has no loops).
\end{itemize}
The (undirected) edge set is
\[
  \Ec:=\bigl\{\{x,y\}\subset\Vc:\ m_1(x,y)>0\bigr\}.
\]
We write $x\sim y$ if $m_1(x,y)>0$ and call $x,y$ \emph{neighbors}. The neighbor set of $x$ is
\[
  \Nc_\Gc(x):=\{y\in\Vc:\ x\sim y\}.
\]
The graph is \emph{locally finite} if $\#\Nc_\Gc(x)<\infty$ for all $x\in\Vc$. The \emph{normalized degree} --- the total edge weight at $x$ divided by the vertex weight --- is
\[
  d_{\Vc}(x):=\frac{1}{m_0(x)}\sum_{y\in\Vc} m_1(x,y).
\]
\end{definition}

Fix $n\in\N$ with $n\ge 1$. We next formalize the $n$-simplicial structure generated by complete subgraphs (cliques) of size at most $n\!+\!1$. Two remarks on the role of $n$. First, it is a choice, not a property of $\Gc$: a graph containing a triangle has cliques of size $2$ and of size $3$ alike, and fixing $n$ \emph{truncates} the clique complex at dimension $n$. Thus $\Sc_1$ carries no $2$-simplex even when $\Gc$ is full of triangles, and $\Sc_2$ carries no tetrahedron even when $\Gc$ contains $K_4$. Second, the case $n=1$ is allowed and is not degenerate: $\Sc_1$ is then the weighted graph $\Gc$ itself, and it is the setting in which the sharp edge bounds of Section~\ref{sharp} are obtained, by applying the top-degree results with $n=1$.

\begin{definition}[Cliques and oriented simplices]\label{def:Fk-clique}
For $k\in\{0,1,\dots,n\}$ let
\[
  \Fc_k:=\Bigl\{(x_0,\dots,x_k)\in\Vc^{k+1}:\ \#\{x_0,\dots,x_k\}=k+1\ \text{and}\ \{x_i,x_j\}\in\Ec\ \forall i\neq j\Bigr\}
\]
be the set of ordered $(k\!+\!1)$-tuples forming a clique (a complete subgraph) in $\Gc$.
Let $\approx$ be the equivalence relation on $\Fc_k$ identifying two tuples that differ by an \emph{even} permutation (the symbol $\sim$ is reserved for adjacency). The set of oriented $k$-simplices is
\[
  T_k:=\Fc_k/\!\approx,
\]
whose elements are written $\sigma=[x_0,\dots,x_k]$. Writing $\varepsilon(\pi)\in\{0,1\}$ for the parity of a permutation $\pi$ of $\{0,1,\dots,k\}$, one has
\[
[x_{\pi(0)},\dots,x_{\pi(k)}]=[x_0,\dots,x_k]\ \text{ if }\varepsilon(\pi)=0,
\]
while for $\varepsilon(\pi)=1$ the left-hand side is the \emph{opposite} class; for $k\ge1$ each clique therefore carries exactly two classes. Note that $T_k$ is a set, so the sign $(-1)^{\varepsilon(\pi)}$ does not act on its elements: it enters only through the values of cochains, see the definition of $\Cc^k_{\mathrm{skew}}(\Vc)$ in Section~\ref{sec:prelim-analytic}. We write $|\sigma|$ for the unoriented $k$-simplex underlying $\sigma\in T_k$, i.e.\ for the underlying $(k+1)$-clique.
\end{definition}

\begin{definition}[Higher weights and face degrees]\label{def:weights}
For each $k\in\{0,1,\dots,n\}$ let $m_k:\Vc^{k+1}\to[0,\infty)$ be a symmetric function (defined on tuples in $\Vc^{k+1}$, not on equivalence classes in $T_k$) such that, \emph{for every permutation $\pi$ of $\{0,1,\dots,k\}$},
\[
  m_k(x_{\pi(0)},\dots,x_{\pi(k)})=m_k(x_0,\dots,x_k),\qquad
  m_k(x_0,\dots,x_k)>0\ \Longleftrightarrow\ (x_0,\dots,x_k)\in\Fc_k.
\]
Thus $m_0$ and $m_1$ are the given vertex and edge weights, and for $k\ge 2$ the values $m_k$ weight the $(k\!+\!1)$-cliques. In particular $m_k$ is strictly positive on every $k$-simplex. By the symmetry assumption, $m_k$ descends to unoriented simplices: we write $m_k(|\sigma|)$, or simply $m_k(\sigma)$, for the common value of $m_k$ on the ordered representatives of $\sigma\in T_k$.
For $(x_1,\dots,x_k)\in \Fc_{k-1}$ we set the common neighbor set
\[
  F_{(x_1,\dots,x_k)}:=\bigcap_{j=1}^{k}\Nc_\Gc(x_j),
\]
and define the $(k\!-\!1)$-face degree by
\[
  d_{k-1}(x_1,\dots,x_k):=
  \frac{1}{m_{k-1}(x_1,\dots,x_k)}\sum_{x_{k+1}\in F_{(x_1,\dots,x_k)}} m_k(x_1,\dots,x_k,x_{k+1}).
\]
When all higher weights $m_k$ are $\{0,1\}$-valued (simple complex), this reduces to $d_{k-1}(x_1,\dots,x_k)=\#F_{(x_1,\dots,x_k)}$; and for $k=1$ one recovers the normalized degree of Definition~\ref{def:weighted-graph}, $d_0=d_{\Vc}$. By symmetry of the weights, $d_{k-1}$ descends to unoriented $(k-1)$-simplices, and its supremum is the raw down-degree $D_k^{\downarrow,\mathrm{HJ}}$ of~\eqref{eq:raw-degrees} used in the Schur bounds of Section~\ref{sec:normalization}.
\end{definition}

\emph{A note on the letters $\Delta$ and $D$.} The letter $\Delta$ denotes the Hodge blocks $\Delta_k$ and, in Section~\ref{sharp} only, the maximum degree $\Delta(G)$ of a graph and $\Delta(L(G))$ of its line graph; the two are always distinguished by their argument, a degree index for the former and a graph for the latter. The letter $D$ denotes the degrees $D_k^{\downarrow},D_k^{\uparrow}$ and the Schur degree $D(\sigma)$, and, in the standard notation for domains, $D[\mathfrak q]$ and $D(A)$ for a form and an operator respectively; the Gauss--Bonnet operator is the calligraphic $\mathcal D$.

\emph{A note on the letter $d$.} Degrees carry a subscript ($d_{\Vc}$ and $d_{k-1}$ above) and the coboundary a superscript ($d^k$, Section~\ref{sec:prelim-analytic}), so the two never collide. The one exception is Section~\ref{sharp}, where an unadorned $d$ denotes the vertex degree of a periodic lattice; the coboundary does not appear there, and the point is recalled at the start of that section.

\begin{definition}[Weighted oriented $n$-simplicial complex]\label{def:simplicial-complex}
$\Sc_n:=(\Vc,(m_k)_{0\le k\le n})$ associated with $\Gc=(\Vc,m_0,m_1)$ is called a (weighted, oriented) $n$-simplicial complex. Since, by the preceding definitions, $\Fc_k$ consists of all $(k+1)$-cliques of $\Gc$ and $m_k$ is positive precisely on these cliques, $\Sc_n$ is entirely determined by its $1$-skeleton $\Gc$: this is a \emph{flag complex} (also called a \emph{clique complex}). In particular, ``hollow'' faces -- e.g., three mutually adjacent vertices without an accompanying $2$-simplex -- cannot occur, in contrast with general abstract simplicial complexes.
For $n=2$ one speaks of a (weighted) triangulation; see, e.g., \cite{BSJ1,Che}.
\end{definition}

\paragraph{Standing assumptions on local finiteness.} Local finiteness is \emph{not} a blanket hypothesis in this paper, and it is worth saying precisely what is needed where.

The combinatorial constructions of the present section --- cochain spaces, coboundaries, formal adjoints, Green's identity --- require nothing beyond countability of $\Vc$. The form-theoretic material of Section~\ref{section3} --- closability, Friedrichs realizations, the energy identities --- requires in addition the \emph{pointwise finiteness of the up-degrees},
\begin{equation}\label{eq:pointwise-updeg}
d_k(\sigma)\;=\;\frac{1}{m_k(\sigma)}\sum_{w\in F_\sigma}m_{k+1}(\sigma,w)\;<\;\infty
\qquad\text{for every $k$-simplex $\sigma$ and every }0\le k\le n-1 .
\end{equation}
Countability alone would not suffice. Indeed, for $f$ the skew indicator of a single $k$-simplex $\sigma$ one computes $\|d^{k+1}f\|_{k+1}^2=\sum_{w\in F_\sigma}m_{k+1}(\sigma,w)=m_k(\sigma)\,d_k(\sigma)$. Thus~\eqref{eq:pointwise-updeg} is exactly the condition under which the form $\mathfrak q$ of Section~\ref{section3} is finite on the core $\Cc^\bullet_c(\Vc)$; without it $\mathfrak q$ is not densely defined, and there is nothing to close. The down part carries no such requirement, $\delta^k$ always mapping compactly supported cochains to compactly supported cochains.

Condition~\eqref{eq:pointwise-updeg} is implied by, but strictly weaker than, local finiteness, and it is implied a fortiori by the uniform bound $D_k^{\uparrow}<\infty$ of Lemma~\ref{lem:form-bounds}. That lemma is therefore also proved without local finiteness: the passage from the form bound to the operator bound is arranged so as not to require that $\widetilde\Delta_k$ preserve the core (see the end of that proof).

By contrast, local finiteness \emph{is} assumed throughout Section~\ref{sharp}, where it is stated explicitly. It is \emph{not} needed in Proposition~\ref{prop:schur}, where the hypothesis $D\in L^\infty$ alone makes every sum absolutely convergent. Nor is it needed in Theorem~\ref{thm:kn-esa}, where assumptions (A1)--(A2) make it superfluous; we retain it there in the statement, for the reasons given in the paragraph following that proof.

\begin{definition}[Line-complex of $\Sc_n$ at top degree]\label{def:line-complex}
The \emph{line-complex} of $\Sc_n$ (at top degree $n$) is the graph $(\widehat{\Vc},\widehat{\Ec})$ whose vertices are the unoriented $n$-simplices $|\sigma|$ of $\Sc_n$. Two \emph{distinct} such vertices $|\sigma|\neq|\sigma'|$ are adjacent in $\widehat{\Ec}$ if and only if $\sigma$ and $\sigma'$ share an $(n-1)$-face in $\Sc_n$; we then write $\sigma\sim\sigma'$. The symbol $\sim$ thus serves for adjacency both in $\Gc$ (Definition~\ref{def:weighted-graph}) and in the line-complex, the ambient graph being determined by its arguments. There is no loop at $|\sigma|$: the diagonal contribution is carried by the potential $V_n$ of Theorem~\ref{thm:kn-esa}, not by $\widehat\Ec$. Note that two distinct $n$-simplices have at most $n$ common vertices, hence at most one common $(n-1)$-face, so the shared face is unique when it exists. For $n=1$ this reduces to the classical line graph $L(G)$ of $G$.
\end{definition}

\subsection{Hilbert cochains and Hodge operators}
\label{sec:prelim-analytic}

\subsubsection{Function spaces}

We endow the cochain spaces with natural Hilbert structures.

\medskip
\noindent\emph{$0$-cochains.}
Set
\[
  \Cc^0(\Vc):=\{f:\Vc\to\C\},\qquad
  \Cc_c^0(\Vc):=\{f\in\Cc^0(\Vc):\ \#\,\mathrm{supp}\,f<\infty\}.
\]
The weighted space
\[
  \ell^2(m_0):=\Bigl\{f\in\Cc^0(\Vc):\ \|f\|_0^2:=\sum_{x\in\Vc} m_0(x)\,|f(x)|^2<\infty\Bigr\}
\]
is a Hilbert space with inner product
\(
  \langle f,g\rangle_0:=\sum_{x\in\Vc} m_0(x)\,f(x)\overline{g(x)}.
\)

\medskip
\noindent\emph{$k$-cochains, $k\ge 1$.}
For $k\in\{1,\dots,n\}$ we define, the conditions below being required for every permutation $\pi$ of $\{0,1,\dots,k\}$,
\begin{align*}
  \Cc^k_{\mathrm{skew}}(\Vc)
  &:=\Bigl\{f:T_k\to\C\ :\ f([x_{\pi(0)},\dots,x_{\pi(k)}])=(-1)^{\varepsilon(\pi)} f([x_0,\dots,x_k])\Bigr\},\\
  \Cc^k_{\mathrm{sym}}(\Vc)
  &:=\Bigl\{f:T_k\to\C\ :\ f([x_{\pi(0)},\dots,x_{\pi(k)}])= f([x_0,\dots,x_k])\Bigr\}.
\end{align*}
The subspaces of compactly supported cochains, i.e.\ those $f$ with $\#\,\mathrm{supp}\,f<\infty$, are denoted $\Cc^{k}_{c,\mathrm{skew}}(\Vc)$ and $\Cc^{k}_{c,\mathrm{sym}}(\Vc)$ respectively. The skew space is the default: unless an explicit label says otherwise, $\Cc^k(\Vc):=\Cc^k_{\mathrm{skew}}(\Vc)$ and $\Cc^{k}_{c}(\Vc):=\Cc^{k}_{c,\mathrm{skew}}(\Vc)$. The symmetric space is recorded only so that $d^k_{\mathrm{sym}}$ can be written down in Remark~\ref{rem:dsq-zero}(ii).

Both spaces carry the same weighted inner product:
\begin{equation}
  \label{eq:2-2-1}
  \langle f,g\rangle_k
  := \frac{1}{(k+1)!}\sum_{(x_0,\dots,x_k)\in \Fc_k}
     m_k(x_0,\dots,x_k)\, f([x_0,\dots,x_k])\,\overline{g([x_0,\dots,x_k])}.
\end{equation}
The corresponding Hilbert spaces are denoted $\ell^2(m_k)$ (the variant notations $\ell^2_{\mathrm{skew}}(m_k)$ and $\ell^2_{\mathrm{can}}(m_k)$ are used in Section~\ref{sec:color}, where distinguishing the skew space from the canonical-basis space is essential). We write $\ell^2$, without argument, for the \emph{unweighted} space, that is, for $\ell^2(m_k)$ with $m_k$ replaced by the constant weight $1$ and the same normalization $\frac{1}{(k+1)!}\sum_{\Fc_k}$. With this convention multiplication by $m_k^{-1/2}$ is a unitary from $\ell^2$ onto $\ell^2(m_k)$, which is what makes the normalized blocks of Section~\ref{sec:normalization} unitary transports.

\medskip
\noindent\emph{Total space.}
We set
\[
  \Hc:=\bigoplus_{k=0}^{n} \ell^2(m_k),\qquad
  \|f\|_{\Hc}^2:=\sum_{k=0}^{n}\|f^{(k)}\|_k^2,\quad f=(f^{(0)},\dots,f^{(n)}).
\]

\begin{lemma}
\label{lem:density-Cc}
For each $k\in\{0,\dots,n\}$, the space $\Cc^{k}_{c}(\Vc)$ is dense in $\ell^2(m_k)$. (Recall from Definition~\ref{def:weights} that $m_k$ is automatically strictly positive on every $k$-simplex.)
\end{lemma}

\begin{proof}
Since $\Vc$ is countable, so is the set of unoriented $k$-simplices; fix an exhaustion $K_1\subset K_2\subset\cdots$ of it by finite sets. For $f\in\ell^2(m_k)$ put $f_N:=f\,\mathbf 1_{K_N}$, which is again a skew cochain, since the truncation is by a set of \emph{unoriented} simplices and therefore does not interfere with the behaviour under permutations. Then $f_N\in\Cc^k_c(\Vc)$ and $\|f-f_N\|_k^2=\sum_{|\sigma|\notin K_N}m_k(|\sigma|)|f(\sigma)|^2\to0$ by dominated convergence, the series $\sum_{|\sigma|}m_k(|\sigma|)|f(\sigma)|^2$ being convergent.
\end{proof}

\subsubsection{Cochain operators and the relation $d^{k+1}\circ d^k=0$}

Let $k\in\{1,\dots,n\}$. The coboundary maps --- which carry no weight, the weights entering only through the inner products~\eqref{eq:2-2-1} and hence through the adjoints $\delta^k$ --- are defined on finitely supported cochains by
\begin{align*}
  &d^{k}_{\mathrm{skew}} : \Cc^{k-1}_{c,\mathrm{skew}}(\Vc)\longrightarrow \Cc^{k}_{\mathrm{skew}}(\Vc),\\
  &\hspace*{2cm}d^{k}_{\mathrm{skew}}f([x_0,\dots,x_k]) := \sum_{i=0}^{k} (-1)^i f([x_0,\dots,\widehat{x_i},\dots,x_k]),\\
 & d^{k}_{\mathrm{sym}} : \Cc^{k-1}_{c,\mathrm{sym}}(\Vc)\longrightarrow \Cc^{k}_{\mathrm{sym}}(\Vc),\\
  &\hspace*{2cm} d^{k}_{\mathrm{sym}}f([x_0,\dots,x_k]) := \sum_{i=0}^{k} f([x_0,\dots,\widehat{x_i},\dots,x_k]).
\end{align*}
(Here $\widehat{x_i}$ means that $x_i$ is omitted.) Note the codomains: $d^kf$ is supported on the cofaces of $\mathrm{supp}\,f$, of which there may be infinitely many, so $d^k$ maps $\Cc^{k-1}_c(\Vc)$ into $\Cc^{k}(\Vc)$, and into $\Cc^{k}_c(\Vc)$ exactly when every $(k-1)$-simplex has finitely many cofaces, i.e.\ $\#F_\tau<\infty$ for all $\tau\in T_{k-1}$. For $k=1$ this is local finiteness; for a fixed $k\ge2$ it is strictly weaker, and local finiteness is the condition that makes it hold simultaneously in all degrees. The formal adjoint $\delta^k$ below carries no such caveat: $\delta^kg$ is supported on the faces of $\mathrm{supp}\,g$, finitely many, so $\delta^k$ always maps $\Cc^{k}_c(\Vc)$ into $\Cc^{k-1}_c(\Vc)$.

Both displayed formulas in fact make sense pointwise on \emph{arbitrary} cochains, each value being a sum of $k+1$ terms, and we use them in that generality whenever convenient. In particular the composition $d^{k+1}\circ d^k$, and the identity $d^{k+1}d^k=0$ of Remark~\ref{rem:dsq-zero} below, are to be read pointwise on $\Cc^{k-1}(\Vc)$. Restricting the coboundaries to compactly supported cochains matters only for the Hilbert-space discussion of Section~\ref{section3}, and it is there that condition~\eqref{eq:pointwise-updeg} enters.

\medskip

The relation $d^{k+1}\circ d^k=0$ requires care, since it holds in only one of the two conventions --- which is precisely why only the skew one is used in this paper.

\begin{remark}\label{rem:dsq-zero}
The identity $d^{k+1}\circ d^k=0$ holds in different generalities depending on the convention:
\begin{enumerate}[leftmargin=2em,label=(\roman*)]
\item In the \emph{skew} convention, the standard cochain calculation shows that
\(
d^{k+1}_{\mathrm{skew}}\,d^{k}_{\mathrm{skew}}=0
\)
holds without any further hypothesis (this is the classical simplicial cochain complex; see, e.g., \cite{EckmannHodge,HJ}).
\item In the \emph{symmetric} convention the identity $d^{k+1}_{\mathrm{sym}}\,d^{k}_{\mathrm{sym}}=0$ does \emph{not} hold. The reason is immediate from the definitions: for a symmetric $f$ and $\sigma=[x_0,\dots,x_{k+1}]$,
\[
d^{k+1}_{\mathrm{sym}}d^{k}_{\mathrm{sym}}f(\sigma)=\sum_{j}\sum_{i\neq j}f\bigl([\dots,\widehat{x_i},\dots,\widehat{x_j},\dots]\bigr)=2\sum_{i<j}f\bigl([\dots,\widehat{x_i},\dots,\widehat{x_j},\dots]\bigr),
\]
each unordered pair being counted twice \emph{with the same sign} instead of with opposite signs. The right-hand side is a positive combination of the values of $f$ on the $(k-1)$-faces and does not vanish in general. This failure persists on colorable complexes when $d^k_{\mathrm{sym}}$ is read as acting on arbitrary symmetric cochains: colorability by itself does \emph{not} restore the identity.

What does hold is the following. Fix an ordered coloring, which by Definition~\ref{def:coloring} is no restriction on a countable complex. Each simplex then has a canonical color-ordered representative $\sigma_{\mathrm{can}}$, and the skew coboundary read on those representatives is, directly from its definition,
\[
\bigl(d^k_{\mathrm{skew}}f\bigr)(\sigma_{\mathrm{can}})\;=\;\sum_{i=0}^k(-1)^i f\bigl((\tau_i)_{\mathrm{can}}\bigr),
\]
where $(\tau_i)_{\mathrm{can}}$ is the canonical representative of the $i$-th $(k-1)$-face (the order inherited from a color-ordered tuple is again color-ordered). We \emph{define} $d^k_{\mathrm{sym},\mathrm{can}}$ to be the operator on $\ell^2_{\mathrm{can}}$ given by this same expression (Definition~\ref{def:canonical-basis} below), and $d^{k+1}_{\mathrm{sym},\mathrm{can}}\circ d^k_{\mathrm{sym},\mathrm{can}}=0$ then follows from (i) by transport through the unitary intertwining of Theorem~\ref{thm:unitary-color}. Note that $d^k_{\mathrm{sym},\mathrm{can}}$ differs from the unsigned symmetric coboundary $d^k_{\mathrm{sym}}$ defined above: the alternating signs $(-1)^i$ are intrinsic to the canonical-basis transport.
\end{enumerate}
Throughout this paper we work in the skew convention, for which case~(i) applies with no hypothesis on the complex. In particular $\Delta_{k}=d^k\delta^k+\delta^{k+1}d^{k+1}$ is a genuine Hodge Laplacian, and the cross-terms $d^2$ and $\delta^2$ in the expansion of $\mathcal D^{\,2}=(d+\delta)^2$ vanish. For $\delta^2$ this follows by duality: for compactly supported $f,g$, Green's identity gives $\langle\delta^{k}\delta^{k+1}g,f\rangle_{k-1}=\langle g,d^{k+1}d^{k}f\rangle_{k+1}=0$, and since the $(k-1)$-simplex indicators are compactly supported and $m_{k-1}>0$, $\delta^{k}\delta^{k+1}g$ vanishes identically. Case~(ii) transports all of this to the canonical-basis operators, once an ordered coloring is fixed (Corollary~\ref{cor:transfer}).
\end{remark}

The (formal) adjoints $\delta^k$ are defined by
\[
  \delta^k:\Cc^{k}_{c}(\Vc)\longrightarrow \Cc^{k-1}_{c}(\Vc),\qquad
  \langle d^k f, g\rangle_k = \langle f, \delta^k g\rangle_{k-1}
  \ \ \text{for all compactly supported } f,g.
\]
Both sides are finite sums, $g$ and $f$ being compactly supported, so the relation is meaningful with no hypothesis on the complex. It determines $\delta^kg$ uniquely: if $\langle f,h\rangle_{k-1}=\langle f,h'\rangle_{k-1}$ for every $f\in\Cc^{k-1}_c(\Vc)$, then testing against the indicators of single $(k-1)$-simplices and using $m_{k-1}>0$ gives $h=h'$. Existence is not an assumption either: it is the content of Lemma~\ref{lem:2-3-1}, whose proof computes $\langle d^kf,g\rangle_k$ directly and exhibits the operator, thereby establishing the displayed relation rather than presupposing it.

\begin{lemma}\label{lem:2-3-1}
Let $k\in\{1,\dots,n\}$. For every $g\in \Cc^{k}_{c}(\Vc)$ and every $(x_1,\dots,x_k)\in \Fc_{k-1}$,
\begin{equation}\label{eq:2-3-1}
\begin{split}
 & \delta^k(g)([x_1,\dots,x_k])
  = \frac{1}{m_{k-1}(x_1,\dots,x_k)}
    \sum_{x_{0}\in F_{(x_1,\dots,x_k)}}
    m_k(x_{0},x_1,\dots,x_k)\\
  &\hspace*{4cm}\times  g([x_{0},x_1,\dots,x_k]).
  \end{split}
\end{equation}
The newly summed vertex $x_0$ appears in \emph{first} position inside $g$: this is what captures the correct sign in the formal adjoint relation $\langle d^k f,g\rangle_k=\langle f,\delta^k g\rangle_{k-1}$.
\end{lemma}

\begin{proof}
Using \eqref{eq:2-2-1} and the definition of $d^k$, expand
\begin{align*}
  &\langle d^k f, g\rangle_k
  = \frac{1}{(k+1)!}\sum_{(x_0,\dots,x_k)\in \Fc_k}
    m_k(x_0,\dots,x_k)\Bigl(\sum_{i=0}^k (-1)^i\, f([x_0,\dots,\widehat{x_i},\dots,x_k])\Bigr)\\
    &\hspace*{5cm}\times\overline{g([x_0,\dots,x_k])}.
\end{align*}
We regroup the terms by the $(k-1)$-face $\tau_i:=[x_0,\dots,\widehat{x_i},\dots,x_k]$: the $k$-simplex $[x_0,\dots,x_k]$ is a coface of $\tau_i$, obtained by inserting the additional vertex $z:=x_i$ at position $i$. By the antisymmetry of $g$,
\[
g([x_0,\dots,x_k])=(-1)^i\,g([z,x_0,\dots,\widehat{x_i},\dots,x_k])=(-1)^i\,g([z,\tau_i]),
\]
since moving $z$ from position $i$ to position $0$ requires exactly $i$ adjacent transpositions, each contributing a sign $-1$; equivalently $\overline{g([x_0,\dots,x_k])}=(-1)^i\,\overline{g([z,\tau_i])}$. Multiplying by the sign $(-1)^i$ carried by $f(\tau_i)$ gives $(-1)^i\cdot(-1)^i=1$: the alternating sign of the coboundary is exactly cancelled by the sign acquired by $\overline g$ under this reordering, so that the summand indexed by $((x_0,\dots,x_k),i)$ equals
\[
f(\tau_i)\;m_k(z,\tau_i)\;\overline{g([z,\tau_i])},
\]
in which the inserted vertex stands in \emph{first} position inside $g$, consistently with the statement of the lemma.

It remains to account for the combinatorial factor. The summand just displayed depends only on the ordered $k$-tuple representing $\tau_i$ and on the inserted vertex $z$, not on the position $i$ at which $z$ was inserted. Conversely, a given pair $\bigl((x_1,\dots,x_k),z\bigr)$ with $(x_1,\dots,x_k)\in\Fc_{k-1}$ and $z\in F_{(x_1,\dots,x_k)}$ arises from exactly $k+1$ pairs $\bigl((x_0,\dots,x_k),i\bigr)$, namely the $k+1$ positions at which $z$ may be inserted into $(x_1,\dots,x_k)$. Hence the prefactor becomes $(k+1)/(k+1)!=1/k!$, and we obtain
\begin{align*}
 & \langle d^k f, g\rangle_k
  = \frac{1}{k!}\sum_{(x_1,\dots,x_k)\in \Fc_{k-1}}
    f([x_1,\dots,x_k])\\
    &\hspace*{3cm} \times\overline{\sum_{x_{0}\in F_{(x_1,\dots,x_k)}}
      m_k(x_{0},x_1,\dots,x_k)\, g([x_{0},x_1,\dots,x_k])}.
\end{align*}
Comparing with
$$\langle f,\delta^k g\rangle_{k-1} = \frac{1}{k!}\sum_{(x_1,\dots,x_k)\in\Fc_{k-1}} m_{k-1}(x_1,\dots,x_k)\, f([x_1,\dots,x_k])\,\overline{\delta^k g([x_1,\dots,x_k])}$$
gives \eqref{eq:2-3-1}.

Finally, $\delta^k$ does map skew cochains to skew cochains. Let $\pi$ permute $(x_1,\dots,x_k)$. In~\eqref{eq:2-3-1} the weights $m_{k-1}(x_1,\dots,x_k)$ and $m_k(x_0,x_1,\dots,x_k)$ are unchanged, being symmetric functions of their arguments, and the common-neighbor set $F_{(x_1,\dots,x_k)}$ is unchanged as well; only $g([x_0,x_1,\dots,x_k])$ is affected, and it picks up the sign $(-1)^{\varepsilon(\pi)}$ since $g$ is alternating and $\pi$ acts on the last $k$ entries. Hence $\delta^kg$ is alternating, so that $\delta^k:\Cc^{k}_{\mathrm{skew}}(\Vc)\to\Cc^{k-1}_{\mathrm{skew}}(\Vc)$.
\end{proof}

\section{Laplacians via the Gauss--Bonnet operator}\label{section3}

\subsection{Definition and block structure}\label{subsec:block-structure}

\paragraph{Convention for Sections~\ref{section3}--\ref{sharp}.} All operators below are the skew ones, written without any convention subscript; $d^{k+1}\,d^{k}=0$ holds unconditionally (Remark~\ref{rem:dsq-zero}(i)), so the cross-terms $d^2$ and $\delta^2$ in the Hodge Laplacian below vanish. On a colorable complex, every statement of these sections transfers verbatim to the canonical-basis operators of Section~\ref{sec:color} (Corollary~\ref{cor:transfer}); we do not repeat this in the body of the subsequent sections.

\subsubsection{Local coordinate formulas}\label{subsubsec:local-formulas}
Let $\Sc_n=(\Vc,(m_k)_{0\le k\le n})$ be a weighted oriented $n$-simplicial complex.

Define the Gauss--Bonnet operator
\[
  \mathcal D := d + \delta \quad\text{on }\ \Cc_{c}^0(\Vc)\oplus\Cc_{c}^1(\Vc)\oplus\cdots\oplus \Cc_{c}^n(\Vc),
\]
where 
\begin{align*}
    d(f_0, \ldots , f_{n}) = (0, d^1f_0, d^2f_1, \ldots , d^nf_{n-1}),
\end{align*}
for all $(f_0, \ldots , f_{n})\in\Cc_{c}^0(\Vc)\oplus\Cc_{c}^1(\Vc)\oplus\cdots\oplus\Cc_{c}^n(\Vc)$,
and $\delta$ is given by
\begin{align*}
    \delta(f_0, \ldots , f_{n})=(\delta^1f_1, \delta^2f_2, \ldots , \delta^nf_{n},0).
\end{align*}
As in Section~\ref{sec:prelim-comb}, $\delta$ maps $\Cc^\bullet_c(\Vc)$ into itself, whereas $d$ --- and hence $\mathcal D$ --- maps it into $\Cc^\bullet(\Vc)$, the image being compactly supported in all degrees as soon as the complex is locally finite. All the expressions below are read pointwise, and under the standing assumption~\eqref{eq:pointwise-updeg} they moreover define elements of $\Hc$ when applied to compactly supported cochains, which is what makes the form $\mathfrak q$ below finite on the core.
The (Hodge) Laplacian is
\[
  L \ :=\ \mathcal D^{\,2} \ =\ d\,\delta + \delta\, d.
\]
Since $d^{k+1}\,d^k=0$ in the skew convention (Remark~\ref{rem:dsq-zero}(i)), the cross terms $d^2$ and $\delta^2$ vanish, and $L$ is block-diagonal by degree:
\begin{equation}\label{eq:2-4-1}
  L \ =\ \bigoplus_{k=0}^{n} \Delta_{k},\qquad
  \Delta_{k} \ :=\ d^{k}\,\delta^{k} \ +\ \delta^{k+1}\,d^{k+1},
\end{equation}
with the conventions $\delta^{0}=0$ and $d^{n+1}=0$.
Explicitly, for $f=(f^{(0)},\dots,f^{(n)})$,
\[
  L f \ =\ \bigl(\delta^1 d^1 f^{(0)},\ (d^1 \delta^1+\delta^2 d^2) f^{(1)},\ \dots,\ d^n \delta^n f^{(n)}\bigr).
\]

\begin{proposition}[Local coordinate formulas for the block Laplacian]\label{rmk:2-4-1}
Fix $k\in\{0,\dots,n\}$ and write a $k$-cochain as $f([x_0,\dots,x_k])$. The two displays below are to be read with the conventions of Corollary~\ref{cor:3-1}: at $k=0$ the down term is absent, $\delta^0:=0$, and no weight $m_{-1}$ occurs; at $k=n$ the up term is absent, $d^{n+1}:=0$, and no weight $m_{n+1}$ occurs.
Recall the extension set of Definition~\ref{def:weights}, written here with the index range shifted by one:
\begin{align*}
 & F_{(x_0,\dots,x_{k-1})}:=\bigcap_{j=0}^{k-1} \Nc_\Gc(x_j).\end{align*}
Its elements are exactly the vertices forming a $k$-simplex with $[x_0,\dots,x_{k-1}]$ --- this is where the flag hypothesis of Definition~\ref{def:simplicial-complex} enters, a common neighbor being automatically distinct from the $x_j$ since $\Gc$ has no loops.
Then $\Delta_{k}=\mathbf L_{k}^-+\mathbf L_{k}^+$ with
\[
  \mathbf L_{k}^-:= d^{k} \delta^{k},\qquad
  \mathbf L_{k}^+:= \delta^{k+1} d^{k+1}.
\]
For skew cochains (the default; the canonical-basis counterpart is discussed after the proof),
\begin{align*}
(\mathbf L_{k,\mathrm{skew}}^- f)([x_0,\dots,x_k])
&= \sum_{i=0}^{k}
   (-1)^{i}\,
   \frac{1}{m_{k-1}(x_0,\dots,\widehat{x_i},\dots,x_k)}
   \,
\\[-1mm]
&\hspace*{-2cm}\times
\sum_{z\in F_{(x_0,\dots,\widehat{x_i},\dots,x_k)}}
   m_k(z,x_0,\dots,\widehat{x_i},\dots,x_k)
   f([z,x_0,\dots,\widehat{x_i},\dots,x_k]),
\\[1mm]
(\mathbf L_{k,\mathrm{skew}}^+ f)([x_0,\dots,x_k])
&= \frac{1}{m_{k}(x_0,\dots,x_k)}
   \sum_{z\in F_{(x_0,\dots,x_k)}}\;
   m_{k+1}(z,x_0,\dots,x_k)\\
 &\hspace*{-1.5cm} \times \biggl(f([x_0,\dots,x_k])\;+\;\sum_{i=0}^{k} (-1)^{i+1}\,
   f([z,x_0,\dots,\widehat{x_i},\dots,x_k])\biggr),
\end{align*}
where in the second formula the $(k+2)$ terms of the alternating $d^{k+1}$-sum on $[z,x_0,\dots,x_k]$ have been split: the term corresponding to omitting $z$ (position $0$ in the augmented tuple, sign $+1$) gives back $f([x_0,\dots,x_k])$, while omitting $x_i$ (position $i+1$, sign $(-1)^{i+1}$) gives $f([z,x_0,\dots,\widehat{x_i},\dots,x_k])$.
\end{proposition}

\begin{proof}
Direct expansion using Lemma~\ref{lem:2-3-1} for $\delta^k$ and the definitions of $d^{k+1}, d^k$, $\delta^{k+1}$ from Section~\ref{sec:prelim-analytic}.
\end{proof}

\smallskip
\noindent\emph{The canonical-basis counterpart.}
The canonical-basis Laplacian $\Delta^{\mathrm{can}}_{k,\mathrm{sym}}$ of Definition~\ref{def:canonical-basis} is the transport of the above through the unitary $U_k$ of Theorem~\ref{thm:unitary-color}. We do not restate the formulas for it here. Since $U_k$ replaces cochains on oriented simplices by functions on unoriented ones, the transported formula carries, besides the signs $(-1)^i$ of the coboundary, the signs incurred in returning each term to its canonical representative. At top degree, where this local form is actually used, it is derived in the proof of Theorem~\ref{thm:kn-esa}; see~\eqref{eq:Delta-n-local-sym}.
\smallskip

\subsubsection{Friedrichs form and first representation theorem}

Define on $\Cc_{c}^\bullet(\Vc):=\Cc_c^0(\Vc)\oplus \Cc^1_{c}(\Vc)\oplus\cdots\oplus \Cc^n_{c}(\Vc)$ the (nonnegative) sesquilinear form
\[
  \mathfrak q(f,g):=\langle d f, d g\rangle_{\Hc}+\langle \delta f,\delta g\rangle_{\Hc},\qquad
  \mathfrak q(f):=\mathfrak q(f,f).
\]

\begin{proposition}[Green's identity and closability]\label{prop:2-5-1}
For all compactly supported $u,v$ in $\Cc_{c}^\bullet(\Vc)$ one has the Green identity
\[
\langle d u, v\rangle_{\Hc}=\langle u,\delta v\rangle_{\Hc},
\]
and the form $\mathfrak q$ is closable on $\Hc$.
\end{proposition}

\begin{proof}
The Green identity is immediate from Lemma~\ref{lem:2-3-1} and finite support. For closability, suppose $f_j\to 0$ in $\Hc$ and $\mathfrak q(f_j-f_\ell)\to 0$. Then $d f_j$ and $\delta f_j$ are Cauchy in $\Hc$ and converge to some $u,v$.
For any compactly supported $g$,
\(
\langle u,g\rangle=\lim_j\langle d f_j,g\rangle=\lim_j\langle f_j,\delta g\rangle=0,
\)
since $\delta g\in\Cc^\bullet_c(\Vc)\subset\Hc$ and $f_j\to0$; as $\Cc^\bullet_c(\Vc)$ is dense, $u=0$. The same argument applied to $\langle v,g\rangle=\lim_j\langle\delta f_j,g\rangle=\lim_j\langle f_j,dg\rangle$ gives $v=0$, this time using $dg\in\Hc$, which is where the standing assumption~\eqref{eq:pointwise-updeg} enters. Consequently $\mathfrak q(f_j)=\|df_j\|^2_{\Hc}+\|\delta f_j\|^2_{\Hc}\to0$. Hence $\mathfrak q$ is closable.
\end{proof}

We denote by $\mathfrak q$ also its closure, with domain $D[\mathfrak q]$.

\begin{proposition}\label{prop:friedrichs}
There exists a unique nonnegative self-adjoint operator $L^{\mathrm F}$ on $\Hc$ such that
\[
\mathfrak q(f,g)=\big\langle (L^{\mathrm F})^{1/2}f,\ (L^{\mathrm F})^{1/2}g\big\rangle_{\Hc},\qquad
D\!\big((L^{\mathrm F})^{1/2}\big)=D[\mathfrak q].
\]
Its decomposition into blocks is recorded in Corollary~\ref{cor:3-1-friedrichs}. If moreover the complex is locally finite, so that $L$ maps $\Cc^\bullet_c(\Vc)$ into $\Hc$ and is therefore a symmetric nonnegative operator on that core, then $L^{\mathrm F}$ is the Friedrichs extension of $L$. Without local finiteness $Lf$ need not lie in $\Hc$ for $f$ in the core --- $\delta^{k+1}d^{k+1}f$ has no reason to be square summable when $d^{k+1}$ is unbounded --- and $L^{\mathrm F}$ should then be read simply as the operator associated with the closed form $\mathfrak q$.
\end{proposition}

\begin{proof}
The form $\mathfrak q$ is nonnegative (clearly), densely defined (since $\Cc^\bullet_{c}(\Vc)$ is dense in $\Hc$ by Lemma~\ref{lem:density-Cc}), and closable (by Proposition~\ref{prop:2-5-1}). Its closure is therefore a closed, densely defined, nonnegative sesquilinear form on the complex Hilbert space $\Hc$. Kato's first representation theorem~\cite[Theorem~VI.2.1]{Kato} associates with such a form a unique nonnegative self-adjoint operator $L^{\mathrm F}$ satisfying
\[
\mathfrak q(f,g)=\big\langle (L^{\mathrm F})^{1/2}f,\ (L^{\mathrm F})^{1/2}g\big\rangle_{\Hc},\qquad
D\!\big((L^{\mathrm F})^{1/2}\big)=D[\mathfrak q].
\]
Under local finiteness, $Lf\in\Hc$ for $f$ in the core, and the Green identity of Proposition~\ref{prop:2-5-1} gives $\mathfrak q(f,f)=\langle Lf,f\rangle_{\Hc}$ there; that $L^{\mathrm F}$ is then the Friedrichs extension of $L$ follows from the construction in~\cite[\S VI.2.3]{Kato}.
\end{proof}

\begin{corollary}\label{cor:3-1}
For each degree $k\in\{0,\dots,n\}$, the sesquilinear form
\[
  \mathfrak q_{k}(u,v)
  := \langle d^{k+1} u, d^{k+1} v\rangle_{k+1} + \langle \delta^{k} u, \delta^{k} v\rangle_{k-1},
  \qquad u,v\in \Cc_{c}^k(\Vc),
\]
is closable on $\ell^2(m_k)$, with closure denoted $\mathfrak q_{k}$ and domain $D[\mathfrak q_{k}]$. (For $k=0$ the term $\langle\delta^{0}u,\delta^{0}v\rangle_{-1}$ is set to $0$ by the convention $\delta^{0}=0$; for $k=n$ the term $\langle d^{n+1}u,d^{n+1}v\rangle_{n+1}$ is set to $0$ by $d^{n+1}=0$).
\end{corollary}

\begin{proof}
Same argument as in Proposition~\ref{prop:2-5-1}, applied degree by degree; the standing assumption~\eqref{eq:pointwise-updeg} enters at the same point, namely to ensure $d^{k+1}g\in\ell^2(m_{k+1})$ for compactly supported $g$.
\end{proof}

\begin{corollary}\label{cor:3-1-friedrichs}
There exists a unique nonnegative self-adjoint operator
$\Delta_{k}^{\mathrm F}$ on $\ell^2(m_k)$ such that
\[
  \mathfrak q_{k}(u,v) = \big\langle (\Delta_{k}^{\mathrm F})^{1/2}u,\ (\Delta_{k}^{\mathrm F})^{1/2}v\big\rangle_k,
  \qquad
  D\!\big((\Delta_{k}^{\mathrm F})^{1/2}\big) = D[\mathfrak q_{k}].
\]
If the complex is locally finite, $\Delta_{k}^{\mathrm F}$ is moreover the Friedrichs extension of the block Laplacian $\Delta_{k}$ initially defined on $\Cc_{c}^k(\Vc)$; without that hypothesis $\Delta_ku$ need not lie in $\ell^2(m_k)$ for $u$ in the core, and $\Delta_k^{\mathrm F}$ is to be read as the operator associated with $\mathfrak q_k$. In all cases one has the orthogonal sum
\(
  L^{\mathrm F} = \bigoplus_{k=0}^n \Delta_{k}^{\mathrm F}.
\)
\end{corollary}

\begin{proof}
The form $\mathfrak q_{k}$ is closed, densely defined and nonnegative on $\ell^2(m_k)$ (Corollary~\ref{cor:3-1}), so Proposition~\ref{prop:friedrichs} applies verbatim in each degree and yields $\Delta_{k}^{\mathrm F}$ with the stated form representation, together with the identification as the Friedrichs extension of $\Delta_k$ under local finiteness. The orthogonal block sum follows from~\eqref{eq:2-4-1} and the block-diagonal structure of $\mathfrak q$ across degrees, the sum over $k$ being finite: the closure of a finite orthogonal sum of closable forms is the sum of the closures, so $D[\mathfrak q]=\bigoplus_k D[\mathfrak q_k]$ and the associated operators decompose accordingly.
\end{proof}

\subsection{Up-/down-Laplacians and Hodge identities}\label{sec:hodge}
\subsubsection{Energy identities}\label{subsubsec:energy}
Recall the up-/down-Laplacian decomposition from Proposition~\ref{rmk:2-4-1}: for each degree $k\in\{0,\dots,n\}$,
\[
  \mathbf L_{k}^- \ :=\ d^{k}\,\delta^{k},
  \qquad
  \mathbf L_{k}^+ \ :=\ \delta^{k+1}\,d^{k+1},
\]
(with the conventions $\delta^{0}=0$ and $d^{n+1}=0$), so that
\[
  \Delta_{k} \ =\ \mathbf L_{k}^- \ +\ \mathbf L_{k}^+,
  \qquad
  L \ =\ \bigoplus_{k=0}^{n}\Delta_{k}.
\]

\begin{lemma}\label{lem:2-6-1}
For $u\in\Cc_{c}^k(\Vc)$,
\[
  \langle \mathbf L_{k}^- u, u\rangle_k \ =\ \|\delta^{k}u\|_{k-1}^2,
  \qquad
  \langle \mathbf L_{k}^+ u, u\rangle_k \ =\ \|d^{k+1}u\|_{k+1}^2,
\]
and hence
\[
  \langle \Delta_{k} u, u\rangle_k \ =\ \|\delta^{k}u\|_{k-1}^2 \ +\ \|d^{k+1}u\|_{k+1}^2.
\]
Consequently, on the core $\Cc^k_c(\Vc)$ one has $\ker \Delta_{k}=\ker d^{k+1}\cap\ker \delta^{k}$; for the Friedrichs realization the corresponding statement is the one at the level of forms,
\[
\ker\Delta_{k}^{\mathrm F}\;=\;\{u\in D[\mathfrak q_k]:\ \mathfrak q_k[u]=0\}\;=\;D[\mathfrak q_k]\cap\ker\overline{d^{k+1}}\cap\ker\overline{\delta^{k}} .
\]
The intersection with $D[\mathfrak q_k]$ cannot be dropped in general, since only the inclusion $D[\mathfrak q_k]\subseteq D[\mathfrak q_k^-]\cap D[\mathfrak q_k^+]$ is available (Remark~\ref{rem:4-7-forms}); it is superfluous as soon as $\Cc^k_c(\Vc)$ is a form core for the sum $\mathfrak q_k^-+\mathfrak q_k^+$, that is as soon as $D[\mathfrak q_k]=D[\mathfrak q_k^-]\cap D[\mathfrak q_k^+]$ --- an equality that is sufficient, though not necessary: what the statement really needs is only $\ker\overline{d^{k+1}}\cap\ker\overline{\delta^{k}}\subseteq D[\mathfrak q_k]$. Being a core for each of the two forms separately does not suffice, and is in any case automatic. The equality does hold throughout the bounded regime of Section~\ref{sec:normalization}, and automatically in degree $0$, where $\mathfrak q_0=\mathfrak q_0^+$.
\end{lemma}

\begin{proof}
The first identity is Proposition~\ref{prop:2-5-1} applied to the pair $(\delta^ku,u)$, both compactly supported:
$
\langle d^{k}\delta^{k}u,u\rangle_k=\langle \delta^{k}u,\delta^{k}u\rangle_{k-1}.
$
The second requires a word, since $d^{k+1}u$ need not be compactly supported. Expanding both sides by Lemma~\ref{lem:2-3-1} gives the same double sum over pairs (a $k$-simplex in $\mathrm{supp}\,u$, an adjacent $(k+1)$-simplex), which converges absolutely: the inner sums are bounded by $m_k(\sigma)d_k(\sigma)<\infty$ under the standing assumption~\eqref{eq:pointwise-updeg}, and $\mathrm{supp}\,u$ is finite. Fubini then gives
$
\langle \delta^{k+1}d^{k+1}u,u\rangle_k=\langle d^{k+1}u,d^{k+1}u\rangle_{k+1}.
$
The last statement is immediate.
\end{proof}

\begin{remark}\label{rem:4-7-forms}
The form identity of Lemma~\ref{lem:2-6-1} extends to the Friedrichs realization: for all $u\in D[\mathfrak q_{k}]$,
\[
  \mathfrak q_{k}[u]\ =\ \bigl\|(\Delta_{k}^{\mathrm F})^{1/2}u\bigr\|_k^2 \ =\ \|\overline{\delta^{k}}u\|_{k-1}^2+\|\overline{d^{k+1}}u\|_{k+1}^2 ,
\]
the middle expression being the one that makes sense on all of $D[\mathfrak q_{k}]$: $\Delta_{k}^{\mathrm F}u$ is defined only on $D(\Delta_{k}^{\mathrm F})$, which is in general a proper subspace of $D[\mathfrak q_{k}]$.
More precisely, $\Delta_{k}^{\mathrm F}$ is, by definition (Corollary~\ref{cor:3-1-friedrichs}), the unique self-adjoint operator associated with the closed form $\mathfrak q_{k}$, that is, with the closure of the sum form
\[
  \mathfrak q_{k}^{-}+\mathfrak q_{k}^{+}\quad\text{restricted to the core } \Cc_{c}^k(\Vc),
  \qquad\text{where } \mathfrak q_{k}^{-}[u]:=\|\overline{\delta^{k}}u\|_{k-1}^2,\quad \mathfrak q_{k}^{+}[u]:=\|\overline{d^{k+1}}u\|_{k+1}^2
\]
are the closed forms associated with $(\mathbf L_{k}^-)^{\mathrm F}$ and $(\mathbf L_{k}^+)^{\mathrm F}$ respectively.

The closures $\overline{\delta^{k}}$ and $\overline{d^{k+1}}$ exist: each of the two operators, taken on the core, has a densely defined adjoint --- containing $d^k$, respectively $\delta^{k+1}$, on the core by Green's identity --- and is therefore closable. We stress that only the inclusion
\[
  D[\mathfrak q_{k}]\ \subseteq\ D[\mathfrak q_{k}^{-}]\cap D[\mathfrak q_{k}^{+}]
\]
holds in general: $\mathfrak q_{k}$ is the closure of a restriction of the sum form to the core, and equality of the two form domains amounts to the (nontrivial) statement that $\Cc_{c}^k(\Vc)$ is a form core for $\mathfrak q_{k}^{-}+\mathfrak q_{k}^{+}$; we do not use such an equality anywhere below.

We emphasize that, at this stage of the paper -- before boundedness has been established -- one should \emph{not} assert the operator identity $\Delta_{k}^{\mathrm F}=(\mathbf L_{k}^-)^{\mathrm F}+(\mathbf L_{k}^+)^{\mathrm F}$ as an unqualified statement. The sum of two unbounded self-adjoint operators is in general only symmetric, and only on the intersection $D\bigl((\mathbf L_{k}^-)^{\mathrm F}\bigr)\cap D\bigl((\mathbf L_{k}^+)^{\mathrm F}\bigr)$ of their operator domains, which may be strictly smaller. It may even fail to be densely defined, let alone self-adjoint or equal to the operator associated with the closed sum form $\mathfrak q_{k}$. What does hold unconditionally is the operator identity on that intersection,
\[
\Delta_{k}^{\mathrm F}u \ =\ (\mathbf L_{k}^-)^{\mathrm F}u+(\mathbf L_{k}^+)^{\mathrm F}u,\qquad u\in D\bigl((\mathbf L_{k}^-)^{\mathrm F}\bigr)\cap D\bigl((\mathbf L_{k}^+)^{\mathrm F}\bigr),
\]
and this identity should not be asserted beyond this intersection. Two hypotheses remove the difficulty: the finiteness $D_{k}^{\downarrow}+D_{k}^{\uparrow}<\infty$ of Lemma~\ref{lem:form-bounds} below, and uniform ellipticity of the weights, under which $\Delta_{k}$ is similar to the normalized block $\widetilde\Delta_{k}$ (Lemma~\ref{lem:elliptic-similarity}). Under both, all three operators are everywhere defined and bounded. The domain issue then disappears, and the operator sum does coincide with $\Delta_{k}^{\mathrm F}$ on all of $\ell^2(m_k)$. No statement below relies on the operator identity outside that regime.
\end{remark}

\subsection{Normalization conventions and spectral footprint}\label{sec:normalization} 
The distinction between energy (or Friedrichs) and degree-normalized Laplacians has a long history in discrete Hodge theory. In the finite setting, \cite{HJ} treats weighted simplicial complexes and clarifies how normalization choices affect spectral properties. Our normalization conventions extend those of \cite{HJ} and fit within the framework of \cite{BaKe}; our operator bounds provide geometry-free criteria ensuring that the hypotheses of that abstract theory are satisfied.

\subsubsection{Two parallel normalizations}
(1) \emph{Energy/Friedrichs normalization (default).} We work in $\ell^2(m_k)$ and set
$\Delta_{k}=\mathbf L_{k}^-+\mathbf L_{k}^+$, so that
\(
\langle \Delta_{k}u,u\rangle_k=\|\delta^{k}u\|_{k-1}^2+\|d^{k+1}u\|_{k+1}^2
\)
(Lemma~\ref{lem:2-6-1}). This yields closed nonnegative forms and the Friedrichs realization.

\noindent(2) \emph{Symmetric (degree) normalization.} On the unweighted $\ell^2$ we define
\begin{equation}\label{eq:norm-blocks}
\mathbb L_{k}^- \ :=\ M_k^{1/2}\,\mathbf L_{k}^-\,M_k^{-1/2},\
\mathbb L_{k}^+ \ :=\ M_k^{1/2}\,\mathbf L_{k}^+\,M_k^{-1/2},\
\widetilde{\Delta}_{k}\ :=\ \mathbb L_{k}^-+\mathbb L_{k}^+.
\end{equation}
Here $M_k$ denotes multiplication by $m_k$. Note that $M_k^{-1/2}$ is \emph{unitary} from the unweighted $\ell^2$ onto $\ell^2(m_k)$, so that~\eqref{eq:norm-blocks} is the transport of the energy blocks through that unitary: $\mathbb L_k^{\pm}$ are again symmetric on the core, and again come from closed nonnegative forms. Being a unitary transport, $\widetilde\Delta_{k}$ is isospectral to $\Delta_{k}^{\mathrm F}$ and has the same norm, with no assumption on the weights (Lemma~\ref{lem:elliptic-similarity}).

\paragraph{Form bounds with explicit constants.}
Define the \emph{down/up degrees} (with the multiplicities $(k+1),(k+2)$ absorbed; cf.\ the convention~\eqref{eq:Dk-conv} below)
\begin{equation}\label{eq:deg-du}
\begin{split}
&D_{k}^{\downarrow}\ :=\ (k+1)\sup_{|\tau|}\ \frac{1}{m_{k-1}(|\tau|)}\sum_{z\in F_{|\tau|}} m_k(|\tau|,z)
\qquad(1\le k\le n),\\
&D_{k}^{\uparrow}\ :=\ (k+2)\sup_{|\sigma|}\ \frac{1}{m_{k}(|\sigma|)}\sum_{z\in F_{|\sigma|}} m_{k+1}(|\sigma|,z)
\qquad(0\le k\le n-1),
\end{split}
\end{equation}
where $|\tau|$ runs over the unoriented $(k-1)$-simplices and $|\sigma|$ over the unoriented $k$-simplices of $\Sc_n$. At the two extreme degrees the quantities above involve $m_{-1}$, resp.\ $m_{n+1}$, which do not exist; consistently with the conventions $\delta^{0}=0$ and $d^{n+1}=0$ we set
\begin{equation}\label{eq:deg-extremes}
D_{0}^{\downarrow}:=0,\qquad D_{n}^{\uparrow}:=0 ,
\end{equation}
so that $\mathbb L_0^-=0$ and $\mathbb L_n^+=0$ and the corresponding bounds in~\eqref{eq:form-bounds} are trivially satisfied. Here $m_{k-1}(|\tau|), m_k(|\sigma|)$ are well-defined by the symmetry of the weights (Definition~\ref{def:weights}), and $F_{|\tau|}$ denotes the common-neighbor set of any (hence every) ordered representative of $|\tau|$. It is convenient to name separately the two \emph{raw} degrees, without the multiplicities:
\begin{equation}\label{eq:raw-degrees}
D_k^{\downarrow,\mathrm{HJ}}:=\sup_{|\tau|}\frac{1}{m_{k-1}(|\tau|)}\sum_{z\in F_{|\tau|}} m_k(|\tau|,z),
\qquad
D_k^{\uparrow,\mathrm{HJ}}:=\sup_{|\sigma|}\frac{1}{m_{k}(|\sigma|)}\sum_{z\in F_{|\sigma|}} m_{k+1}(|\sigma|,z),
\end{equation}
so that
\begin{equation}\label{eq:Dk-conv}
D_k^\downarrow\;=\;(k+1)\, D_k^{\downarrow,\mathrm{HJ}},\qquad
D_k^\uparrow\;=\;(k+2)\, D_k^{\uparrow,\mathrm{HJ}}.
\end{equation}

\paragraph{Comparison with the Horak--Jost normalization.} Throughout this paper we absorb the multiplicities $(k+1)$ and $(k+2)$ --- which arise from counting the ordered representatives of unoriented $k$- and $(k+1)$-simplices --- into the symbols $D_k^\downarrow$ and $D_k^\uparrow$; the raw quantities are~\eqref{eq:raw-degrees}. The quantity $D_k^{\uparrow,\mathrm{HJ}}$ is the up-degree in the convention of~\cite[Definition~3.1]{HJ} (where it is simply called \emph{degree}), so a reader comparing with~\cite{HJ} should multiply the Horak--Jost degree by $(k+2)$ to obtain our $D_k^\uparrow$. The dual quantity is not a new one: comparing the two definitions in~\eqref{eq:raw-degrees} gives $D_k^{\downarrow,\mathrm{HJ}}=D_{k-1}^{\uparrow,\mathrm{HJ}}$, and hence $D_k^{\downarrow}=D_{k-1}^{\uparrow}$, so the down-degree in degree $k$ \emph{is} the Horak--Jost degree read one degree lower (Remark~\ref{rem:degrees-coincide}). The multiplicity $(k+1)$ arises from the bijection $P_k\leftrightarrow\Fc_k$ used in the proof below; the multiplicity $(k+2)$ arises from the analogous bijection $P_{k+1}\leftrightarrow\Fc_{k+1}$ together with the Cauchy--Schwarz step on the $k+2$ terms of $d^{k+1}$.

The comparison is sharpest in the normalized case. When the weights are chosen so that $D_k^{\uparrow,\mathrm{HJ}}=1$ --- the higher-degree analogue of the vertex weight $m_0(x)=\sum_y m_1(x,y)$ of the Consequences below --- the up part of~\eqref{eq:form-bounds} reduces to $\|\mathbb L_k^{+}\|\le k+2$. This is exactly the bound of~\cite{HJ} for the normalized up-Laplacian. Its equality case, and the refinements of~\cite{FWW,SWF} through the balancedness of the incidence-signed graph, are discussed in Section~\ref{sec:related}.

\begin{remark}[The two degree families coincide]\label{rem:degrees-coincide}
Directly from~\eqref{eq:deg-du}, the two suprema in $D_k^{\downarrow}$ and $D_{k-1}^{\uparrow}$ are taken over the same set of $(k-1)$-simplices, of the same quantity, with the same multiplicity $(k+1)=((k-1)+2)$. Hence
\[
D_k^{\downarrow}\;=\;D_{k-1}^{\uparrow}\qquad(1\le k\le n),
\]
so~\eqref{eq:form-bounds} is a function of the single sequence $(D_k^{\uparrow})_{0\le k\le n-1}$: it reads $\|\widetilde\Delta_k\|\le D_{k-1}^{\uparrow}+D_k^{\uparrow}$. The coincidence is not an accident of the definitions. Writing $A:=d^k$, one has $\mathbf L_k^-=AA^\ast$ and $\mathbf L_{k-1}^+=A^\ast A$, so these two blocks have the same non-zero spectrum and the same norm; the two estimates of Lemma~\ref{lem:form-bounds} bound the same quantity and produce the same constant. At $k=1$ this is Lemma~\ref{lem:iso-01}, used in Section~\ref{sharp}; the general case is the same argument. In particular the top-degree block $\Delta_n=\mathbf L_n^-$ is spectrally the up-Laplacian $\mathbf L_{n-1}^+$ in degree $n-1$, which is the operator estimated in the finite setting by~\cite{HJ,FWW,SWF}; Theorem~\ref{thm:kn-esa} should be read as the $\ell^2$, infinite, weighted counterpart of those bounds.
\end{remark}

\begin{lemma}[Schur-type form bound]\label{lem:form-bounds}
Let $\Sc_n$ be a weighted $n$-simplicial complex as in Definition~\ref{def:simplicial-complex}, let $k\in\{0,\dots,n\}$, and let $D_{k}^{\downarrow},D_{k}^{\uparrow}\in[0,\infty]$ be the down/up degrees defined in~\eqref{eq:deg-du}. Then, for all $u\in\Cc_{c}^k(\Vc)$,
\begin{equation}\label{eq:form-bounds}
\langle \mathbb L_{k}^- u, u\rangle \ \le\ D_{k}^{\downarrow}\,\|u\|^2,\quad
\langle \mathbb L_{k}^+ u, u\rangle \ \le\ D_{k}^{\uparrow}\,\|u\|^2,\quad
\langle \widetilde{\Delta}_{k} u, u\rangle \ \le\ (D_{k}^{\downarrow}+D_{k}^{\uparrow})\,\|u\|^2.
\end{equation}
Consequently, if $D_{k}^{\downarrow}+D_{k}^{\uparrow}<\infty$, then $\widetilde\Delta_{k}$ is bounded on the core $\Cc_{c}^k(\Vc)$ and extends uniquely to a bounded self-adjoint operator on $\ell^2$ with
\[
\|\widetilde\Delta_{k}\|\ \le\ D_{k}^{\downarrow}+D_{k}^{\uparrow}.
\]
\end{lemma}

\begin{remark}[On the quality of the constant]\label{rem:form-bounds-sharpness}
No sharpness is claimed in~\eqref{eq:form-bounds} for general $k$. In degree $0$ the bound is attained on a non-trivial class of graphs: with $m_0\equiv m_1\equiv1$ on a $d$-regular graph one has $D_0^{\downarrow}+D_0^{\uparrow}=2d$, and $\|\widetilde\Delta_0\|=2d$ precisely when $\inf\mathrm{spec}(Q_G)=0$, which holds in particular on the bipartite amenable ones (Proposition~\ref{prop:sharpness-bipartite}). It is not attained on all $d$-regular graphs: on the $d$-regular tree with $d\ge3$ one has $\|\widetilde\Delta_0\|=d+2\sqrt{d-1}<2d$ (Remark~\ref{rem:tree-counterexample}). In intermediate degrees, by contrast, the constant should be read as a boundedness criterion with an explicit --- not an optimal --- constant.

Its slack is not localized in the Cauchy--Schwarz step on the $k+2$ terms of $d^{k+1}$ alone. Remark~\ref{rem:degrees-coincide} shows that this step produces the same number as the weighted Cauchy--Schwarz of the down estimate one degree higher; this is an observation, not a proof that the step is lossless, since the two may lose identically. The slack comes from replacing, in $D_k^{\downarrow}=(k+1)D_k^{\downarrow,\mathrm{HJ}}$, the sum over the $k+1$ faces of a given simplex by $(k+1)$ times a single supremum over all faces, and --- in the weighted case --- from the comparison between a supremum and a geometric mean; see the end of Section~\ref{sec:esa}. At top degree, Theorem~\ref{thm:kn-esa} recasts the same constant in a structurally more informative form; the comparison of the two routes is carried out at the end of Section~\ref{sec:esa}.
\end{remark}

Before proving Lemma~\ref{lem:form-bounds}, we record an enumeration identity that is the main combinatorial step in the argument.

\begin{lemma}\label{lem:bijection-Fk}
Define the set
\[
P_k\;:=\;\bigl\{((x_1,\dots,x_k),z):\ (x_1,\dots,x_k)\in\Fc_{k-1},\ z\in F_{(x_1,\dots,x_k)}\bigr\}.
\]
Then, for each fixed insertion position $j\in\{0,\dots,k\}$, the map $\Phi_j:P_k\to\Fc_k$ defined by
\[
\Phi_j\bigl((x_1,\dots,x_k),z\bigr):=(x_1,\dots,x_j,z,x_{j+1},\dots,x_k)
\]
is a bijection. We write $\Phi:=\Phi_k$ for the map that appends $z$ in last position.
\end{lemma}

\begin{proof}
$\Phi_j$ is well-defined, since $z$ is adjacent to every $x_i$, so the resulting tuple is again a clique; and it is injective by construction, the position $j$ being fixed. For surjectivity, given $(y_0,\dots,y_k)\in\Fc_k$, put $z:=y_j$ and $(x_1,\dots,x_k):=(y_0,\dots,\widehat{y_j},\dots,y_k)$. Any sub-tuple of a clique is a clique, so $(x_1,\dots,x_k)\in\Fc_{k-1}$, and $y_j$ is adjacent to all its entries, so $z\in F_{(x_1,\dots,x_k)}$; by construction $\Phi_j((x_1,\dots,x_k),z)=(y_0,\dots,y_k)$, and this preimage is clearly the only one.
\end{proof}

\begin{proof}[Proof of Lemma~\ref{lem:form-bounds}]
The two extreme degrees are trivial and are set aside first: for $k=0$ one has $\delta^{0}=0$, hence $\mathbb L_0^-=0$, and the first bound holds with $D_0^{\downarrow}=0$ as in~\eqref{eq:deg-extremes}; for $k=n$ one has $d^{n+1}=0$, hence $\mathbb L_n^+=0$ and the second bound holds with $D_n^{\uparrow}=0$. In what follows we therefore assume $1\le k\le n$ for the down estimate and $0\le k\le n-1$ for the up estimate, which are exactly the ranges on which the corresponding degrees in~\eqref{eq:deg-du} are defined.

We first record the fundamental identity that links the normalized blocks (acting on the unweighted $\ell^2$) to weighted norms of cochains. By definition,
\[
\mathbb L_{k}^- \;=\; M_k^{1/2}\,d^{k}\delta^{k}\,M_k^{-1/2}\quad\text{on the unweighted }\ell^2.
\]
Setting $v:=M_k^{-1/2}u$ (so that $u\in\ell^2$ unweighted iff $v\in\ell^2(m_k)$, and $\|v\|_{\ell^2(m_k)}=\|u\|_{\ell^2}$), we compute
\begin{equation}\label{eq:Lminus-norm-id}
\langle \mathbb L_{k}^- u, u\rangle_{\ell^2}
\;=\;\langle d^{k}\delta^{k}\,M_k^{-1/2}u,\,M_k^{-1/2}u\rangle_{\ell^2(m_k)}
\;=\;\|\delta^{k}v\|_{\ell^2(m_{k-1})}^2.
\end{equation}
Indeed, the bracket $\langle\cdot,\cdot\rangle_{\ell^2(m_k)}$ already absorbs the $M_k$ factor in the definition of the weighted inner product~\eqref{eq:2-2-1}, and $\langle d^k\delta^kv,v\rangle_{\ell^2(m_k)}=\|\delta^kv\|_{\ell^2(m_{k-1})}^2$ is Lemma~\ref{lem:2-6-1}.

By Lemma~\ref{lem:2-3-1},
$$
(\delta^{k} v)([x_1,\dots,x_k])=\frac{1}{m_{k-1}(x_1,\dots,x_k)}\sum_{z\in F_{(x_1,\dots,x_k)}} m_k(z,x_1,\dots,x_k)\, v([z,x_1,\dots,x_k]).
$$

By Lemma~\ref{lem:bijection-Fk}, summing over $((x_1,\dots,x_k),z)\in P_k$ is the same as summing over $(x_0,\dots,x_k)\in\Fc_k$:
\begin{equation}\label{eq:sum-bijection}
\sum_{(x_1,\dots,x_k)\in\Fc_{k-1}}\sum_{z\in F_{(x_1,\dots,x_k)}}h(x_1,\dots,x_k,z)\;=\;\sum_{(x_0,\dots,x_k)\in\Fc_k}h(x_0,\dots,x_k),
\end{equation}
for any function $h$ on $\Fc_k$.

\emph{Main estimate (down-Laplacian).}
The weighted inner product on $\ell^2(m_{k-1})$ is, by~\eqref{eq:2-2-1}, the sum over ordered $k$-tuples $(x_1,\dots,x_k)\in\Fc_{k-1}$ divided by $k!$, and similarly for $\ell^2(m_k)$ with factor $1/(k+1)!$ on $\Fc_k$. Abbreviate $\vec x:=(x_1,\dots,x_k)$ and $F_{\vec x}:=F_{(x_1,\dots,x_k)}$. We use $m_k(z,\vec x)=m_k(\vec x,z)$ by symmetry of the weights, and $|v([z,\vec x])|=|v([\vec x,z])|$ because $v$ is alternating, so that the sign disappears in the modulus. Then
\begin{align*}
\|\delta^{k} v\|_{k-1}^2
&= \frac{1}{k!}\sum_{\vec x\in\Fc_{k-1}} \frac{1}{m_{k-1}(\vec x)}\,\biggl|\sum_{z\in F_{\vec x}} m_k(\vec x,z)\, v([\vec x,z])\biggr|^2\\
& \le\ \frac{1}{k!}\sum_{\vec x\in\Fc_{k-1}} \frac{1}{m_{k-1}(\vec x)}\biggl(\sum_{z\in F_{\vec x}} m_k(\vec x,z)\biggr)\biggl(\sum_{z\in F_{\vec x}} m_k(\vec x,z)\,\bigl|v([\vec x,z])\bigr|^2\biggr)\\
&\le\ D_{k}^{\downarrow,\mathrm{HJ}}\,\frac{1}{k!}\sum_{(\vec x,z)\in P_k} m_k(\vec x,z)\,\bigl|v([\vec x,z])\bigr|^2,
\end{align*}
the last step using $\frac{1}{m_{k-1}(\vec x)}\sum_{z\in F_{\vec x}} m_k(\vec x,z)\le D_{k}^{\downarrow,\mathrm{HJ}}$ from~\eqref{eq:raw-degrees}.
Each unoriented $k$-simplex has exactly $(k+1)!$ ordered representatives in $\Fc_k$, so by~\eqref{eq:sum-bijection}, the right-hand side equals $D_k^{\downarrow,\mathrm{HJ}}\cdot(k+1)\cdot\|v\|_{\ell^2(m_k)}^2$ (the factor $(k+1)=(k+1)!/k!$ being the ratio between the normalizations of $\ell^2(m_k)$ and the $\Fc_k$-sum). Hence
\[
\|\delta^k v\|_{k-1}^2\;\le\;(k+1)\,D_k^{\downarrow,\mathrm{HJ}}\,\|v\|_{\ell^2(m_k)}^2\;=\;D_k^\downarrow\,\|v\|_{\ell^2(m_k)}^2,
\]
using~\eqref{eq:Dk-conv}. Combined with~\eqref{eq:Lminus-norm-id}, this proves the first bound in~\eqref{eq:form-bounds}.

\smallskip
\emph{Main estimate (up-Laplacian).} The argument for $\langle \mathbb L_{k}^+ u,u\rangle\le D_k^\uparrow\|u\|^2$ is structurally parallel to the down case but uses the dual bijection $P_{k+1}\leftrightarrow\Fc_{k+1}$, where
\[
P_{k+1}\;:=\;\bigl\{((y_0,\dots,y_k),w):\ (y_0,\dots,y_k)\in\Fc_k,\ w\in F_{(y_0,\dots,y_k)}\bigr\}.
\]
Starting from $\langle\mathbb L_{k}^+ u,u\rangle_{\ell^2}=\|d^{k+1} v\|_{\ell^2(m_{k+1})}^2$ (the analogue of~\eqref{eq:Lminus-norm-id} for the up block) with $v:=M_k^{-1/2}u$, we estimate $|(d^{k+1} v)(y_0,\dots,y_k,w)|^2$. By the formula for $d^{k+1}$ (Section~\ref{sec:prelim-comb}), this is the squared modulus of an alternating sum of $k+2$ terms, each of the form $\pm v$ evaluated on a codimension-one face of the $(k+1)$-tuple. By Cauchy--Schwarz (in the unweighted form, on $k+2$ terms), 
\[
\bigl|(d^{k+1} v)(y_0,\dots,y_k,w)\bigr|^2 \;\le\; (k+2)\sum_{j=0}^{k+1}\bigl|v(\widehat{j})\bigr|^2,
\]
where $\widehat{j}$ denotes the $j$-th codimension-one face of $(y_0,\dots,y_k,w)$. Summing over $(y_0,\dots,y_k,w)\in P_{k+1}$ with the weight $m_{k+1}(y_0,\dots,y_k,w)/(k+2)!$ (the normalization of $\ell^2(m_{k+1})$, by~\eqref{eq:2-2-1}) and using the bijection $P_{k+1}\leftrightarrow\Fc_{k+1}$, the right-hand side becomes
\begin{align*}
&\frac{k+2}{(k+2)!}\sum_{(z_0,\dots,z_{k+1})\in\Fc_{k+1}}
m_{k+1}(z_0,\dots,z_{k+1})
\sum_{j=0}^{k+1}\bigl|v(\widehat{j})\bigr|^2
\\
&=\;
\frac{(k+2)^2}{(k+2)!}
\sum_{(y_0,\dots,y_k)\in\Fc_k}
\Bigl(
\sum_{w\in F_{(y_0,\dots,y_k)}}
m_{k+1}(y_0,\dots,y_k,w)
\Bigr)
|v(y_0,\dots,y_k)|^2,
\end{align*}
The last step exchanges the sums on $j$ and on $\Fc_{k+1}$. Fix $j\in\{0,\dots,k+1\}$ and delete the $j$-th entry: the map $(z_0,\dots,z_{k+1})\mapsto\bigl((z_0,\dots,\widehat{z_j},\dots,z_{k+1}),z_j\bigr)$ identifies the sum over $\Fc_{k+1}$ with the sum over $P_{k+1}$ --- it is $\Phi_j^{-1}$ of Lemma~\ref{lem:bijection-Fk}, applied in degree $k+1$. Since $m_{k+1}$ is symmetric in its arguments, the resulting contribution does not depend on $j$, and summing over the $k+2$ values of $j$ produces the factor $k+2$ in the numerator. By the definition~\eqref{eq:raw-degrees} of the raw up-degree $D_k^{\uparrow,\mathrm{HJ}}$,
the inner sum is bounded by $D_k^{\uparrow,\mathrm{HJ}}\cdot m_k(y_0,\dots,y_k)$. Since $\|v\|^2_{\ell^2(m_k)}=\frac{1}{(k+1)!}\sum_{(y_0,\dots,y_k)\in\Fc_k}m_k(y_0,\dots,y_k)|v(y_0,\dots,y_k)|^2$, the prefactor of the last display becomes $\frac{(k+2)^2}{(k+2)!}\cdot(k+1)!=\frac{(k+2)^2}{k+2}=k+2$, and one obtains
\[
\|d^{k+1} v\|^2_{\ell^2(m_{k+1})}\le (k+2)\,D_k^{\uparrow,\mathrm{HJ}}\,\|v\|_{\ell^2(m_k)}^2=D_k^\uparrow\|v\|_{\ell^2(m_k)}^2,
\]
which is the second bound in~\eqref{eq:form-bounds}. The third bound follows by adding the down and up bounds. Finally, assume $D_k^\downarrow+D_k^\uparrow<\infty$ and set $C:=D_k^\downarrow+D_k^\uparrow$. Consider the sesquilinear form $\mathfrak t(u,w):=\langle\widetilde\Delta_{k}u,w\rangle$ for $u,w$ in the core $\Cc_{c}^k(\Vc)$. It is well defined: $\widetilde\Delta_{k}u$ is a genuine cochain by the standing assumption~\eqref{eq:pointwise-updeg}, and pairing it with the compactly supported $w$ is a finite sum, so no integrability question arises. It is Hermitian and nonnegative by Green's identity and Lemma~\ref{lem:2-6-1}. The third bound states that $\mathfrak t(u,u)\le C\|u\|^2$. By the generalized Cauchy--Schwarz inequality for the nonnegative form $\mathfrak t$,
\[
|\langle\widetilde\Delta_{k}u,w\rangle|\;=\;|\mathfrak t(u,w)|\;\le\;\mathfrak t(u,u)^{1/2}\,\mathfrak t(w,w)^{1/2}\;\le\;C\,\|u\|\,\|w\|
\qquad(u,w\in\Cc^k_c(\Vc)).
\]
Since $\langle\cdot,\cdot\rangle$ is conjugate-linear in its second argument, $w\mapsto\overline{\mathfrak t(u,w)}$ is a bounded \emph{linear} functional on the dense subspace $\Cc^k_c(\Vc)$, of norm $\le C\|u\|$. By the Riesz representation theorem it is represented by a unique vector $Au\in\ell^2$ with $\|Au\|\le C\|u\|$. Moreover $Au=\widetilde\Delta_{k}u$, and in particular $\widetilde\Delta_{k}u\in\ell^2$: for every $w\in\Cc^k_c(\Vc)$ one has $\langle Au,w\rangle=\mathfrak t(u,w)=\langle\widetilde\Delta_{k}u,w\rangle$, both pairings being finite sums, so the cochain $\widetilde\Delta_{k}u$ and the vector $Au$ agree at every $k$-simplex, as is seen by taking for $w$ the indicator of that simplex. The map $u\mapsto Au$ is linear, densely defined, bounded by $C$ and symmetric, hence extends by continuity to a bounded self-adjoint operator on $\ell^2$ with $\|\widetilde\Delta_{k}\|\le D_k^\downarrow+D_k^\uparrow$. This proves the stated consequence.

A word on why the argument is phrased through $\mathfrak t$. The shorter route would be to write
\[
\|\widetilde\Delta_{k}u\|^2\ \le\ \mathfrak t(u,u)^{1/2}\,\mathfrak t\bigl(\widetilde\Delta_ku,\widetilde\Delta_ku\bigr)^{1/2},
\]
but this presupposes that $\widetilde\Delta_{k}$ maps the core into itself, which requires local finiteness; that hypothesis is not made here.
\end{proof}

\begin{lemma}\label{lem:elliptic-similarity}
Assume $0<c_k\le m_k(\sigma)\le C_k<\infty$ for all $k$-simplices $\sigma$.
Then $M_k^{\pm1/2}$ are boundedly invertible on $\ell^2$, and
\(
\widetilde{\Delta}_{k}=M_k^{1/2}\,\Delta_{k}^{\mathrm F}\,M_k^{-1/2}
\)
holds as an identity of self-adjoint operators with
\(
\|M_k^{1/2}\|\le C_k^{1/2}\) and \(\|M_k^{-1/2}\|\le c_k^{-1/2}.
\)
In particular $\mathrm{spec}(\widetilde{\Delta}_{k})=\mathrm{spec}(\Delta_{k}^{\mathrm F})$ and $e^{-t\widetilde{ \Delta}_{k}}=M_k^{1/2}e^{-t\Delta_{k}^{\mathrm F}}M_k^{-1/2}$. Let us stress that the spectral identity, and hence the equality of norms, in fact needs no ellipticity at all: by the remark following~\eqref{eq:norm-blocks}, $M_k^{-1/2}$ is a unitary from the unweighted $\ell^2$ onto $\ell^2(m_k)$ and $\widetilde\Delta_k$ is the transport of $\Delta_k^{\mathrm F}$ through it. Uniform ellipticity is what allows the conjugation to be read \emph{within} a single fixed space, with the quantitative bounds $\|M_k^{\pm1/2}\|$ above.
\end{lemma}

\begin{proof}
Bounded invertibility is immediate from the inequalities $$c_k\|u\|^2\le \langle M_k u,u\rangle\le C_k\|u\|^2.$$
Let $\mathfrak q_{k}[u]=\|d^{k+1}u\|_{k+1}^2+\|\delta^{k}u\|_{k-1}^2$ on $\ell^2(m_k)$ with core $\Cc_{c}^k$; by closability its closure defines the Friedrichs realization $\Delta_{k}^{\mathrm F}$ (Corollary~\ref{cor:3-1-friedrichs}).
Define on unweighted $\ell^2$ the form
\(
\widehat{\mathfrak q}_{k}[w]:=\mathfrak q_{k}[M_k^{-1/2}w].
\)
Then $\widehat{\mathfrak q}_{k}$ is closed, nonnegative, and its operator is precisely $\widetilde{\Delta}_{k}$ by \eqref{eq:norm-blocks}.
Since $M_k^{\pm1/2}$ are boundedly invertible, $w\mapsto M_k^{-1/2}w$ is a topological isomorphism between the form domains, and Kato's first representation theorem~\cite[Theorem~VI.2.1]{Kato} yields
$\widetilde{\Delta}_{k}=M_k^{1/2}\Delta_{k}^{\mathrm F}M_k^{-1/2}$ with the asserted semigroup and spectral relations.
\end{proof}

\subsubsection{Consequences}
By Lemma~\ref{lem:form-bounds}, whenever $D_{k}^{\downarrow}+D_{k}^{\uparrow}<\infty$ the operator $\widetilde{\Delta}_{k}$ is bounded, self-adjoint, and extends uniquely from the core to all of $\ell^2$ with $\|\widetilde\Delta_{k}\|\le D_{k}^{\downarrow}+D_{k}^{\uparrow}$; essential self-adjointness on the core follows immediately from this boundedness.
Take $k=0$ and the normalizing vertex weight $m_0(x)=\sum_{y}m_1(x,y)$, for which $d_{\Vc}\equiv1$ and hence $D_0^{\uparrow}=2$. One recovers the classical normalized graph Laplacian, with $\mathrm{spec}\subset[0,2]$. In the weighted setting this standard bound follows from the form representation $0\le \langle \widetilde\Delta_{0}u,u\rangle\le 2\|u\|^2$ via Cauchy--Schwarz on the coboundary; see, e.g., \cite[Lemma~1.7(v)]{Chung} or \cite{KLW-monograph}. For a general vertex weight the corresponding bound is $\|\widetilde\Delta_0\|\le D_0^{\uparrow}=2\sup_{x}d_{\Vc}(x)$, which need not be $\le 2$.
At top degree $k{=}n$, $\widetilde{\Delta}_{n}$ reduces to a discrete Schr\"odinger operator on the line-complex: an adjacency operator with a \emph{signed} kernel, plus a nonnegative diagonal potential $V_n$ carrying the diagonal contribution (Theorem~\ref{thm:kn-esa} in Section~\ref{sec:esa}, and Section~\ref{sharp}). It is \emph{not} a graph Laplacian on the line-complex: the sign in front of the adjacency term is $+$, not $-$, as~\eqref{eq:Delta-n-local-skew} makes explicit.

\section{Unitary equivalence on colorable complexes}\label{sec:color}

\begin{definition}[Proper colorings and ordered colorings]\label{def:coloring}
A graph (or the $1$-skeleton of a simplicial complex) is called \emph{$p$-partite} if its vertex set can be written as a disjoint union
\[
\mathcal{V} = \Vc_1 \sqcup \Vc_2 \sqcup \cdots \sqcup \Vc_p,
\]
such that no two vertices in the same subset $\Vc_i$ are adjacent. Equivalently, $G$ admits a \emph{proper vertex-coloring} with $p$ colors $c:\Vc\to\{1,\dots,p\}$, in which case the chromatic number satisfies $\chi(G)\le p$; for $p=2$ one recovers bipartite graphs.

More generally, by an \emph{ordered coloring} of $G$ we mean a proper vertex-coloring $c\colon\Vc\to C$ with values in a totally ordered set $(C,\le)$. We call $G$ (and the complex $\Sc_n$) \emph{colorable} if it admits an ordered coloring. This is the only structure used in this section. The construction below requires two things: that $c$ be injective on each simplex, which properness guarantees since the vertices of a simplex are pairwise adjacent; and that the colors be comparable, so that each simplex has a well-defined color-increasing ordering.

\emph{The hypothesis is not restrictive.} Every countable graph admits an ordered coloring with $C=\N$. Indeed, enumerate $\Vc=\{v_1,v_2,\dots\}$ and color $v_j$ with the least element of $\N$ not used by the neighbors of $v_j$ among $v_1,\dots,v_{j-1}$ --- a finite set, whatever the degree of $v_j$. (The index $j$ is unrelated to the dimension $n$ of the complex.) Since our vertex sets are countable (Section~\ref{sec:prelim-comb}), the colorability hypothesis of Theorem~\ref{thm:unitary-color} below is satisfied by \emph{every} complex considered in this paper.

We nevertheless state it explicitly, for two reasons. The construction depends on a \emph{fixed} choice of ordered coloring. And the quantitative examples of Section~\ref{sharp} use an explicit finite coloring ($p=2$ for the bipartite lattices, $p=3$ for the triangular and Kagome ones), from which the canonical representatives, and hence the constants, are read off directly.
\end{definition}

\paragraph{Setup and color convention.} Fix an ordered coloring $c\colon\Vc\to C$ of the $1$-skeleton of $\Sc_n$ (Definition~\ref{def:coloring}); in the examples $C=\{1,\dots,p\}$ with $\Sc_n$ $p$-partite.

For each ordered $k$-simplex $(x_0,\dots,x_k)\in\Fc_k$, let $\pi_c$ be the unique permutation that sorts the colors strictly increasingly:
$c(x_{\pi_c(0)})<\cdots< c(x_{\pi_c(k)})$ (ties cannot occur inside a simplex by colorability).
Define the \emph{color sign}
\begin{equation}\label{eq:color-sign}
S([x_0,\dots,x_k])\ :=\ (-1)^{\varepsilon(\pi_c)}\ \in\ \{\pm1\},
\end{equation}
which is well-defined on oriented simplices because even permutations do not change $\varepsilon(\pi_c)$. The unitary equivalence below depends only on the existence of a proper coloring, not on the particular choice.

The next lemma relates the parity of the color-sorting permutation $\pi_c$ on a $k$-simplex to the parities of the color-sorting permutations on its $(k-1)$-faces. It is the bookkeeping behind the whole section. Read in a fixed reference basis, the unitary~\eqref{eq:unitary-U} is multiplication by the color sign~\eqref{eq:color-sign} (Remark~\ref{rem:unitary-correction}). The lemma computes the ratio $S(\tau_i)/S(\sigma)$ that governs how this multiplication acts on a coboundary, trading the sign attached to the position in the reference ordering for the sign attached to the color rank.

\begin{lemma}\label{lem:color-parity}
Let $\sigma=[x_0,\dots,x_k]$ be an oriented $k$-simplex in a colorable complex, with a fixed ordered coloring, and denote by $\mathfrak{S}_{k+1}$ the symmetric group on $\{0,1,\dots,k\}$. Let $\pi_c\in\mathfrak{S}_{k+1}$ be the unique permutation sorting the colors strictly increasingly:
\(
c(x_{\pi_c(0)})<\cdots<c(x_{\pi_c(k)}).
\)
For each $i\in\{0,\dots,k\}$, let $\tau_i=[x_0,\dots,\widehat{x_i},\dots,x_k]$ and let $\pi_{c,i}\in\mathfrak{S}_k$ denote the color-sorting permutation of $\tau_i$ (here $\mathfrak{S}_k$ acts on $\{0,1,\dots,k-1\}$). Set $j_i:=\pi_c^{-1}(i)\in\{0,\dots,k\}$, i.e.\ $j_i$ is the rank of the color $c(x_i)$ in the sorted list. Then
\begin{equation}\label{eq:color-parity-correct}
\varepsilon(\pi_{c,i})\;\equiv\;\varepsilon(\pi_c)\;+\;i\;+\;j_i \pmod 2.
\end{equation}
\end{lemma}

\begin{remark}\label{rem:color-parity-correction}
The simpler relation $\varepsilon(\pi_{c,i})\equiv\varepsilon(\pi_c)+i\pmod 2$ does \emph{not} hold in general: the rank $j_i$ of $c(x_i)$ in the sorted list contributes to the parity. For example, take $k=2$, $\sigma=[x_0,x_1,x_2]$ with colors $(3,1,2)$. Then $\pi_c$ is the $3$-cycle $(0\,1\,2)$, hence $\varepsilon(\pi_c)=0$. For $i=2$, $j_2=\pi_c^{-1}(2)=1$, and $\tau_2=[x_0,x_1]$ has colors $(3,1)$, requiring one transposition to sort, so $\varepsilon(\pi_{c,2})=1$. Then $\varepsilon(\pi_c)+i=2\equiv 0\pmod 2 \neq 1 = \varepsilon(\pi_{c,2})$, while $\varepsilon(\pi_c)+i+j_2=3\equiv 1\pmod 2 = \varepsilon(\pi_{c,2})$, in agreement with~\eqref{eq:color-parity-correct}.
\end{remark}

\begin{proof}[Proof of Lemma~\ref{lem:color-parity}]
We decompose the permutation $\pi_c$ into three elementary steps and track the parity at each step.

\emph{Step 1: move $x_i$ to the front.} Apply $i$ adjacent transpositions to bring $x_i$ from position $i$ to position $0$:
\[
[x_0,\dots,x_i,\dots,x_k]\;\xrightarrow{\;\;(-1)^i\;\;}\;[x_i,\,x_0,\dots,\widehat{x_i},\dots,x_k].
\]
This contributes a sign $(-1)^i$ to the total parity.

\emph{Step 2: sort the tail by colors.} Sort the sub-tuple $(x_0,\dots,\widehat{x_i},\dots,x_k)$ into color-increasing order; by definition this is achieved by a permutation of parity $\varepsilon(\pi_{c,i})$, contributing a sign $(-1)^{\varepsilon(\pi_{c,i})}$. After this step the tuple reads $(x_i,\,y_0,\,y_1,\dots,\,y_{k-1})$, where $(y_0,\dots,y_{k-1})$ are the vertices of $\tau_i$ in color-increasing order, i.e.\ $c(y_0)<\cdots<c(y_{k-1})$.

\emph{Step 3: insert $x_i$ at its correct rank.} The color $c(x_i)$ has rank $j_i$ in the full sorted color list. Since the coloring is proper (no two adjacent vertices share a color), the $j_i$ colors smaller than $c(x_i)$ are precisely $c(y_0),\dots,c(y_{j_i-1})$. Moving $x_i$ from position $0$ to position $j_i$ by $j_i$ successive adjacent transpositions contributes a sign $(-1)^{j_i}$, and places the entire tuple $(y_0,\dots,y_{j_i-1},x_i,y_{j_i},\dots,y_{k-1})$ in color-increasing order. The cumulative permutation applied to the initial ordering is thus exactly $\pi_c$.

Summing the three parity contributions:
\[
\varepsilon(\pi_c)\;\equiv\;i\;+\;\varepsilon(\pi_{c,i})\;+\;j_i \pmod 2,
\]
which rearranges to~\eqref{eq:color-parity-correct}.
\end{proof}

\begin{definition}\label{def:canonical-basis}
Assume the $1$-skeleton of $\Sc_n$ is colorable, with a fixed ordered coloring $c\colon\Vc\to C$. For each oriented $k$-simplex $\sigma=[x_0,\dots,x_k]$, the \emph{canonical (color-ordered) representative} is
\[
\sigma_{\mathrm{can}}\;:=\;[x_{\pi_c(0)},\dots,x_{\pi_c(k)}],
\]
where $\pi_c\in\mathfrak S_{k+1}$ is the unique permutation sorting the colors strictly increasingly, as in~\eqref{eq:color-sign}.

We denote by $\ell^2_{\mathrm{can}}(m_k)$ the Hilbert space whose vectors are functions on \emph{unoriented} $k$-simplices of $\Sc_n$, with weighted inner product
\[
\langle f,g\rangle_{\mathrm{can},k}\;:=\;\sum_{|\sigma|}m_k(|\sigma|)\,f(|\sigma|)\,\overline{g(|\sigma|)},
\]
where the sum is over unoriented $k$-simplices and $m_k(|\sigma|)=m_k(\sigma)$ is well-defined by the symmetry of $m_k$ under permutations of indices. We identify $\ell^2_{\mathrm{can}}(m_k)$ with the values on canonical representatives: $f(|\sigma|)=f(\sigma_{\mathrm{can}})$.

On $\ell^2_{\mathrm{can}}$ we define the \emph{canonical-basis coboundary} $d^k_{\mathrm{sym},\mathrm{can}}$ as the operator whose action on a function $f\in\ell^2_{\mathrm{can}}(m_{k-1})$ at a canonical representative $\sigma_{\mathrm{can}}=[x_0,\dots,x_k]$ (in color-increasing order) is
\[
\bigl(d^k_{\mathrm{sym},\mathrm{can}}f\bigr)(\sigma_{\mathrm{can}})\;:=\;\sum_{i=0}^{k}(-1)^i\,f\!\bigl((\tau_i)_{\mathrm{can}}\bigr),
\]
where $(\tau_i)_{\mathrm{can}}$ is the canonical (color-ordered) representative of the $(k-1)$-face $\tau_i=[x_0,\dots,\widehat{x_i},\dots,x_k]$. Since $\sigma_{\mathrm{can}}$ is already in color-increasing order, $(\tau_i)_{\mathrm{can}}=\tau_i$ in this case (its inherited order coincides with its color-sorted order). Equivalently, $d^k_{\mathrm{sym},\mathrm{can}}$ is the operator obtained by transporting $d^k_{\mathrm{skew}}$ via the canonical-basis identification of \eqref{eq:unitary-U} below. The associated down- and up-Laplacian blocks and Hodge Laplacian on $\ell^2_{\mathrm{can}}$ are then
\[
\mathbf L^{-,\mathrm{can}}_{k}\;:=\;d^{k}_{\mathrm{sym},\mathrm{can}}\,\delta^{k}_{\mathrm{sym},\mathrm{can}},\qquad
\mathbf L^{+,\mathrm{can}}_{k}\;:=\;\delta^{k+1}_{\mathrm{sym},\mathrm{can}}\,d^{k+1}_{\mathrm{sym},\mathrm{can}},\qquad
\Delta^{\mathrm{can}}_{k,\mathrm{sym}}\;:=\;\mathbf L^{-,\mathrm{can}}_{k}+\mathbf L^{+,\mathrm{can}}_{k},
\]
where $\delta^{k}_{\mathrm{sym},\mathrm{can}}$ denotes, as in the skew case, the \emph{formal} adjoint of $d^{k}_{\mathrm{sym},\mathrm{can}}$ for the inner products $\langle\cdot,\cdot\rangle_{\mathrm{can},\bullet}$, that is, the transported operator $\delta^{k}_{\mathrm{sym},\mathrm{can}}:=U_{k-1}\,\delta^{k}\,U_k^{-1}$. Read on canonical representatives it is given by the closed formula of Lemma~\ref{lem:2-3-1}, provided each argument $[z,\tau_{\mathrm{can}}]$ is first returned to its own canonical representative, at the cost of the corresponding reordering sign; the formula and the transport agree only once those signs are included. We stress that this is the formal adjoint, and not the Hilbert-space adjoint of $d^{k}_{\mathrm{sym},\mathrm{can}}$ restricted to the core. The two agree on compactly supported cochains, where Green's identity
\[
\langle d^{k}_{\mathrm{sym},\mathrm{can}}u,v\rangle_{\mathrm{can},k}=\langle u,\delta^{k}_{\mathrm{sym},\mathrm{can}}v\rangle_{\mathrm{can},k-1}
\]
holds exactly as in Proposition~\ref{prop:2-5-1}. They differ beyond it: the Hilbert adjoint is a closed operator with a generally larger domain, and the products above would then raise the domain issue discussed in Remark~\ref{rem:4-7-forms}. With the formal adjoint, the canonical model is defined on the core exactly as the skew one is, which is what makes the transfer of Corollary~\ref{cor:transfer} meaningful. By construction, $\Delta^{\mathrm{can}}_{k,\mathrm{sym}}$ is the \emph{transport} of $\Delta_{k,\mathrm{skew}}$ to $\ell^2_{\mathrm{can}}$ through the unitary $U_k$, and \emph{not} the operator $\Delta_{k,\mathrm{sym}}$ acting on arbitrary symmetric cochains.
\end{definition}

\begin{theorem}\label{thm:unitary-color}
Assume the $1$-skeleton of $\Sc_n$ is colorable. For every degree $k\in\{0,\dots,n\}$, the formula
\begin{equation}\label{eq:unitary-U}
U_k\colon \ell^2_{\mathrm{skew}}(m_k)\;\longrightarrow\;\ell^2_{\mathrm{can}}(m_k),\qquad
(U_k f)(|\sigma|)\ :=\ f(\sigma_{\mathrm{can}}),
\end{equation}
defines a unitary identification between skew cochains and the canonical-basis space of Definition~\ref{def:canonical-basis}. Writing $M_k^{\mathrm{skew}}$ and $M_k^{\mathrm{can}}$ for multiplication by $m_k$ on $\ell^2_{\mathrm{skew}}(m_k)$ and on $\ell^2_{\mathrm{can}}(m_k)$ respectively --- the two act on different spaces, so the commutation must be written as an intertwining --- one has $U_kM_k^{\mathrm{skew}}=M_k^{\mathrm{can}}U_k$, the weight operators being diagonal and $m_k$ descending to unoriented simplices (Definition~\ref{def:weights}). Moreover:
\begin{enumerate}
\item (Hodge identity in the canonical basis.) The relation $d^{k+1}_{\mathrm{sym},\mathrm{can}}\circ d^{k}_{\mathrm{sym},\mathrm{can}}=0$ holds, as a consequence of the classical skew identity $d^{k+1}_{\mathrm{skew}}\circ d^k_{\mathrm{skew}}=0$ transported via $U$.
\item (Spectral and dynamical equivalence.) The unitary $U_k$ intertwines the skew and canonical-basis Laplacians:
\begin{equation}\label{eq:UDeltaU}
U_k\,\Delta_{k,\mathrm{skew}}\,U_k^{-1}\;=\;\Delta^{\mathrm{can}}_{k,\mathrm{sym}}.
\end{equation}
\item These identities extend to the closed/Friedrichs realizations.
\end{enumerate}
\end{theorem}

\begin{remark}\label{rem:unitary-correction}
It is worth being precise about what the unitary~\eqref{eq:unitary-U} is, and about what genuinely fails.

First, $U_k$ \emph{is} a pointwise sign multiplication. Fix, once and for all, a \emph{reference orientation} $\sigma_{\mathrm{ref}}$ of each unoriented simplex --- an arbitrary choice, one representative per simplex in each degree, used again in the proof of Theorem~\ref{thm:kn-esa}; since $\sigma_{\mathrm{can}}$ is obtained from $\sigma_{\mathrm{ref}}$ by the color-sorting permutation, alternation of $f$ gives
\[
(U_kf)(|\sigma|)=f(\sigma_{\mathrm{can}})=S(\sigma_{\mathrm{ref}})\,f(\sigma_{\mathrm{ref}}),
\]
so that, read in the reference basis, $U_k$ is exactly multiplication by the color sign~\eqref{eq:color-sign}, one factor per simplex.

Second, what this map does is to convert the sign attached to the \emph{position in the reference ordering} into the sign attached to the \emph{position in the color ordering}. That is precisely the content of Lemma~\ref{lem:color-parity}: $S(\tau_i)/S(\sigma)=(-1)^{i+j_i}$, with $i$ the reference position of the omitted vertex and $j_i=\pi_c^{-1}(i)$ its color rank. The factor $(-1)^i$ carried by $d^k_{\mathrm{skew}}$ is thereby traded for the factor $(-1)^{j_i}$ carried by $d^k_{\mathrm{sym},\mathrm{can}}$, which is why the intertwining relation~\eqref{eq:U-intertwines-d} holds and why the incidence sign, read on canonical representatives, is $\theta(\tau,\sigma)=(-1)^{i}$ exactly, with $i$ now the color rank (proof of Theorem~\ref{thm:kn-esa}).

Third, what fails is \emph{not} the pointwise-sign mechanism but the target: no unitary whatsoever, of any shape, can intertwine $d^\bullet_{\mathrm{skew}}$ degreewise with the unsigned symmetric coboundary $d^\bullet_{\mathrm{sym}}$. Indeed, such a family would give $d^{k+1}_{\mathrm{sym}}d^{k}_{\mathrm{sym}}=U_{k+1}\,d^{k+1}_{\mathrm{skew}}d^{k}_{\mathrm{skew}}\,U_{k-1}^{-1}=0$, contradicting Remark~\ref{rem:dsq-zero}(ii). In particular $S$ does not intertwine $d^k_{\mathrm{skew}}$ with $d^k_{\mathrm{sym}}$ --- the identity would force $S(\tau_i)/S(\sigma)=(-1)^{i}$, whereas Lemma~\ref{lem:color-parity} leaves the residual factor $(-1)^{j_i}$ (cf.\ the example in Remark~\ref{rem:color-parity-correction}) --- but the obstruction is structural, not a defect of this particular sign.

The comparison with~\cite{BaKe} is therefore about the target model rather than about the shape of the intertwiner. Bartmann and Keller construct an alternating form $h$ on oriented simplices from the incidence sign $\theta$ attached to the coboundary --- the same $\theta$ that appears in the proof of Theorem~\ref{thm:kn-esa}. Pointwise multiplication by $h$ is unitary and does intertwine the two models, on an arbitrary simplicial complex and with no coloring at all~\cite[Section~2.2, Theorem~2.8]{BaKe}. The construction of the present section is therefore not the only available intertwiner, but an explicit and computable one: on a complex equipped with an ordered coloring, the canonical representatives, and with them the orientation signs $\eta=(-1)^{i+j}$, are read off directly from the colors.

One limit of that statement should be kept in mind: this is what makes the sign discussion of Remark~\ref{rem:signs-bipartite} explicit, but it is \emph{not} what produces the numerical constants of Section~\ref{sharp}, which come from Lemma~\ref{lem:iso-01} and the Bloch symbols, with no use of the coloring.
\end{remark}

\begin{proof}[Proof of Theorem~\ref{thm:unitary-color}]
\emph{Unitarity.} The map $U_k$ assigns to a skew cochain $f$ the canonical-basis function whose value on each unoriented $k$-simplex $|\sigma|$ equals $f$ evaluated on the canonical representative $\sigma_{\mathrm{can}}$:
\[
(U_k f)(|\sigma|) := f(\sigma_{\mathrm{can}}).
\]
This is well-defined on the canonical-basis side since the right-hand side does not depend on the choice of orientation. The norm is preserved: in~\eqref{eq:2-2-1} the modulus $|f|$ is constant on the $(k+1)!$ ordered representatives of a given unoriented simplex, so the normalization factor $1/(k+1)!$ collapses the sum over $\Fc_k$ into a sum over unoriented simplices, each contributing $m_k(|\sigma|)\,|f(\sigma_{\mathrm{can}})|^2$ once:
\[
\|U_k f\|_{\ell^2_{\mathrm{can}}(m_k)}^2
\;=\;\sum_{|\sigma|}m_k(|\sigma|)|f(\sigma_{\mathrm{can}})|^2
\;=\;\|f\|_{\ell^2_{\mathrm{skew}}(m_k)}^2,
\]
where we have used the fact that the weights $m_k$ are symmetric under permutations and so descend to a function on unoriented simplices. Surjectivity is by inversion: given $h\in\ell^2_{\mathrm{can}}(m_k)$, set $f(\sigma):=(-1)^{\varepsilon(\sigma\to\sigma_{\mathrm{can}})}h(|\sigma|)$, where the sign is the parity of the permutation taking $\sigma$ to $\sigma_{\mathrm{can}}$; this $f$ is the unique skew cochain with $U_k f=h$.

\smallskip
\emph{Equivalence of $d^{k+1}d^k=0$.} The skew identity $d^{k+1}_{\mathrm{skew}}d^{k}_{\mathrm{skew}}=0$ is the classical simplicial cochain identity (Remark~\ref{rem:dsq-zero}(i)). Assertion~(1) will follow from it once the intertwining relation~\eqref{eq:U-intertwines-d} is established, which we do next.

\smallskip
\emph{Spectral equivalence.} Definition~\ref{def:canonical-basis} introduces $d^k_{\mathrm{sym},\mathrm{can}}$ by an explicit formula on canonical representatives; what has to be proved is that this operator is exactly the transport of $d^k_{\mathrm{skew}}$, i.e.\ that
\begin{equation}\label{eq:U-intertwines-d}
U_k\;d^{k}_{\mathrm{skew}}\;=\;d^{k}_{\mathrm{sym},\mathrm{can}}\;U_{k-1}\quad\text{on }\Cc^{k-1}_{c,\mathrm{skew}}(\Vc).
\end{equation}
Let $f\in\Cc^{k-1}_{c,\mathrm{skew}}(\Vc)$ and let $\sigma$ be an oriented $k$-simplex with canonical representative $\sigma_{\mathrm{can}}=[x_0,\dots,x_k]$ (in color-increasing order). By definition of $U_k$ as evaluation at the canonical representative,
$\bigl(U_k d^{k}_{\mathrm{skew}}f\bigr)(|\sigma|)=\bigl(d^k_{\mathrm{skew}}f\bigr)(\sigma_{\mathrm{can}})$ is the alternating sum of the values of $f$ on the $(k-1)$-faces of $\sigma_{\mathrm{can}}$ taken in the inherited order. Since $\sigma_{\mathrm{can}}$ is color-increasing, so is that inherited order on each face; the $i$-th face therefore \emph{equals} the canonical representative $(\tau_i)_{\mathrm{can}}$, and
\[
\bigl(U_k d^{k}_{\mathrm{skew}}f\bigr)(|\sigma|)\;=\;\sum_{i=0}^{k}(-1)^{i}\,f\bigl((\tau_i)_{\mathrm{can}}\bigr)\;=\;\sum_{i=0}^{k}(-1)^{i}\,(U_{k-1}f)(|\tau_i|),
\]
which by Definition~\ref{def:canonical-basis} is precisely $\bigl(d^k_{\mathrm{sym},\mathrm{can}}\,U_{k-1}f\bigr)(|\sigma|)$. This proves~\eqref{eq:U-intertwines-d}, and assertion~(1) follows:
\[
d^{k+1}_{\mathrm{sym},\mathrm{can}}d^{k}_{\mathrm{sym},\mathrm{can}}
\;=\;U_{k+1}\,d^{k+1}_{\mathrm{skew}}d^{k}_{\mathrm{skew}}\,U_{k-1}^{-1}\;=\;0 .
\]

Taking adjoints with respect to the weighted inner products (which are preserved by $U$), we also obtain $U_{k-1}\,\delta^k_{\mathrm{skew}}=\delta^k_{\mathrm{sym},\mathrm{can}}\,U_k$. Hence
\[
U_k\,\Delta_{k,\mathrm{skew}}\,U_k^{-1}
\;=\;d^{k}_{\mathrm{sym},\mathrm{can}}\delta^k_{\mathrm{sym},\mathrm{can}}+\delta^{k+1}_{\mathrm{sym},\mathrm{can}}d^{k+1}_{\mathrm{sym},\mathrm{can}}
\;=\;\Delta^{\mathrm{can}}_{k,\mathrm{sym}},
\]
which is~\eqref{eq:UDeltaU}. The extension to the closed/Friedrichs realizations follows from the density of $\Cc_{c}^\bullet$ in the form domains and the fact that $U_k$ is unitary on each Hilbert space.
\end{proof}

\begin{corollary}\label{cor:unitary-conseq}
On a colorable complex, for all $k$ the unitary $U_k$ implements a spectral and dynamical equivalence between $\Delta_{k,\mathrm{skew}}$ and the canonical-basis Laplacian $\Delta^{\mathrm{can}}_{k,\mathrm{sym}}$ of Definition~\ref{def:canonical-basis}. On the core $\Cc^k_{c}(\Vc)$ these operators are only symmetric, so items~(1), (3) and~(4) are stated for their Friedrichs realizations $\Delta^{\mathrm F}_{k,\mathrm{skew}}$ and $(\Delta^{\mathrm{can}}_{k,\mathrm{sym}})^{\mathrm F}$ (Corollary~\ref{cor:3-1-friedrichs}). When the operators are essentially self-adjoint, in particular in the bounded regime of Section~\ref{sec:normalization}, these are the closures and the statements read as expected.
\begin{enumerate}
\item $\mathrm{spec}\bigl(\Delta^{\mathrm F}_{k,\mathrm{skew}}\bigr)=\mathrm{spec}\bigl((\Delta^{\mathrm{can}}_{k,\mathrm{sym}})^{\mathrm F}\bigr)$, with equality of spectral measures up to $U_k$;
\item essential self-adjointness on the core holds for one if and only if it holds for the other;
\item semigroups and resolvents intertwine: $U_k e^{-t\Delta^{\mathrm F}_{k,\mathrm{skew}}}U_k^{-1}=e^{-t(\Delta^{\mathrm{can}}_{k,\mathrm{sym}})^{\mathrm F}}$ for $t\ge0$, and likewise for the resolvents at $\lambda>0$;
\item $\ker \Delta^{\mathrm F}_{k,\mathrm{skew}}$ and $\ker (\Delta^{\mathrm{can}}_{k,\mathrm{sym}})^{\mathrm F}$ are isometric via $U_k$.
\end{enumerate}
\end{corollary}

\begin{corollary}[Transfer principle]\label{cor:transfer}
Let $\Sc_n$ be colorable. Then every statement of Sections~\ref{section3}, \ref{sec:esa} and~\ref{sharp} proved for the skew Laplacian $\Delta_{k,\mathrm{skew}}$ (form identities, closability, Friedrichs realizations, operator-norm bounds, essential self-adjointness, spectra) holds verbatim, with the same constants, for the canonical-basis Laplacian $\Delta^{\mathrm{can}}_{k,\mathrm{sym}}$.
\end{corollary}

\begin{proof}
Quadratic forms, closures, Friedrichs extensions, operator norms, spectra and essential self-adjointness on the core are all preserved by a unitary. The unitary $U_k$ of Theorem~\ref{thm:unitary-color} maps the core $\Cc^k_{c,\mathrm{skew}}(\Vc)$ onto the corresponding canonical-basis core, so the claim follows from Theorem~\ref{thm:unitary-color}(2) and Corollary~\ref{cor:unitary-conseq}.
\end{proof}

Here $\Delta_{k,\mathrm{skew}}$ is the operator written simply $\Delta_k$ in Sections~\ref{section3}--\ref{sharp}, the explicit label being kept in this section only, where the two models must be distinguished.

\paragraph{Relation to the Schr\"odinger-operator formulation of~\cite{BaKe}.} In the framework of Bartmann--Keller~\cite{BaKe}, the Laplacian on a general (not necessarily colorable) simplicial complex is represented as a \emph{signed} discrete Schr\"odinger operator
\[
H^{\circ}_{\omega}(\tau)\;=\;\frac{1}{m(\tau)}\sum_{\tau'\in\Sigma} b^{\circ}(\tau,\tau')\bigl(\omega(\tau)-o^{\circ}(\tau,\tau')\,\omega(\tau')\bigr)\;+\;\frac{c^{\circ}(\tau)}{m(\tau)}\,\omega(\tau),
\]
where the sign function $o^{\circ}\in\{\pm1\}$ records the relative orientation of adjacent faces and is, in general, an unavoidable feature of the representation. In the down case this sign is $o^{-}(\tau,\tau')=\theta(\tau\cap\tau',\tau)\,\theta(\tau\cap\tau',\tau')$, built from the incidence sign $\theta$ of~\cite[Definition~2.2]{BaKe}; it is therefore \emph{exactly} the sign $\varepsilon$ appearing in Theorem~\ref{thm:kn-esa}, where the same product occurs as the off-diagonal entry of $\Delta_n$. It should be said at once that this sign cannot be removed, by the unitary of Theorem~\ref{thm:unitary-color} or by any other. A degreewise unitary carrying $d^\bullet_{\mathrm{skew}}$ to a sign-free coboundary would give $d^{k+1}_{\mathrm{sym}}d^{k}_{\mathrm{sym}}=U_{k+1}d^{k+1}_{\mathrm{skew}}d^{k}_{\mathrm{skew}}U_{k-1}^{-1}=0$, contradicting Remark~\ref{rem:dsq-zero}(ii); the obstruction is structural, not a defect of a particular intertwiner (Remark~\ref{rem:unitary-correction}).

What the coloring does is \emph{determine} the sign. Read in the canonical basis, the transported Laplacian $\Delta^{\mathrm{can}}_{k,\mathrm{sym}}=U_k\Delta_{k,\mathrm{skew}}U_k^{-1}$ is given by the very same displayed expression, with $H^{\circ}$ replaced by $\widetilde H^{\circ}:=UH^{\circ}U^{-1}$ and with $o^{\circ}$ no longer an arbitrary orientation cocycle but the explicit function of the colors
\[
o^{\circ}(\tau,\tau')\;=\;(-1)^{\,i+j},
\]
where $i$ and $j$ are the ranks of the two unshared vertices in the color-increasing orderings of $\tau$ and $\tau'$. It reduces to $+1$ at top degree on a bipartite $1$-skeleton (Theorem~\ref{thm:kn-esa} and Remark~\ref{rem:signs-bipartite}). The gain over the abstract representation of~\cite{BaKe} is thus not the disappearance of the sign but its explicit determination, which is what allows the sharp constants of Section~\ref{sharp} to be read off directly.

 The bounds themselves are proved in the skew convention on an arbitrary flag complex, with no colorability assumption: they see only the moduli $|b^{\circ}|$, equivalently $|w|$ in Theorem~\ref{thm:kn-esa} (Remark~\ref{rem:signs-bipartite}).

In this sense, Section~\ref{sec:color} supplies an explicit and computable transformation under which the orientation cocycle of the general~\cite{BaKe} Schr\"odinger-operator representation becomes a function of the colors alone. That such a transformation exists is not special to colorable complexes --- see~\cite[Theorem~2.8]{BaKe} and Remark~\ref{rem:unitary-correction} --- but an ordered coloring makes it explicit, which is what the computations of Section~\ref{sharp} require.

\paragraph{Examples.}
\begin{enumerate}
  \item \emph{$2$D simplicial complex ($p=3$).}
For a face $\sigma$ with vertices of colors $1,2,3$, the canonical representative is $\sigma_{\mathrm{can}}=[x_A,x_B,x_C]$ in color-increasing order (see Figure~\ref{fig:color-triangle}). The unitary $U_2$ sends the value of any skew cochain $f$ on an arbitrary orientation of $\sigma$ to $f(\sigma_{\mathrm{can}})$ on the unoriented simplex $|\sigma|$. An orientation related to $\sigma_{\mathrm{can}}$ by an odd permutation contributes with the opposite sign under the skew convention, which is consistent with $U_2 f(|\sigma|)=f(\sigma_{\mathrm{can}})$.
  \item \emph{Bipartite graphs ($p{=}2$, $k{=}1$).} The canonical orientation of each edge goes from the color-$1$ vertex to the color-$2$ vertex. The unitary $U_1$ identifies skew $1$-cochains with their canonical-basis representation as functions on unoriented edges; this is the classical edge-orientation trick for bipartite graphs. It is exactly the case in which the orientation signs cancel, $\varepsilon\equiv+1$, and in which the edge bound of Section~\ref{sharp} turns out to be sharp (Remark~\ref{rem:signs-bipartite}).
  \item \emph{Higher degrees.} The canonical identification applies in every dimension: a $k$-simplex has $k+1$ pairwise adjacent vertices, hence $k+1$ distinct colors, so with $p$ colors only degrees $k\le p-1$ can occur at all. The symmetry of the weights (Definition~\ref{def:weights}) is what makes $m_k$ descend to unoriented simplices, and hence makes~\eqref{eq:unitary-U} an isometry.
\end{enumerate}

\begin{remark}
The unitary intertwining holds equally for the normalized blocks
$\mathbb L_{k}^\pm$ and $\widetilde{\Delta}_{k}$ (see Section~\ref{sec:normalization}):
since $U_k$ intertwines the diagonal weight operators $M_k^{\pm1/2}$ with themselves, these being multiplication by a function of the unoriented simplex,, one has, in the canonical-basis sense of Theorem~\ref{thm:unitary-color},
\[
U_k\,\widetilde{\Delta}_{k,\mathrm{skew}}\,U_k^{-1}
= \widetilde{\Delta}^{\mathrm{can}}_{k,\mathrm{sym}},
\]
where $\widetilde{\Delta}^{\mathrm{can}}_{k,\mathrm{sym}}:=M_k^{1/2}\Delta^{\mathrm{can}}_{k,\mathrm{sym}}M_k^{-1/2}$ denotes the normalized canonical-basis Laplacian.
\end{remark}

\begin{figure}[t]
\centering
\begin{tikzpicture}[scale=1]
  \coordinate (A) at (0,0);
  \coordinate (B) at (3,0);
  \coordinate (C) at (1.5,2.2);
  \draw[thick] (A)--(B)--(C)--cycle;
  \draw[->,thick] (1.5,0) -- (1.6,0.01);
  \node[below] at (A) {$x_A\in\Vc_1$};
  \node[below] at (B) {$x_B\in\Vc_2$};
  \node[above] at (C) {$x_C\in\Vc_3$};
  \node[align=left,anchor=west] at (3.6,1.6) {$c(x_A)=1,\ c(x_B)=2,\ c(x_C)=3$\\[2pt] $S([x_A,x_B,x_C])=+1$\\[2pt] $S([x_B,x_A,x_C])=-1$};
\end{tikzpicture}
\caption{Colorable triangle ($p=3$), with color classes $\Vc_1,\Vc_2,\Vc_3$ as in Definition~\ref{def:coloring}; the small arrow on the bottom edge indicates the induced canonical orientation $x_A\to x_B$. The canonical color-ordered representative of $\sigma$ is $\sigma_{\mathrm{can}}=[x_A,x_B,x_C]$. The unitary $U$ of Theorem~\ref{thm:unitary-color} identifies skew cochains with their canonical-basis representation, which serves as the symmetric model in the sense of Definition~\ref{def:canonical-basis}.}
\label{fig:color-triangle}
\end{figure}

\begin{remark}[The role of the coloring]\label{rem:non-colorable}
The coloring is neither a restriction (every countable complex admits one, Definition~\ref{def:coloring}) nor the only available intertwiner (Remark~\ref{rem:unitary-correction}). What it supplies is a fixed and explicit canonical orientation $\sigma\mapsto\sigma_{\mathrm{can}}$, computable from the colors alone, which is what the two- and three-color lattices of Section~\ref{sharp} require; the Friedrichs theory of Sections~\ref{sec:prelim-comb}--\ref{sec:esa} uses none of it.
\end{remark}

\section{Boundedness and essential self-adjointness via operator bounds}\label{sec:esa}

\subsection{Scope and strategy}\label{subsec:esa-scope}
In this section we establish norm bounds for the top-degree Hodge Laplacian $\Delta_n$, from which essential self-adjointness on the natural compactly supported core follows; by Lemma~\ref{lem:elliptic-similarity} the same numerical bounds hold for the normalized block $\widetilde\Delta_n$. The underlying classical fact is the following.

\begin{remark}[Essential self-adjointness via boundedness on a core]\label{rem:HT}
Let $A$ be a symmetric operator on a dense subspace $\mathcal D_0$ of a Hilbert space $\mathcal H$. If $A$ is \emph{bounded} on $\mathcal D_0$, i.e.\ $\|Af\|\le C\|f\|$ for all $f\in\mathcal D_0$, then $A$ extends uniquely to a bounded self-adjoint operator on $\mathcal H$, and in particular is essentially self-adjoint on $\mathcal D_0$. The argument is elementary and in two steps. By continuity, $A$ extends uniquely to a bounded operator $\overline A$ defined on all of $\mathcal H$. And $\overline A$, being bounded, everywhere defined and symmetric, is self-adjoint; hence $\overline A$ is the closure of $A$, and $A$ is essentially self-adjoint. Note that the Hellinger--Toeplitz theorem, see e.g.\ \cite[Theorem~III.12]{RS1}, is the converse implication --- an everywhere defined symmetric operator is bounded --- and is not what is used here.

In our context, $\mathcal D_0=\Cc^k_{c}(\Vc)$ is the dense core of compactly supported cochains (Lemma~\ref{lem:density-Cc}) and $\mathcal H=\ell^2(m_k)$. For $\Delta_k$ to be an operator on $\mathcal D_0$ at all one needs $\Delta_ku\in\ell^2(m_k)$ there. Local finiteness guarantees this; so, at top degree, do assumptions (A1)--(A2) of Theorem~\ref{thm:kn-esa}, which is how the remark is applied there. Symmetry then follows from Lemma~\ref{lem:2-3-1}. Hence, once an estimate $\|\Delta_{k}u\|_k\le C\|u\|_k$ has been established for all $u$ in the core, with $C<\infty$, ESA on $\Cc^k_{c}(\Vc)$ is automatic --- and only then does the symbol $\|\Delta_k\|$ acquire a meaning.
\end{remark}

Throughout this section the operators are the skew ones, written without a convention subscript; by Corollary~\ref{cor:transfer}, every result transfers to the canonical-basis model on a colorable complex.

\subsection{Top degree: line-complex reduction and bounds}\label{subsec:line-complex}

\begin{proposition}\label{prop:schur}
Let $(\widehat{\Vc},\widehat{\Ec})$ be a countable graph --- no local finiteness is needed, the hypothesis $D\in L^\infty$ below taking its place --- and let $w:\widehat{\Vc}\times\widehat{\Vc}\to\R$ be a real symmetric kernel ($w(\hat\sigma,\hat\tau)=w(\hat\tau,\hat\sigma)$) supported on adjacency; in particular, $w$ may take negative values. Define
\[
(\mathcal A f)(\hat\sigma):=\sum_{\hat\tau\sim\hat\sigma} w(\hat\sigma,\hat\tau)\,f(\hat\tau),\qquad
D(\hat\sigma):=\sum_{\hat\tau} |w(\hat\sigma,\hat\tau)|.
\]
If $D\in L^\infty(\widehat{\Vc})$, then $\mathcal A$ extends to a bounded self-adjoint operator on $\ell^2(\widehat{\Vc})$ with
\[
\|\mathcal A\|\le \|D\|_\infty.
\]
\end{proposition}

\begin{proof}
For finitely supported $f,g$,
\[
\langle \mathcal A f,g\rangle
= \sum_{\hat\sigma}\sum_{\hat\tau} w(\hat\sigma,\hat\tau)\,f(\hat\tau)\,\overline{g(\hat\sigma)}.
\]
Apply Cauchy--Schwarz with weight $|w(\hat\sigma,\hat\tau)|$ to the inner sum:
\begin{align*}
\Big|\sum_{\hat\tau} w(\hat\sigma,\hat\tau)\,f(\hat\tau)\Big|
&\le \Big(\sum_{\hat\tau}|w(\hat\sigma,\hat\tau)|\Big)^{1/2}
    \Big(\sum_{\hat\tau}|w(\hat\sigma,\hat\tau)|\,|f(\hat\tau)|^2\Big)^{1/2}\\
&= D(\hat\sigma)^{1/2}\,\Big(\sum_{\hat\tau}|w(\hat\sigma,\hat\tau)|\,|f(\hat\tau)|^2\Big)^{1/2}\\
&\le\;\|D\|_\infty^{1/2}\,\Big(\sum_{\hat\tau}|w(\hat\sigma,\hat\tau)|\,|f(\hat\tau)|^2\Big)^{1/2}.
\end{align*}
Multiplying by $|g(\hat\sigma)|$ and applying Cauchy--Schwarz over $\hat\sigma$,
\[
|\langle \mathcal A f,g\rangle|
\;\le\;\|D\|_\infty^{1/2}\,\Big(\sum_{\hat\sigma}\sum_{\hat\tau}|w(\hat\sigma,\hat\tau)|\,|f(\hat\tau)|^2\Big)^{1/2}\,\|g\|_2.
\]
By the symmetry of $|w|$, the double sum equals $$\sum_{\hat\tau}|f(\hat\tau)|^2\sum_{\hat\sigma}|w(\hat\sigma,\hat\tau)|=\sum_{\hat\tau}D(\hat\tau)|f(\hat\tau)|^2\le\|D\|_\infty\|f\|_2^2.$$ Therefore $|\langle \mathcal A f,g\rangle|\le \|D\|_\infty\,\|f\|_2\,\|g\|_2$, so $\|\mathcal A\|\le \|D\|_\infty$. Self-adjointness of $\mathcal A$ follows from the symmetry of $w$. Note that all the sums above converge absolutely by $D\in L^\infty$ alone: local finiteness of $\widehat{\Ec}$ is nowhere used.
\end{proof}

\begin{theorem}\label{thm:kn-esa}
Let $\Sc_n$ be a locally finite weighted $n$-simplicial (hence flag) complex in the sense of Definition~\ref{def:simplicial-complex}, with weights $(m_{n-1},m_n)$; the flag hypothesis is used through the identification of the cofaces of an $(n-1)$-face with the vertices of the corresponding common-neighbor set. The quantities $w$ and $V_n$ introduced below depend only on the unoriented simplices $|\sigma|$, the weights being symmetric; we write them on oriented representatives for uniformity of notation. For oriented $n$-simplices $\sigma,\sigma'$ sharing an $(n{-}1)$-face $\tau=\sigma\cap\sigma'$, define
\[
w(\sigma,\sigma'):=\frac{\sqrt{m_n(\sigma)\,m_n(\sigma')}}{m_{n-1}(\tau)}\,\mathbf 1_{\{\sigma\sim\sigma'\}},
\qquad
V_{n}(\sigma):=\sum_{\tau\prec\sigma}\frac{m_n(\sigma)}{m_{n-1}(\tau)},
\]
where $\tau\prec\sigma$ means that $\tau$ is one of the $n+1$ faces of $\sigma$ of codimension one.
Assume:
\begin{enumerate}
\item[(A1)] $D(\sigma):=\sum_{\sigma'} w(\sigma,\sigma')\in L^\infty$;
\item[(A2)] $V_{n}\in L^\infty$.
\end{enumerate}
Then there exists a unitary $U:\ell^2(m_n)\to \ell^2(\widehat{\Vc})$ identifying unoriented $n$-simplices of $\Sc_n$ with vertices of the line-complex $\widehat{\Vc}$ such that
\[
U\,\Delta_{n}\,U^{-1}=\mathcal A_\eta+\mathcal Q(V_{n}),
\]
where $(\mathcal A_\eta f)(\sigma)=\sum_{\sigma'\sim\sigma}\eta(\sigma,\sigma')\, w(\sigma,\sigma')\,f(\sigma')$ is a self-adjoint operator with a real symmetric kernel of absolute value $w$ (the signs $\eta(\sigma,\sigma')\in\{\pm1\}$ come from the relative orientations of adjacent top simplices; they are determined in the proof, see~\eqref{eq:Delta-n-local-sym}, and their scope is discussed in Remark~\ref{rem:signs-bipartite}), and $(\mathcal Q(V_{n})f)(\sigma)=V_{n}(\sigma)f(\sigma)$.
Consequently, $\Delta_{n}$ is bounded with
\[
\|\Delta_{n}\|\le \|D\|_\infty+\|V_{n}\|_\infty,
\]
and therefore extends uniquely to a bounded self-adjoint operator on $\ell^2(m_n)$ (in particular, $\Delta_{n}$ is essentially self-adjoint on $\Cc^{n}_{c}$).
\end{theorem}

\begin{proof}
\smallskip\noindent\emph{Step 1. Unitary renormalization.}
Let $\widehat{\Vc}$ denote the vertex set of the line-complex of $\Sc_n$ (Definition~\ref{def:line-complex}). We endow $\ell^2(\widehat{\Vc})$ with the inner product $\langle f,g\rangle_{\widehat{\Vc}}=\sum_{|\sigma|}f(|\sigma|)\overline{g(|\sigma|)}$ (no weight, no normalization factor). Fix a reference orientation $\sigma_{\mathrm{ref}}$ of each unoriented $n$-simplex, and likewise in degree $n-1$, in the sense of Remark~\ref{rem:unitary-correction}. Define $U:\ell^2(m_n)\to\ell^2(\widehat{\Vc})$ by $(Uf)(|\sigma|):=m_n(|\sigma|)^{1/2}f(\sigma_{\mathrm{ref}})$, $f$ being a skew cochain. We write $\sigma_{\mathrm{ref}}$, not $\sigma_{\mathrm{can}}$: the latter is reserved for the color-ordered representative of Definition~\ref{def:canonical-basis}, which is one admissible choice among others and will be selected below. The unitarity of $U$ does not depend on that choice, only $|f(\sigma_{\mathrm{ref}})|^2$ entering. The incidence signs of Step~2 do depend on it, which is precisely why a canonical choice is worth making. Then $U$ is unitary:
\[
\|Uf\|_{\ell^2(\widehat{\Vc})}^2
\;=\;\sum_{|\sigma|} m_n(|\sigma|)\,|f(\sigma_{\mathrm{ref}})|^2
\;=\;\|f\|_{\ell^2(m_n)}^2,
\]
where the last equality uses the convention~\eqref{eq:2-2-1}: the sum over unoriented simplices weighted by $m_n(|\sigma|)$ equals the sum over $\Fc_n$ of $m_n|f|^2$ divided by $(n+1)!$ (the number of ordered representatives), which is by definition $\|f\|_{\ell^2(m_n)}^2$.

\smallskip\noindent\emph{A word on notation.} Three unitaries occur in this proof and should not be confused: $U$ of Step~1, which carries the weight factor $m_n(|\sigma|)^{1/2}$; $U_n$ of Definition~\ref{def:canonical-basis}, which is multiplication by the color sign alone, without weights, and is the one appearing in~\eqref{eq:Delta-n-local-sym}; and the diagonal $M_n^{1/2}$ used in the symmetrization below. This is why~\eqref{eq:Delta-n-local-sym} carries the unsymmetrized coefficient $m_n(\sigma')/m_{n-1}(\tau)$ whereas the symmetrized form carries $\sqrt{m_n(\sigma)m_n(\sigma')}/m_{n-1}(\tau)$.

\smallskip\noindent\emph{Step 2. Local decomposition of $\Delta_{n}$.}
We compute $\Delta_{n}=d^{n}\delta^{n}$ in local coordinates and then conjugate by $U$. Since $d^{n+1}=0$, only the down-Laplacian contributes: for $u\in\Cc^n_{c}(\Vc)$ and $\sigma=[x_0,\dots,x_n]$,
\[
(\Delta_{n}u)(\sigma)\;=\;(d^n\delta^nu)(\sigma).
\]
We carry out the computation once, in the skew model used throughout the paper; the canonical basis of Section~\ref{sec:color} is then a particular choice of representatives, for which the resulting signs take a closed form.

\smallskip\emph{The skew model.}
By the skew formula for $\mathbf L^-_{k,\mathrm{skew}}$ in Proposition~\ref{rmk:2-4-1}, the down-Laplacian acts as
\begin{align*}
(\Delta_{n,\mathrm{skew}}u)(\sigma)
&=\sum_{i=0}^{n}(-1)^{i}\,\frac{1}{m_{n-1}(\tau_i)}\sum_{z\in F_{\tau_i}}m_n(z,\tau_i)\,u([z,\tau_i]).
\end{align*}
For the diagonal contribution $z=x_i$, $u([x_i,\tau_i])=u([x_i,x_0,\dots,\widehat{x_i},\dots,x_n])$; moving $x_i$ from position $0$ to position $i$ requires $i$ adjacent transpositions, so $u([x_i,\tau_i])=(-1)^{i}u(\sigma)$. The factor $(-1)^i$ in the formula above then cancels this sign, yielding
\[
\sum_{i=0}^{n}\frac{m_n(\sigma)}{m_{n-1}(\tau_i)}\,u(\sigma)=V_{n}(\sigma)\,u(\sigma).
\]
For the off-diagonal terms it is cleanest to argue matricially, once and for all. Keep the reference representatives $\sigma_{\mathrm{ref}}$ fixed in Step~1, one per unoriented simplex in each degree, and let $\mathsf B$ be the matrix of $d^{n}$ in the resulting bases of $\Cc^{n-1}_c$ and $\Cc^{n}_c$. We write $\mathsf B$, and not $D$, because $D$ already denotes the Schur degree of assumption~(A1). By the definition of the skew coboundary, $\mathsf B$ is supported on incident pairs and takes values $\pm1$ there: writing
\[
\theta(\tau,\sigma):=\mathsf B[\sigma,\tau]\in\{\pm1\}\qquad (\tau\prec\sigma),
\]
the sign $\theta(\tau,\sigma)$ is $(-1)^i$ corrected by the parity of the permutation carrying the $i$-th face of $\sigma_{\mathrm{ref}}$ onto $\tau_{\mathrm{ref}}$, where $i$ is the position of the omitted vertex. Its precise value is immaterial for what follows; only $\theta\in\{\pm1\}$ is used. (This $\theta$ is the incidence sign of~\cite[Definition~2.2]{BaKe}.)

Let $M_{n-1}$ and $M_{n}$ denote multiplication by the weights $m_{n-1}$ and $m_n$ in these bases. Since $\delta^{n}$ is the adjoint of $d^{n}$ for the weighted inner products, its matrix is $M_{n-1}^{-1}\mathsf B^{\mathsf T}M_{n}$, and therefore
\[
\Delta_{n}=d^{n}\delta^{n}\ \text{ has matrix }\ \mathsf B\,M_{n-1}^{-1}\,\mathsf B^{\mathsf T}\,M_{n}.
\]
(This matrix is not symmetric, the bases being orthogonal but not orthonormal; what is symmetric, as it must be, is $M_{n}\Delta_{n}$, which equals $M_n\mathsf BM_{n-1}^{-1}\mathsf B^{\mathsf T}M_n$.) Two distinct $n$-simplices have at most one common $(n-1)$-face, so the product contains at most one term and we read off, with $\tau=\sigma\cap\sigma'$,
\[
\Delta_{n}[\sigma,\sigma']\;=\;\theta(\tau,\sigma)\,\theta(\tau,\sigma')\,\frac{m_n(\sigma')}{m_{n-1}(\tau)}\;=\;\varepsilon(\sigma,\sigma')\,a(\sigma,\sigma'),\qquad
\varepsilon(\sigma,\sigma'):=\theta(\tau,\sigma)\,\theta(\tau,\sigma').
\]
The symmetry $\varepsilon(\sigma,\sigma')=\varepsilon(\sigma',\sigma)$ is now immediate: the right-hand side is a product of two numbers attached to the pair $(\tau,\sigma)$ and $(\tau,\sigma')$, in which $\sigma$ and $\sigma'$ play interchangeable roles. Likewise the diagonal entry is $\sum_{\tau\prec\sigma}\theta(\tau,\sigma)^2\,m_n(\sigma)/m_{n-1}(\tau)=V_{n}(\sigma)$, since $\theta^2=1$; this recovers the computation above.

Hence, in the skew convention,
\begin{equation}\label{eq:Delta-n-local-skew}
(\Delta_{n,\mathrm{skew}}u)(\sigma)\;=\;V_{n}(\sigma)\,u(\sigma)\;+\;\sum_{\sigma'\sim\sigma}\widetilde a(\sigma,\sigma')\,u(\sigma'),\qquad
\widetilde a(\sigma,\sigma')=\varepsilon(\sigma,\sigma')\,a(\sigma,\sigma').
\end{equation}
Note that the off-diagonal sum carries the alternating signs $\varepsilon(\sigma,\sigma')\in\{\pm1\}$, but the global sign in front of the sum is $+$ (not $-$): the orientation-induced minus signs are absorbed into $\widetilde a$ via $\varepsilon$, and the diagonal still adds $+V_{n}(\sigma)u(\sigma)$, consistent with $\Delta_{n,\mathrm{skew}}\ge 0$.

\smallskip\emph{The canonical basis: the sign in closed form.}
So far the representatives were an arbitrary fixed choice, and $\theta$ --- hence $\varepsilon$ --- depends on it. Take now for $\sigma_{\mathrm{ref}}$ the color-ordered representative of Section~\ref{sec:color}, i.e.\ $\sigma_{\mathrm{ref}}=\sigma_{\mathrm{can}}$. Deleting a vertex from a color-ordered tuple leaves a color-ordered tuple, so the $i$-th face of $\sigma_{\mathrm{can}}$ \emph{is} the canonical representative $\tau_{\mathrm{can}}$, with no permutation to correct: $\theta(\tau,\sigma)=(-1)^{i}$ exactly. Therefore
\begin{equation}\label{eq:Delta-n-local-sym}
(\Delta^{\mathrm{can}}_{n,\mathrm{sym}}u)(\sigma)\;=\;V_{n}(\sigma)\,u(\sigma)+\sum_{\sigma'\sim\sigma}(-1)^{i+j}\,a(\sigma,\sigma')\,u(\sigma'),
\end{equation}
where $i$ and $j$ are now the \emph{color ranks}, in $\sigma$ and in $\sigma'$ respectively, of the two vertices not lying in $\tau$. The coloring therefore does not remove the sign; what it does is make it canonical --- determined by the coloring alone, with no choice of orientations left --- and computable. It does vanish in one case of interest, discussed in Remark~\ref{rem:signs-bipartite}.

\smallskip\emph{Symmetrization via $U$.}
The unitary $U$ defined in Step 1 symmetrizes the off-diagonal coefficients: when we conjugate $\Delta_{n}$ by $U$ (and read the result on the unoriented line-complex $\widehat\Vc$), each off-diagonal term $a(\sigma,\sigma')=m_n(\sigma')/m_{n-1}(\tau)$ acquires a factor $m_n(|\sigma|)^{1/2}/m_n(|\sigma'|)^{1/2}$ from the change of basis, becoming $\sqrt{m_n(|\sigma|)m_n(|\sigma'|)}/m_{n-1}(|\tau|)$. Setting $g=Uu$, we obtain
\[
(U\Delta_{n}U^{-1}g)(\sigma)
=V_{n}(\sigma)\,g(\sigma)\;+\;\sum_{\sigma'\sim\sigma}\widetilde w(\sigma,\sigma')\,g(\sigma'),
\]
where the off-diagonal kernel
\[
\widetilde w(\sigma,\sigma'):=\varepsilon(\sigma,\sigma')\,\frac{\sqrt{m_n(\sigma)m_n(\sigma')}}{m_{n-1}(\sigma\cap\sigma')}
\]
is symmetric in $\sigma,\sigma'$: the modulus is manifestly symmetric, and the sign $\varepsilon$ was just shown to be.

Summarizing: there exist signs $\eta(\sigma,\sigma')\in\{\pm 1\}$ with $\eta(\sigma,\sigma')=\eta(\sigma',\sigma)$ and a non-negative kernel $w(\sigma,\sigma')=\sqrt{m_n(\sigma)m_n(\sigma')}/m_{n-1}(\sigma\cap\sigma')\,\mathbf 1_{\{\sigma\sim\sigma'\}}$ such that
\[
U\,\Delta_{n}\,U^{-1}\;=\;\mathcal Q(V_{n})+\mathcal A_\eta,
\]
where $\mathcal A_\eta$ has signed kernel $w_\eta(\sigma,\sigma'):=\eta(\sigma,\sigma')\,w(\sigma,\sigma')$, so that $|w_\eta(\sigma,\sigma')|=w(\sigma,\sigma')$. By Proposition~\ref{prop:schur}, $\|\mathcal A_\eta\|\le\|D\|_\infty$ where $D(\sigma)=\sum_{\sigma'}w(\sigma,\sigma')$.

\smallskip\noindent\emph{Step 3. Operator norm estimates.}
Assumption (A2) gives $\|\mathcal Q(V_{n})\|=\|V_{n}\|_\infty<\infty$. Hence
\[
\|U\Delta_{n}U^{-1}\|\;\le\;\|\mathcal A_\eta\|+\|\mathcal Q(V_{n})\|\;\le\;\|D\|_\infty+\|V_{n}\|_\infty.
\]

\smallskip\noindent\emph{Step 4. Self-adjoint extension.}
For $u\in\Cc^n_c$ the cochain $\Delta_nu$ is defined pointwise without any finiteness hypothesis: $\delta^nu$ is a finite sum, since $u$ has finite support, and $d^n$ then contributes $n+1$ terms. Note, however, that $\Delta_nu$ need \emph{not} have finite support: while $\mathrm{supp}(\delta^nu)$ consists of faces of $\mathrm{supp}\,u$ and is finite, $\mathrm{supp}(\Delta_nu)$ contains every $n$-simplex through one of those faces, and there may be infinitely many of them when the complex is not locally finite. Square summability is therefore not automatic, and it is here that Steps~1--3 are used: the identity $U\Delta_nu=(\mathcal A_\eta+\mathcal Q(V_n))Uu$ holds pointwise, all sums being finite because $Uu$ has finite support, and the right-hand side is bounded on $\ell^2(\widehat{\Vc})$ by Proposition~\ref{prop:schur}, which requires no local finiteness. Hence $\Delta_nu\in\ell^2(m_n)$, with the stated bound. Moreover $\Delta_{n}=d^n\delta^n$ is symmetric on the core $\Cc^n_{c}$, $\delta^n$ being the formal adjoint of $d^n$ (Lemma~\ref{lem:2-3-1}). By Steps 1--3, $\|\Delta_{n}u\|_{\ell^2(m_n)}\le(\|D\|_\infty+\|V_{n}\|_\infty)\|u\|_{\ell^2(m_n)}$ for $u\in\Cc^n_{c}$. By the boundedness criterion of Remark~\ref{rem:HT}, $\Delta_{n}$ extends uniquely to a bounded self-adjoint operator on $\ell^2(m_n)$; in particular, it is essentially self-adjoint on $\Cc^n_{c}$.
\end{proof}

\smallskip
\noindent\emph{On the local finiteness hypothesis.} It is worth recording that local finiteness is not used anywhere in the proof, and may be omitted from the statement. For $u\in\Cc^n_c$ the cochain $\Delta_nu$ is defined pointwise without any finiteness assumption, as noted in Step~4; its support may well be infinite, but the identity $U\Delta_nu=(\mathcal A_\eta+\mathcal Q(V_n))Uu$ holds pointwise and the right-hand side is bounded by Proposition~\ref{prop:schur}, which requires no local finiteness either. Assumptions (A1)--(A2) therefore carry the whole argument. We nevertheless keep the hypothesis in the statement, since every complex to which the theorem is applied below is locally finite, and since it makes the line-complex of Definition~\ref{def:line-complex} locally finite as well.

\begin{corollary}\label{cor:kn-line-degree}
If each $n$-simplex has at most $L<\infty$ neighbors across its $(n{-}1)$-faces, and the weights satisfy
\(
0<m_{-}\le m_{n-1}(\tau)\) and \(m_n(\sigma)\le m^{+}<\infty,
\)
then
\[
\|D\|_\infty\le L\,\tfrac{m^+}{m_-},\qquad
\|V_{n}\|_\infty\le (n+1)\tfrac{m^+}{m_-},
\]
so $\Delta_{n}$ extends uniquely to a bounded self-adjoint operator on $\ell^2(m_n)$.
\end{corollary}

\begin{remark}[Relation to~\cite{BaKe}]\label{rem:BaKe-thm63}
The \emph{representation} underlying Theorem~\ref{thm:kn-esa} --- an adjacency operator with a signed kernel plus a diagonal potential --- is the specialization to the top-degree block $\Delta_{n}$ of the general signed Schr\"odinger-operator representation of Bartmann--Keller~\cite[Section~3.4]{BaKe}, discussed in Section~\ref{sec:color}. What the theorem adds to that representation is the identification of the underlying graph as the line-complex of $\Sc_n$, the closed form of the kernel and of the potential in terms of $(m_{n-1},m_n)$, and the resulting quantitative bound. The adjacency kernel $w(\sigma,\sigma')$ and the potential $V_{n}(\sigma)$ of the theorem are the explicit line-complex incarnations of the abstract kernel $b^{\circ}$ and potential $c^{\circ}$ of~\cite{BaKe}. We also note that~\cite[Theorem~4.1]{BaKe} addresses boundedness of the Hodge Laplacian, giving an equivalence for finite-dimensional complexes such as ours and a sufficient criterion in general.

It should be said plainly that the two boundedness criteria are qualitatively equivalent in our setting. Theorem~\ref{thm:kn-esa}, and more generally Lemma~\ref{lem:form-bounds}, gives a criterion \emph{degree by degree}, for each block $\Delta_{k}$ separately, whereas~\cite[Theorem~4.1]{BaKe} characterizes boundedness of the \emph{full} direct sum $\bigoplus_{k}\Delta_{k}$ across all degrees at once; here $L$ is the operator of Section~\ref{section3}, the letter $L$ in Corollary~\ref{cor:kn-line-degree} above being an unrelated combinatorial bound. Since $n$ is finite here and the norm of a finite orthogonal sum is the maximum of the norms, the two criteria in fact coincide qualitatively; the contribution of the present paper on this point is the explicit constant, and its closed form at top degree, not the boundedness criterion.
\end{remark}

\begin{remark}[Signs, bipartiteness, and sharpness]\label{rem:signs-bipartite}
In the canonical basis the signs are $\eta(\sigma,\sigma')=(-1)^{i+j}$, with $i,j$ the color ranks of the two unshared vertices~\eqref{eq:Delta-n-local-sym}. At top degree $n=1$ on a $d$-regular graph, three facts of this paper then turn out to be one. The signs become $\equiv+1$, for the canonical $2$-coloring, exactly when the graph is bipartite; for a general ordered coloring the coloring-free form of the statement is that the signing $\eta$ is \emph{balanced}.

In one direction, with a $2$-coloring the two endpoints not shared by adjacent edges carry the same color, so $i=j$ and $\eta=(-1)^{i+j}=+1$. In the other, the computation below shows that a non-bipartite graph cannot have $\eta\equiv+1$; it is independent of the coloring: changing the coloring replaces $\eta(e,e')$ by $s(e)\eta(e,e')s(e')$ for some $s$ with values in $\{\pm1\}$ (a \emph{switching}), which leaves the product of the signs along any cycle unchanged, each $s(e)$ occurring twice and $s(e)^2=1$. For an edge $e$ with canonical representative $[u,v]$, $c(u)<c(v)$, the incidence signs are $\theta(\{v\},e)=+1$ and $\theta(\{u\},e)=-1$, so $\theta(\{u\},e)\,\theta(\{v\},e)=-1$. Hence, for a cycle $x_0x_1\cdots x_{\ell-1}x_0$ of $G$ with edges $e_i=\{x_i,x_{i+1}\}$, regrouping the product $\prod_{i}\eta(e_i,e_{i+1})$ edge by edge gives
\[
\prod_{i=0}^{\ell-1}\eta(e_i,e_{i+1})=\prod_{j=0}^{\ell-1}\theta(\{x_j\},e_j)\,\theta(\{x_{j+1}\},e_j)=(-1)^{\ell}.
\]
So the sign product around the corresponding cycle of $L(G)$ is $-1$ on every odd cycle of $G$; $\eta\equiv+1$ therefore forces all cycles of $G$ to be even, i.e.\ $G$ bipartite. The remaining generators of the cycle space of $L(G)$ are the triangles formed by three edges of $G$ through a common vertex $x$ --- the cycle space of $L(G)$ being generated by these triangles together with the cycles inherited from $G$ --- and on those the product is $\bigl(\theta(\{x\},e_1)\theta(\{x\},e_2)\theta(\{x\},e_3)\bigr)^2=+1$; balancedness of $\eta$ is therefore decided by the cycles of $G$ alone. In the language of signed graphs: $\eta$ is a balanced signing of $L(G)$ precisely when $G$ is bipartite.

The same bipartition governs the attainment of the bound $\|\widetilde\Delta_1\|\le2d$, through the spectral criterion of Proposition~\ref{prop:sharpness-bipartite}. The link is one-directional and conditional: a bipartition makes the two endpoints not shared by adjacent edges carry the same color, and it puts $2d$ in the spectrum of $\Delta_0$ \emph{provided the graph is also amenable}. Neither implication is an equivalence --- the $d$-regular tree is bipartite without attaining $2d$ (Remark~\ref{rem:tree-counterexample}), and a non-bipartite graph may attain it (Section~\ref{sharp}). The exact criterion is $\inf\mathrm{spec}(Q_G)=0$.

None of the \emph{estimates} depends on any of this. Proposition~\ref{prop:schur} bounds $\mathcal A_\eta$ through $|w|$ alone, so the signs are irrelevant to every bound proved here, bipartite or not.
\end{remark}

\paragraph{Example ($n=3$, unit weights).}
On a uniform tetrahedral mesh with unit weights $m_{n-1}\equiv m_n\equiv 1$, each interior tetrahedron has four neighbors. In the skew convention, the local formula~\eqref{eq:Delta-n-local-skew} specializes to
\begin{equation}\label{eq:tetra-top}
(\Delta_{3,\mathrm{skew}}u)(\sigma)\;=\;4\,u(\sigma)\ +\ \sum_{\sigma'\sim \sigma}\widetilde a(\sigma,\sigma')\,u(\sigma'),
\end{equation}
where $\widetilde a(\sigma,\sigma')\in\{\pm 1\}$ encode the relative orientation signs (so the off-diagonal sum is signed, even though the global sign in front is $+$). In the canonical basis of Section~\ref{sec:color},
\begin{equation}\label{eq:tetra-top-sym}
(\Delta^{\mathrm{can}}_{3,\mathrm{sym}}u)(\sigma)\;=\;4\,u(\sigma)\ +\ \sum_{\sigma'\sim \sigma}(-1)^{i+j}\,u(\sigma'),
\end{equation}
in line with~\eqref{eq:Delta-n-local-sym}. The contrast with~\eqref{eq:tetra-top} is not that the signs disappear. They do not: the relevant notion here is the balancedness of the signing $\eta$ on the line-complex (Remark~\ref{rem:signs-bipartite}), and there is no reason for it to hold on a tetrahedral mesh: by that remark balancedness is decided by the cycles of the underlying adjacency structure, which here is not bipartite. The contrast is that the signs cease to be a matter of choice. In~\eqref{eq:tetra-top} they depend on the orientations one picks; in~\eqref{eq:tetra-top-sym} they are read off from the colors.

\paragraph{Comparison with intermediate degrees.}
For every $k\in\{0,\dots,n\}$, Lemma~\ref{lem:form-bounds} gives the Schur-type bound
\[
\|\widetilde\Delta_{k}\|\ \le\ D_k^{\downarrow}+D_k^{\uparrow}.
\]
At top degree $k=n$ the down term is the only one present, and Theorem~\ref{thm:kn-esa} reinterprets it in terms of the adjacency structure of the line-complex. In the unweighted case the reinterpretation can only help, and often does not: one has
\[
\|D\|_\infty+\|V_n\|_\infty\;=\;\sup_{\sigma}\ \sum_{\tau\prec\sigma}\#F_{\tau}\ \le\ (n+1)\sup_{\tau}\#F_{\tau}\;=\;D_n^{\downarrow},
\]
with equality precisely when some top simplex has all of its $n+1$ codimension-one faces of maximal coface count --- or, on an infinite complex where the suprema need not be attained, when a sequence of top simplices does so asymptotically. That condition holds on the homogeneous complexes of interest here. At $n=1$ on an unweighted $d$-regular graph, for instance, both routes give exactly $2d$, since $D_1^{\downarrow}=2\sup_x d_{\Vc}(x)=2d$ and $\Delta(L(G))+2=2d$. It is not automatic, however: for two triangles glued along an edge ($n=2$, unit weights) the left-hand side is $1+3=4$, which is the exact norm, against $D_2^{\downarrow}=6$. What the line-complex form does provide is a \emph{structural} decomposition --- adjacency plus potential on an explicit graph. That is what makes the Floquet--Bloch analysis and the sharpness discussion of Section~\ref{sharp} possible; the sharp constants themselves come from that analysis, not from the reduction.

The two routes do differ in two respects. First, the top-degree bound $\|D\|_\infty+\|V_n\|_\infty$ is attained on a non-trivial class of complexes, which is not known for the intermediate-degree estimate (Remark~\ref{rem:form-bounds-sharpness}); it is of course not attained on all of them, the triangular lattice of Table~\ref{tab:lattices-combined} giving $9$ against a bound of $12$. Already at $n=1$ on an unweighted $d$-regular bipartite amenable graph it equals $2(d-1)+2=2d$, the exact norm by Proposition~\ref{prop:sharpness-bipartite}; and it is attained on the two glued triangles above.

Second --- and this deserves to be stated plainly, since the displayed inequality of the previous paragraph might suggest otherwise --- on a \emph{weighted} complex the two bounds are \emph{incomparable}: neither dominates. Both routes use only the down block at $k=n$, so the discrepancy has nothing to do with the Cauchy--Schwarz step on the $k+2$ terms of $d^{k+1}$, which is absent here. It has two independent sources. The first is the counting just described: $D_n^{\downarrow}$ replaces the sum over the $n+1$ faces of a given top simplex by $(n+1)$ times a single supremum. The second is the symmetrization performed by the conjugation of Step~1, which turns the coefficient $m_n(\sigma')/m_{n-1}(\tau)$ into the geometric mean $\sqrt{m_n(\sigma)m_n(\sigma')}/m_{n-1}(\tau)$; this is a gain when $m_n(\sigma)\le m_n(\sigma')$ and a loss otherwise. The two effects pull in opposite directions, and either can win.

\begin{example}[Neither bound dominates]\label{ex:two-routes}
Let $n=2$, let $\sigma_1=\{x_0,x_1,x_2\}$ and $\sigma_2=\{x_0,x_1,x_3\}$ be two triangles glued along $\tau_0=\{x_0,x_1\}$, and take $m_1\equiv1$, $m_2(\sigma_1)=M$, $m_2(\sigma_2)=\varepsilon$ with $M\ge\varepsilon>0$. Then
\[
D_2^{\downarrow}=3\,(M+\varepsilon),\qquad
\|D\|_\infty+\|V_2\|_\infty=3M+\sqrt{M\varepsilon},
\]
so the line-complex bound is the better one exactly when $M<9\varepsilon$. For $M=\varepsilon=1$ it gives $4$, which is the exact norm, against $D_2^{\downarrow}=6$; for $M=1$, $\varepsilon=1/100$ it gives $3.1$ against $D_2^{\downarrow}=3.03$, and is then the worse of the two (the exact norm being $\simeq3.003$). In practice both should be computed and the smaller retained.
\end{example}

\noindent For $0<k<n$ both the up and the down components contribute, and no such reduction is available.

\section{Sharp edge-Laplacian bounds without geometric assumptions}\label{sharp}

\paragraph{Amenability: terminology.} A countable graph $G$ is \emph{amenable} if it admits a \emph{F\o lner sequence}, i.e.\ a sequence $(\Omega_R)$ of finite vertex sets with $|\partial\Omega_R|/|\Omega_R|\to0$, where $\partial\Omega_R:=\{x\in\Omega_R:\ x\sim y\text{ for some }y\notin\Omega_R\}$ is the inner vertex boundary. Equivalently, the isoperimetric (Cheeger) constant of $G$ vanishes. For a connected $d$-regular graph, the discrete Cheeger inequalities give the spectral form of the criterion used below: $\sup\mathrm{spec}(A_G)=d$, i.e.\ $\inf\mathrm{spec}(\Delta_0)=0$, if and only if $G$ is amenable. They are usually stated for the \emph{edge} isoperimetric constant $h_E$, in the form $h_E^2/(2d)\le\inf\mathrm{spec}(\Delta_0)\le h_E$; since $h_V\le h_E\le d\,h_V$, the vertex constant $h_V$ used above vanishes exactly when $h_E$ does, which is all we need. See Dodziuk~\cite{Dod} and Mohar~\cite{Moh}. On a Cayley graph this is Kesten's criterion~\cite{Kesten}; for a general $d$-regular graph the equivalence is the Cheeger-inequality statement, and it is in that form that we use it.

In this section we focus on the edge (top-degree) Hodge Laplacian $\Delta_{1}$ on a locally finite (possibly weighted) graph $G$ with finite maximum degree. Here $n=1$, so the complex is the graph itself and $\Delta_1=d^1\delta^1$, the up-term vanishing by the convention $d^{2}=0$. This deserves emphasis, because several lattices below --- the triangular, Kagome and face-centered cubic ones --- are full of triangles, so that their clique complexes do carry $2$-simplices. Choosing $n=1$ truncates at dimension one and studies the edge block of the truncation. The Hodge Laplacian of the untruncated clique complex would carry the additional up-term $\delta^{2}d^{2}$; it is a different operator, outside the scope of this section.

The bound $\|\widetilde\Delta_{1}\|\le 2d$ for unweighted $d$-regular graphs ($d\ge 2$) is derived in Proposition~\ref{prop:top-bounded-ESA}. Its weighted version, with the constant $C_1/c_0$, is Proposition~\ref{prop:weighted-extension}.

The bound $2d$ is attained on every \emph{finite} unweighted $d$-regular bipartite graph, and on every infinite one that is amenable (Proposition~\ref{prop:sharpness-bipartite} below). The interest of the Floquet--Bloch computations therefore lies elsewhere: they give the \emph{exact} norm on the non-bipartite lattices (triangular, FCC, Kagome), where the bound is strict, and confirm the value $2d$ on the bipartite ones ($\Z^2,\Z^3,$ BCC, $\Z^4$, diamond).

We first record the elementary identity which underlies every exact value in this section, and which is used repeatedly below.

\begin{lemma}[Edge/vertex isospectrality]\label{lem:iso-01}
Let $G$ be a locally finite weighted graph with $\sup_x d_{\Vc}(x)<\infty$, viewed as a $1$-complex ($n=1$). Write $T:=d^1$, a bounded operator from $\ell^2(m_0)$ to $\ell^2(m_1)$, with $\delta^1=T^\ast$. Then $\Delta_1=TT^\ast$ and $\Delta_0=T^\ast T$, so that
\[
\mathrm{spec}(\Delta_1)\setminus\{0\}=\mathrm{spec}(\Delta_0)\setminus\{0\},
\qquad
\|\Delta_{1}\|=\|T\|^2=\|\Delta_0\|.
\]
Moreover $\|\widetilde\Delta_{1}\|=\|\Delta_{1}\|$ and $\mathrm{spec}(\widetilde\Delta_1)=\mathrm{spec}(\Delta_1)$, with \emph{no} ellipticity assumption on the weights: as observed after~\eqref{eq:norm-blocks}, $M_1^{-1/2}$ is a \emph{unitary} from the unweighted $\ell^2$ onto $\ell^2(m_1)$, and $\widetilde\Delta_1$ is by definition the transport of $\Delta_1$ through it. (In the unweighted case $M_1=I$ and the two operators simply coincide.)
\end{lemma}

\begin{proof}
Boundedness of $T$ is Lemma~\ref{lem:form-bounds} in degree $0$ ($\|d^1u\|_1^2=\langle\mathbf L_0^+u,u\rangle_0\le D_0^{\uparrow}\|u\|_0^2$), and $\delta^1=T^\ast$ on the core by Proposition~\ref{prop:2-5-1}, hence everywhere by continuity. The equality of non-zero spectra of $TT^\ast$ and $T^\ast T$ and the identity $\|TT^\ast\|=\|T\|^2=\|T^\ast T\|$ are standard.
\end{proof}

\subsection{Floquet--Bloch symbols for the edge block on periodic lattices}\label{app:bloch}
\paragraph{Setup and normalization}
Let $G=(\Vc,\Ec)$ be a periodic unweighted graph, invariant under a lattice $\Gamma\simeq\Z^r$ ($r$ the rank of the period lattice) acting freely on $\Vc$. Choose a fundamental cell containing finitely many vertices and finitely many edge orbits; let $N_1$ denote the number of edge orbits. The Fourier transform identifies $\ell^2(\Ec)$ with a direct integral over the Brillouin zone $\mathbb T^r$:
\[
\ell^2(\Ec)\ \simeq\ \int_{\mathbb T^r}^{\oplus}\C^{N_1}\,d\theta.
\]
Under this reduction, $\widetilde\Delta_{1}$ is unitarily equivalent to multiplication by a Hermitian $N_1\times N_1$ matrix-valued symbol $\hat\sigma(\theta)$ (here we write $\hat\sigma$ to distinguish it from a simplex), and
\[
\|\widetilde\Delta_{1}\|\ =\ \sup_{\theta\in\mathbb T^r}\ \rho(\hat\sigma(\theta)),
\]
where $\rho$ denotes the spectral radius; the symbols being Hermitian and continuous in $\theta$, the supremum is attained. In practice we do not compute the $N_1\times N_1$ edge symbol but the symbol of the vertex Laplacian $\Delta_0$, which has one row per vertex orbit and is therefore much smaller, and we pass to the edge block through Lemma~\ref{lem:iso-01}. The two are consistent: at each $\theta$ the edge symbol has rank at most the number of vertex orbits, and its non-zero spectrum is that of the vertex symbol. The rank drops at $\theta=0$, where the vertex symbol vanishes on the constants.

Throughout this section, and only here, the symbol $d$ denotes the vertex degree of the lattice --- neither the rank $r$ of the period lattice, nor the coboundary operator $d$ of Sections~\ref{sec:prelim-comb}--\ref{sec:esa}, which does not appear below.

\paragraph{Examples}
\begin{itemize}[leftmargin=2em]
  \item \textbf{Square lattice $\Z^2$ ($d=4$).}
  The fundamental cell contains $1$ vertex and $2$ orbits of edges. The Bloch symbol is a $2\times 2$ Hermitian matrix whose spectrum is $\{0,\ 4-2\cos\theta_1-2\cos\theta_2\}$, with maximum $8$ attained at $\theta=(\pi,\pi)$. Hence $\|\widetilde\Delta_{1}\|=8=2d$ (vs.\ the bound $2d=8$, here attained).

  \item \textbf{Triangular lattice ($d=6$, non-bipartite).}
  The fundamental cell contains $1$ vertex and $3$ orbits of edges. With primitive directions $a_1,a_2,a_1+a_2$, the largest eigenvalue of the Bloch symbol is
  \[
  \lambda_{\max}(\theta)\;=\;6-2\bigl(\cos\theta_1+\cos\theta_2+\cos(\theta_1+\theta_2)\bigr),
  \]
  maximized at the Brillouin-zone corner $\theta=(2\pi/3,2\pi/3)$ (the $K$ point), where each of the three cosines equals $-1/2$, giving $\|\widetilde\Delta_{1}\|=6-2(-3/2)=9$ (vs.\ the bound $2d=12$). The value $2d=12$ is \emph{not} attained, the exact norm being $9$. Note that the strictness is established here by the Bloch computation itself, and not deduced from non-bipartiteness. By Proposition~\ref{prop:sharpness-bipartite}, the attainment of $2d$ is governed by the signless Laplacian condition $\inf\mathrm{spec}(dI+A_G)=0$. For a \emph{finite} connected graph that condition is equivalent to bipartiteness; on an infinite graph it may hold on a non-bipartite graph as well, it being enough that the odd cycles be confined to a region a F\o lner sequence can avoid.

  \item \textbf{Cubic lattice $\Z^3$ ($d=6$).}
  The fundamental cell contains $1$ vertex and $3$ orbits of edges. The Bloch symbol is a $3\times 3$ Hermitian matrix whose largest eigenvalue is $\lambda_{\max}(\theta)=2(3-\cos\theta_1-\cos\theta_2-\cos\theta_3)$, with maximum at $\theta=(\pi,\pi,\pi)$. Hence $\|\widetilde\Delta_{1}\|=12=2d$ (vs.\ the bound $2d=12$, here attained).

  \item \textbf{BCC lattice ($d=8$, bipartite).} Geometrically, the body-centered cubic lattice has vertex set $\Z^3$ together with one additional vertex at the center of each unit cube, two vertices being adjacent when they are at minimal distance $\sqrt3/2$. Equivalently, each center vertex is joined to the $8$ corners of its cube and to nothing else; the cubic bonds of $\Z^3$ are \emph{not} retained. The graph is therefore $8$-regular and bipartite, the two classes being the corners and the centers.
  In the conventional cubic cell ($2$ vertices per cell, hence $8$ edge orbits, one for each of the eight center-to-corner bonds), the Bloch symbol of the vertex Laplacian $\Delta_0$ is a $2\times 2$ Hermitian matrix whose eigenvalues are $8\pm|H_{AB}(\theta)|$, with $H_{AB}(\theta)=8\cos(\theta_1/2)\cos(\theta_2/2)\cos(\theta_3/2)$ the cross-term between the two sublattices. Its maximum modulus $|H_{AB}|=8$ is attained at $\theta=0$, giving $\|\Delta_0\|=8+8=16$. Since $\|\widetilde\Delta_{1}\|=\|\Delta_0\|$ (Lemma~\ref{lem:iso-01}), this yields $\|\widetilde\Delta_{1}\|=16=2d$ (vs.\ the bound $2d=16$, here attained). The value $16$ is achieved because BCC is bipartite (the two sublattices A and B); the harmonic mode $f(x)=(-1)^{\mathrm{sublattice}(x)}$ realizes the value $2d$.

  \item \textbf{FCC lattice ($d=12$, non-bipartite).} Geometrically, the face-centered cubic lattice has vertex set $\Z^3$ together with one additional vertex at the center of each face of each unit cube, two vertices being adjacent when they are at (minimal) distance $1/\sqrt2$; the cubic bonds of $\Z^3$ are again \emph{not} retained. Each face center is thus adjacent to the $4$ corners of its own face and to the $8$ neighboring face centers, and each corner to the $12$ face centers surrounding it, so that the graph is $12$-regular.
  In the primitive cell ($1$ vertex, $6$ edge orbits with primitive vectors $a_1,a_2,a_3,a_1{-}a_2,a_2{-}a_3,a_3{-}a_1$), the edge symbol of $\Delta_1=d^1\delta^1$ is a $6\times6$ Hermitian matrix of rank at most one, so that it has at most one non-zero eigenvalue, which by Lemma~\ref{lem:iso-01} is the scalar symbol of $\Delta_0$ (the rank is zero exactly at $\theta=0$, where that symbol vanishes):
  \[
  \lambda_{\max}(\theta)\;=\;12-A(\theta),\qquad
  A(\theta):=2\bigl(\cos\theta_1+\cos\theta_2+\cos\theta_3+\cos(\theta_1{-}\theta_2)+\cos(\theta_2{-}\theta_3)+\cos(\theta_3{-}\theta_1)\bigr),
  \]
  where $A$ is the adjacency symbol. The minimum is $\min_\theta A=-4$, rather than the value $-12$ that all six cosines being $-1$ simultaneously would give, which is impossible. This is certified by an exact algebraic identity: setting $z_0:=1$ and $z_j:=e^{i\theta_j}$ for $j=1,2,3$, the bracket equals $\sum_{0\le j<k\le3}\Re(z_j\overline{z_k})=\tfrac12\bigl(\bigl|\textstyle\sum_{j=0}^3 z_j\bigr|^2-4\bigr)\ge-2$, with equality exactly when $\sum_{j=0}^3z_j=0$. Such $\theta$ exist, e.g.\ $\theta=(\pi,\pi,0)$, where the first three cosines sum to $-1$ and the three differences to $-1$. Hence
  \[
  \|\widetilde\Delta_{1}\|\;=\;12-(-4)\;=\;16\quad\text{(rigorously certified)},
  \]
  attained for instance at $\theta=(\pi,\pi,0)$. The bound $2d=24$ is \emph{not} attained, the exact norm being $16$; as for the triangular lattice, this is established by the Bloch computation and not inferred from non-bipartiteness. We stress that the value is not read off a finite approximation: on a torus quotient FCC/$(N\Z)^3$ the Bloch phases are constrained to $\tfrac{2\pi}{N}\Z^3$, and the condition $\sum_{j=0}^3 z_j=0$ forces $N$ to be even, so the supremum $16$ is missed by every odd quotient. The algebraic identity above is what makes the value exact.
\end{itemize}

\emph{Geometric description of the Kagome lattice, and the last two entries of the table.} The Kagome lattice appearing in Table~\ref{tab:lattices-combined} is the planar lattice of corner-sharing triangles; equivalently, it can be visualized as a hexagonal tiling with a triangle attached across each edge of every hexagon. It has vertex degree $d=4$ and is non-bipartite (each triangle is a $3$-cycle). In the primitive cell with three sites $A,B,C$, its adjacency symbol is the $3\times3$ Hermitian matrix with off-diagonal entries $1+e^{-i\theta_1}$, $1+e^{-i\theta_2}$ and $1+e^{-i(\theta_2-\theta_1)}$; it carries a flat band at $-2$ and two dispersive bands lying above it, the lower of which touches $-2$ at $\theta=0$, so that $\min_\theta\mathrm{spec}=-2$ and, by Lemma~\ref{lem:iso-01}, $\|\widetilde\Delta_{1}\|=4-(-2)=6$.

The two remaining entries are obtained in the same way. The diamond lattice ($d=4$) has two sites per cell and adjacency symbol with off-diagonal entry $f(\theta)=1+e^{-i\theta_1}+e^{-i\theta_2}+e^{-i\theta_3}$, hence spectrum $\pm|f(\theta)|$ with $\max_\theta|f|=4$ attained at $\theta=0$; the hypercubic lattice $\Z^4$ ($d=8$) has scalar symbol $2\sum_{j=1}^{4}\cos\theta_j$, of minimum $-8$ at $\theta=(\pi,\pi,\pi,\pi)$. Both are bipartite and amenable, so Proposition~\ref{prop:sharpness-bipartite} applies as well and yields $8$ and $16$ respectively, in agreement with these Bloch computations.

\paragraph{Exactly when the bound $2d$ is attained.} The criterion is not bipartiteness itself but the vanishing of the spectral infimum of the signless Laplacian, the two coinciding for a finite connected graph and not in general (Remark~\ref{rem:tree-counterexample}). Throughout this paragraph the graph is \emph{unweighted}, $m_0\equiv m_1\equiv1$: this is where the value $2d$ comes from, the corresponding constant for a general vertex weight being $2\sup_x d_{\Vc}(x)$ (Section~\ref{sec:normalization}).

\begin{proposition}[Exact sharpness of the bound $2d$]\label{prop:sharpness-bipartite}
Let $G$ be a connected locally finite unweighted $d$-regular graph ($d\ge2$), and let $Q_G:=dI+A_G$ denote its signless Laplacian. Then $\|\widetilde\Delta_{1}\|=\|\Delta_0\|=d-\inf\mathrm{spec}(A_G)\le 2d$, and the following are equivalent:
\begin{enumerate}[leftmargin=2em,label=(\roman*)]
\item $\|\widetilde\Delta_{1}\|=2d$;
\item $2d\in\mathrm{spec}(\Delta_0)$;
\item $\inf\mathrm{spec}(Q_G)=0$, i.e.\ $\displaystyle\inf_{\|u\|_{\ell^2}=1}\ \sum_{\{x,y\}\in\Ec}|u(x)+u(y)|^2=0$.
\end{enumerate}
Each of these implies $\sup\mathrm{spec}(A_G)=d$, hence the amenability of $G$ by the Cheeger criterion recalled above. Conversely, they all hold as soon as $G$ is bipartite \emph{and} amenable --- in particular for every finite bipartite $d$-regular graph and for every bipartite periodic lattice, such as $\Z^2$, $\Z^3$, $\Z^4$, BCC and diamond. For $G$ finite and connected, (i)--(iii) hold if and only if $G$ is bipartite.
\end{proposition}

\begin{proof}
The first identity is Lemma~\ref{lem:iso-01} together with $\Delta_0=dI-A_G$. Since $\mathrm{spec}(\Delta_0)\subset[0,2d]$ is closed, $\|\Delta_0\|=2d$ if and only if $2d\in\mathrm{spec}(\Delta_0)$, which gives (i)$\Leftrightarrow$(ii). Next, $2d\in\mathrm{spec}(\Delta_0)$ means $-d\in\mathrm{spec}(A_G)$, i.e.\ $\inf\mathrm{spec}(Q_G)=0$ with $Q_G=dI+A_G\ge0$, and the quadratic form of $Q_G$ is
\[
\langle Q_Gu,u\rangle=\sum_{x\in\Vc}d\,|u(x)|^2+\sum_{x\in\Vc}\sum_{y\sim x}u(y)\overline{u(x)}
=\sum_{\{x,y\}\in\Ec}|u(x)+u(y)|^2 ,
\]
whence (ii)$\Leftrightarrow$(iii). Moreover, each of (i)--(iii) implies amenability. Indeed, let $(u_m)$ be normalized with $\langle Q_Gu_m,u_m\rangle\to0$ and set $v_m:=|u_m|$, so that $\|v_m\|=1$. Since $\bigl||a|-|b|\bigr|\le|a+b|$ for all $a,b\in\C$, the displayed identity applied to $\Delta_0=dI-A_G$ gives
\[
\langle \Delta_0v_m,v_m\rangle=\sum_{\{x,y\}\in\Ec}\bigl|v_m(x)-v_m(y)\bigr|^2
\ \le\ \sum_{\{x,y\}\in\Ec}\bigl|u_m(x)+u_m(y)\bigr|^2=\langle Q_Gu_m,u_m\rangle\longrightarrow0 ,
\]
so $0=\inf\mathrm{spec}(\Delta_0)$, i.e.\ $\sup\mathrm{spec}(A_G)=d$. For a connected $d$-regular graph this is exactly the amenability criterion recalled above, in its spectral form. Note that the weaker conclusion $\|A_G\|=d$, which follows at once from $\inf\mathrm{spec}(A_G)=-d$, would \emph{not} suffice here: the criterion bears on the top of the spectrum of $A_G$, not on its norm.

For the converse, assume $G$ bipartite with color function $c$ and put $f:=(-1)^{c(\cdot)}$, so that $f(x)+f(y)=0$ on every edge. If $G$ is finite, $f\in\ell^2$ is an eigenvector, $\Delta_0f=2df$, and (i)--(iii) follow. If $G$ is infinite and amenable, $f\notin\ell^2$, but a F\o lner sequence $(\Omega_R)$ furnishes a Weyl sequence. The normalized truncations $f_R:=f\chi_{\Omega_R}/\|f\chi_{\Omega_R}\|$ satisfy
\[
\|(\Delta_0-2d)f_R\|^2\ \le\ 4d^2(1+d)\,\frac{|\partial\Omega_R|}{|\Omega_R|}\ \longrightarrow\ 0.
\]
Indeed, $(\Delta_0-2d)f=0$ pointwise at every $x$ all of whose neighbours lie in $\Omega_R$, so the summand at $x$ vanishes unless $x\in\partial\Omega_R$ or $x$ is adjacent to $\partial\Omega_R$ --- at most $(1+d)|\partial\Omega_R|$ vertices. At each of them $|(\Delta_0-2d)f_R(x)|\le 2d\|f_R\|_\infty$ and $\|f_R\|_\infty=|\Omega_R|^{-1/2}$, whence the bound. Hence $2d\in\mathrm{spec}(\Delta_0)$. It is not an eigenvalue: an eigenfunction $u$ for $2d$ would satisfy $u(y)=-u(x)$ along every edge, hence $|u|$ constant on each connected component, which is incompatible with $u\in\ell^2$ on an infinite connected graph. Finally, for $G$ finite and connected, $\inf\mathrm{spec}(Q_G)=0$ means that some non-zero $u$ satisfies $u(y)=-u(x)$ along every edge; connectedness then forces $|u|$ to be constant and the sign of $u$ to define a bipartition, so (iii) is in that case equivalent to bipartiteness.
\end{proof}

\begin{remark}[Bipartiteness alone is not enough in the infinite case]\label{rem:tree-counterexample}
The amenability hypothesis in Proposition~\ref{prop:sharpness-bipartite} cannot be removed. Let $G=T_d$ be the $d$-regular tree, $d\ge3$: it is infinite, $d$-regular and bipartite, yet by Kesten's theorem~\cite{Kesten} $\mathrm{spec}(A_{T_d})=[-2\sqrt{d-1},2\sqrt{d-1}]$, so
\[
\|\widetilde\Delta_{1}\|\;=\;d+2\sqrt{d-1}\;<\;2d\qquad(d\ge3),
\]
e.g.\ $3+2\sqrt2\approx5.83<6$ for $d=3$. What fails is precisely the Weyl-sequence step: on a non-amenable graph $|\partial\Omega_R|/|\Omega_R|$ is bounded away from $0$ for every exhaustion, and the truncations of $(-1)^{c(\cdot)}$ do not become approximate eigenvectors. Bipartiteness alone therefore guarantees only the symmetry of $\mathrm{spec}(\Delta_0)$ about $d$, not the attainment of $2d$.
\end{remark}

Applying this to the examples: the lattices $\Z^2,\Z^3$, BCC, $\Z^4$ and diamond are bipartite (with the sum-of-coordinates parity, resp.\ the sublattice parity, giving the two color classes) and amenable, whence $\|\widetilde\Delta_{1}\|=2d$; for the triangular, Kagome and FCC lattices the strict inequality $\|\widetilde\Delta_{1}\|<2d$ is read off from the Bloch symbols computed above. We stress that non-bipartiteness alone would not give strictness on an infinite graph.

In $\Z^2$, delete the two edges $\{(0,0),(1,0)\}$ and $\{(5,5),(6,5)\}$ and insert instead $\{(0,0),(5,5)\}$ and $\{(1,0),(6,5)\}$. Every degree is preserved, so the graph is still $4$-regular. Each inserted edge joins two vertices of equal coordinate-sum parity, so it closes an odd cycle and the graph is no longer bipartite. Yet the F\o lner sequence may be taken to consist of large squares translated away from the modified region, on which the truncated function $(-1)^{c(x)}$ built from the bipartition of the unmodified lattice is still almost harmonic; the value $2d$ therefore remains in the spectrum.

\subsection{Application: boundedness and ESA for top-degree Laplacians}\label{subsec:application-line}
Let $G$ be a locally finite graph with maximal degree $\Delta(G)<\infty$. The \emph{line graph} $L(G)$ has vertices given by the edges of $G$, with adjacency determined by sharing an endpoint. An edge $\{u,v\}$ meets exactly $\deg u-1$ other edges at $u$ and $\deg v-1$ at $v$, and these are distinct, so
\begin{equation}\label{eq:line-degree}
\deg_{L(G)}(\{u,v\})=\deg u+\deg v-2,
\qquad
\Delta(L(G))=\max_{\{u,v\}\in\Ec}\bigl(\deg u+\deg v\bigr)-2\;\le\;2(\Delta(G)-1).
\end{equation}
The first identity is what makes the bound below coincide with the edge-degree estimate of~\cite{AM}; see Remark~\ref{rem:anderson-morley}.

\begin{proposition}\label{prop:top-bounded-ESA}
Let $G=(\Vc,\Ec)$ be a locally finite \emph{unweighted} graph (i.e.\ $m_0\equiv m_1\equiv 1$) and let $L(G)$ denote its line graph with maximal degree $\Delta(L(G)):=\sup_{e\in\Ec}\deg_{L(G)}(e)$. Assume $\Delta(L(G))<\infty$. Then
\[
\|\Delta_{1}\|\ \le\ \Delta(L(G))+2 ,
\]
where the ``$+2$'' is the contribution of the constant potential $\mathcal Q(V_{1})$, with $V_{1}\equiv 2$. Since $m_0\equiv m_1\equiv1$ gives $M_1=I$, the energy and degree-normalized blocks coincide here, $\widetilde\Delta_{1}=\Delta_{1}$, and the bound holds for both without further argument. In particular $\Delta_{1}$ extends uniquely to a bounded self-adjoint operator on $\ell^2(\Ec,m_1)$. Combined with~\eqref{eq:line-degree}, the universal $d$-regular bound stated in the introduction follows: for any $d$-regular graph with $d\ge 2$, where $\Delta(G)=d$,
\[
\|\widetilde\Delta_{1}\|\ =\ \|\Delta_{1}\|\ \le\ 2\Delta(G).
\]
\end{proposition}

\begin{proof}
We treat the \emph{unweighted} case $m_0\equiv m_1\equiv 1$; the weighted case is then handled by Proposition~\ref{prop:weighted-extension} below.

\smallskip\noindent
\emph{Step 1. Reduction.}
By the line-complex reduction (Theorem~\ref{thm:kn-esa} applied with $n=1$), there exists a unitary $U:\ell^2(m_1)\to \ell^2(\widehat{\Vc})$ such that $U\,\Delta_{1}\,U^{-1}=\mathcal A_\eta+\mathcal Q(V_{1})$, where $\mathcal A_\eta$ is the signed adjacency operator on the line-complex, with kernel of modulus $w$.
In the present unweighted setting, the diagonal coefficient is, for any oriented edge $e=[x,y]$,
\[
V_{1}(e)\;=\;\frac{1}{m_0(x)}+\frac{1}{m_0(y)}\;=\;2,
\]
so $V_{1}\equiv 2$ is constant.

\smallskip\noindent
\emph{Step 2. Schur bound for $\mathcal A$.}
In the unweighted setting $m_0\equiv m_1\equiv1$, the kernel of $\mathcal A_\eta$ has modulus $w(e,e')=1$ exactly when the edges $e,e'$ share an endpoint in $G$, so that its Schur degree is
\[
D(e)=\sum_{e'\sim e}\bigl|\eta(e,e')\,w(e,e')\bigr|=\deg_{L(G)}(e)\le\Delta(L(G)).
\]
Proposition~\ref{prop:schur}, applied on the line graph $L(G)$, therefore gives $\|\mathcal A_\eta\|\le\Delta(L(G))$.

\smallskip\noindent
\emph{Step 3. Norm estimate.}
Since $\|\mathcal Q(V_{1})\|=2$ (constant diagonal in the unweighted case), $\|U\Delta_{1}U^{-1}\|\le\|\mathcal A_\eta\|+\|\mathcal Q(V_1)\|\le \Delta(L(G))+2$, hence the same holds for $\|\Delta_{1}\|$.

\smallskip\noindent
\emph{Step 4. Relation to $\Delta(G)$.}
Since $\Delta(L(G))\le 2(\Delta(G)-1)$, we obtain
\[
\|\widetilde\Delta_{1}\|\ =\ \|\Delta_{1}\|\ \le\ 2(\Delta(G)-1)+2\ =\ 2\Delta(G),
\]
which is the announced bound; as noted in the statement, $\widetilde\Delta_{1}=\Delta_{1}$ here, so a single estimate covers both. The statement is given for $d$-regular graphs, matching the form announced in the introduction; the step itself uses neither regularity nor $d\ge2$, the inequality $\Delta(L(G))\le2(\Delta(G)-1)$ holding for every graph.
\end{proof}

\begin{remark}[Attribution and comparison with spectral graph theory]\label{rem:anderson-morley}
Proposition~\ref{prop:top-bounded-ESA} is not new for finite graphs. Since $\Delta(L(G))+2=\max_{\{u,v\}\in\Ec}\bigl(\deg u+\deg v\bigr)$ by~\eqref{eq:line-degree}, it is exactly the classical edge-degree bound of Anderson and Morley~\cite{AM}. Their proof already proceeds through the same two steps used above. First, the non-zero spectra of $\mathsf B^{\mathsf T}\mathsf B$ and $\mathsf B\mathsf B^{\mathsf T}$ coincide (our Lemma~\ref{lem:iso-01}). Second, the absolute row sums of $\mathsf B\mathsf B^{\mathsf T}$ --- the line graph adjacency with $2$'s on the diagonal --- are bounded by $\deg u+\deg v$ (our Step~2). What Theorem~\ref{thm:kn-esa} adds is the extension of this mechanism to arbitrary top degree $n$ on a weighted flag complex in the infinite, $\ell^2$ setting, together with the explicit determination of the orientation cocycle; the case $n=1$ is included for calibration and for the sharpness discussion, not as a new estimate.

Sharper constants are moreover available in the finite case. Rojo, Soto and Rojo~\cite{RSR} and Das~\cite{Das} refined~\cite{AM} to
\[
\lambda_{\max}(\Delta_0)\ \le\ \max_{\{u,v\}\in\Ec}\bigl|\Nc_\Gc(u)\cup\Nc_\Gc(v)\bigr|
\ =\ \max_{\{u,v\}\in\Ec}\bigl(\deg u+\deg v-\#(\Nc_\Gc(u)\cap\Nc_\Gc(v))\bigr),
\]
which is never larger than $\Delta(L(G))+2$ and is strictly smaller as soon as every edge lies in a triangle. On the finite torus quotients of the lattices of Table~\ref{tab:lattices-combined} this already yields $10$ on the triangular lattice (two common neighbours per edge), $20$ on FCC (four) and $7$ on Kagome (one), against the values $2d$ equal to $12$, $24$ and $8$ respectively --- and it is exact on the bipartite lattices, which are triangle-free. Our contribution at $n=1$ is therefore the \emph{exact} values of Table~\ref{tab:lattices-combined} and the sharpness criterion of Proposition~\ref{prop:sharpness-bipartite}, not the bound $2d$ itself.

In higher degree the picture is analogous: with respect to the finite-complex bounds of~\cite{FWW,SWF} recalled in Section~\ref{sec:related}, what Theorem~\ref{thm:kn-esa} contributes is the infinite, weighted, $\ell^2$ setting, where the issue is boundedness and essential self-adjointness rather than the location of a largest eigenvalue.
\end{remark}

\begin{table}[ht]
\centering
\caption{Lattice norms: effective bound $\Delta(L(G))+2=2d$ (the value established by the proof of Proposition~\ref{prop:top-bounded-ESA}) and the exact Floquet--Bloch constant. All the lattices listed are amenable, so the bound $2d$ is attained on the bipartite ones ($\Z^2$, $\Z^3$, BCC, $\Z^4$, diamond) by Proposition~\ref{prop:sharpness-bipartite}; on the non-bipartite ones (triangular, FCC, Kagome) the strict inequality is established by the Floquet--Bloch computation of Section~\ref{app:bloch}, not by non-bipartiteness alone.}\label{tab:lattices-combined}
\small
\setlength{\tabcolsep}{4pt}
\begin{tabular}{lcccc}
\hline
\textbf{Lattice} & $\boldsymbol{d}$ & $\boldsymbol{\Delta(L)}$ & $\boldsymbol{2d}$ & $\boldsymbol{\|\widetilde\Delta_{1}\|}$ \\
 & & \textbf{(line)} & \textbf{(bound)} & \textbf{(exact)} \\
\hline
Square $\Z^2$ (bipartite)         & $4$  & $6$  & $8$  & $8$ \\
Triangular (non-bipartite)        & $6$  & $10$ & $12$ & $9$ \\
Cubic $\Z^3$ (bipartite)          & $6$  & $10$ & $12$ & $12$ \\
BCC (bipartite)                   & $8$  & $14$ & $16$ & $16$ \\
FCC (non-bipartite)               & $12$ & $22$ & $24$ & $16$ \\
Kagome (planar, non-bipartite)    & $4$  & $6$  & $8$  & $6$ \\
Diamond (3D, bipartite)           & $4$  & $6$  & $8$  & $8$ \\
Hypercubic $\Z^4$ (bipartite)     & $8$  & $14$ & $16$ & $16$ \\
\hline
\end{tabular}

\smallskip
\noindent {\footnotesize Exact values for the bipartite lattices follow from Proposition~\ref{prop:sharpness-bipartite}: all of them are amenable, so the Weyl sequence built from $(-1)^{c(x)}$ along a F\o lner exhaustion gives $2d\in\mathrm{spec}(\Delta_0)$. Bipartiteness alone would not suffice in the infinite case; see Remark~\ref{rem:tree-counterexample}. For the non-bipartite lattices the exact value comes from Floquet--Bloch analysis: maximization of the closed-form Bloch symbol in Section~\ref{app:bloch}, where the maximum is located in closed form and certified by an algebraic identity. Column $\Delta(L)=2(d-1)$ is the maximum degree of the line graph.}
\end{table}

\begin{proposition}\label{prop:weighted-extension}
Let $G=(\Vc,\Ec)$ be a locally finite weighted graph with vertex weights $m_0$ and edge weights $m_1$, and assume $\Delta(L(G))<\infty$. Assume uniform comparability:
\(
0<c_0\le m_0(x)\le C_0\) for all $x\in\Vc$ and \(0<c_1\le m_1(x,y)\le C_1\) for all $x\sim y$. Only $c_0$ and $C_1$ enter the bound~\eqref{eq:nat-bound}; the two-sided form is stated because it is the standard uniform comparability assumption, and because it is the hypothesis under which the quantitative similarity of Lemma~\ref{lem:elliptic-similarity} applies.
Then the top-degree Hodge Laplacian $\Delta_{1}$ on the core $\Cc^1_{c}$ satisfies
\begin{equation}\label{eq:nat-bound}
\|\Delta_{1}\|\;\le\;\frac{C_1}{c_0}\bigl(\Delta(L(G))+2\bigr),
\end{equation}
and extends uniquely to a bounded self-adjoint operator on $\ell^2(\Ec,m_1)$; the same bound holds verbatim for $\widetilde\Delta_{1}$ (Lemma~\ref{lem:iso-01}).
\end{proposition}

\begin{proof}
The line-complex formula of Theorem~\ref{thm:kn-esa} (specialized to $n=1$) gives the decomposition
\[
U\,\Delta_{1}\,U^{-1}\;=\;\mathcal A_\eta\;+\;\mathcal Q(V_{1}),
\]
where $\mathcal A_\eta$ is a (signed) symmetric adjacency operator with kernel of absolute value
\[
w(e,e')\;=\;\frac{\sqrt{m_1(e)\,m_1(e')}}{m_0(e\cap e')}\,\mathbf 1_{\{e\sim e'\}},
\]
and $V_{1}(e)=\sum_{v\subset e}\frac{m_1(e)}{m_0(v)}$ is the diagonal potential. Steps~1 and~2 below verify hypotheses (A1) and (A2) of Theorem~\ref{thm:kn-esa} under the comparability assumption, and yield the norm bound
\begin{equation}\label{eq:weighted-bound-decomp}
\|\Delta_{1}\|\;\le\;\|D\|_\infty\;+\;\|V_{1}\|_\infty.
\end{equation}

\smallskip
\emph{Step 1. Bound on the line-complex Schur degree $D$.}
For every edge $e\in\Ec$,
\[
D(e)\;\le\;\frac{C_1}{c_0}\,\deg_{L(G)}(e),\qquad
\text{hence}\qquad
\|D\|_\infty\;\le\;\frac{C_1}{c_0}\,\Delta(L(G)).
\]

\smallskip
\emph{Step 2. Bound on the diagonal potential $V_{1}$.}
For every edge $e\in\Ec$ with endpoints $x,y$,
\[
V_{1}(e)\;=\;\frac{m_1(e)}{m_0(x)}+\frac{m_1(e)}{m_0(y)}\;\le\;2\,\frac{C_1}{c_0}.
\]

\smallskip
\emph{Step 3. Combination.}
Combining Steps~1 and~2 with~\eqref{eq:weighted-bound-decomp} yields~\eqref{eq:nat-bound}:
\begin{equation}\label{eq:combined-prelim}
\|\Delta_{1}\|\;\le\;\frac{C_1}{c_0}\,\Delta(L(G))+2\,\frac{C_1}{c_0}\;=\;\frac{C_1}{c_0}\bigl(\Delta(L(G))+2\bigr).
\end{equation}

\end{proof}

\textbf{Acknowledgements.} The authors thank the anonymous referees for their careful reading and constructive comments, which improved the paper. We are also grateful to P.~Bartmann and M.~Keller for sharing their preprint and for the ensuing correspondence.

\section*{Declarations}

\subsection*{Funding}
This research received no specific grant from any funding agency in the public, commercial, or not-for-profit sectors.

\subsection*{Ethical approval} Not applicable.

\subsection*{Informed consent} Not applicable.

\subsection*{Author Contributions}
Marwa Ennaceur and Amel Jadlaoui read and approved the final manuscript.

\subsection*{Data Availability Statement}
No datasets were generated or analyzed during the current study. All the constants reported in Section~\ref{app:bloch} are obtained in closed form from the Bloch symbols given there, and no numerical computation is relied upon.

\subsection*{Conflict of Interest}
The authors declare that they have no conflict of interest.

\subsection*{Clinical trial number} Not applicable.

\textbf{Corresponding author:} Amel Jadlaoui (email: amel.jadlaoui@yahoo.com)
\end{document}